\documentclass[a4paper,11pt]{article}

\usepackage{ucs}
\usepackage[utf8x]{inputenc}
\usepackage{graphicx}
\usepackage{amsfonts}
\usepackage{amssymb}
\usepackage{amsmath}
\usepackage{amsthm}
\usepackage{dsfont}
\usepackage{enumerate}
\usepackage{fullpage}
\usepackage{ifthen}
\usepackage{epic}
\usepackage{authblk}
\usepackage{textcomp}
\usepackage[all]{xy}
\usepackage{mathrsfs}
\usepackage{subfigure}
\usepackage{titlesec}
\usepackage[small]{caption}
\usepackage{tikz,tkz-tab}
\usepackage{xcolor}
\usepackage{comment}

\usetikzlibrary{arrows,shapes,positioning}
\usetikzlibrary{decorations.markings}
\tikzstyle arrowstyle=[scale=1.5]
\tikzstyle directed=[postaction={decorate,decoration={markings,
    mark=at position .5 with {\arrow[arrowstyle]{stealth}}}}]
\tikzstyle reverse directed=[postaction={decorate,decoration={markings,
    mark=at position .5 with {\arrowreversed[arrowstyle]{stealth};}}}]

\usepackage[hypertexnames=false,colorlinks=true,linkcolor=blue,citecolor=blue]{hyperref}
\usepackage[numbers,comma,square,sort&compress]{natbib}
\usepackage[letterpaper,text={7in,9in},centering]{geometry}

\usepackage{titlesec}
\graphicspath{{eps/}{pdf/}}

\graphicspath{{Figures_Manuscript/}}

\numberwithin{equation}{section}

\newtheorem{lem}{Lemma}[section]

  {\begin{trivlist}\item[]{\bf Hypothesis #1 }\em}{\end{trivlist}}

\newcommand{\R}{\mathbb{R}}

\newcommand{\N}{\mathbb{N}}
\newcommand{\Z}{\mathbb{Z}}

\newcommand{\W}{\mathcal{W}}

\renewcommand{\Re}{\mathrm{Re}}

\newcommand{\mbi}{\mathbf{i}}

\newcommand{\bqq}{\begin{equation}}
\newcommand{\eqq}{\end{equation}}
\newcommand{\bqs}{\begin{equation*}}
\newcommand{\eqs}{\end{equation*}}

\newenvironment{Proof}[1][.]%
 {\begin{trivlist}\item[]\textbf{Proof#1 }}%
 {\hspace*{\fill}$\rule{0.3\baselineskip}{0.35\baselineskip}$\end{trivlist}}

\title{Wave propagation phenomena in nonlinear hierarchical neural networks with predictive coding feedback dynamics}
\author[1]{\textsc{Andrea Alamia}}
\author[2]{\textsc{L\'ea Dalli\`es}}
\author[2]{\textsc{Gr\'egory Faye}\footnote{ \texttt{gregory.faye@math.univ-toulouse.fr}}}
\author[1,3]{\textsc{Rufin VanRullen}}
\affil[1]{CerCo, CNRS, 31052 Toulouse, France}
\affil[2]{Institut de Math\'ematiques de Toulouse ; UMR5219, Universit\'e de Toulouse ; UPS IMT, F-31062 Toulouse Cedex 9 France}
\affil[3]{ANITI, Universit\'e de Toulouse, 31062, Toulouse, France}
\begin{document}
\maketitle

\begin{abstract}
We propose a mathematical framework to systematically explore the propagation  properties of a class of continuous in time nonlinear neural network models comprising a hierarchy of processing areas, mutually connected according to the principles of predictive coding. We precisely determine the conditions under which upward propagation, downward propagation or even propagation failure can occur in both bi-infinite and semi-infinite idealizations of the model. We also study the long-time behavior of the system when either a fixed external input is constantly presented at the first layer of the network or when this external input consists in the presentation of constant input with large amplitude for a fixed time window followed by a reset to a down state of the network for all later times. In both cases, we numerically demonstrate the existence of threshold behavior for the amplitude of the external input characterizing whether or not a full propagation within the network can occur. Our theoretical results are consistent with predictive coding theories and allow us to identify regions of parameters that could be associated with dysfunctional perceptions.
\end{abstract}

\section{Introduction}

In most sensory systems, cortical areas are organized in a clearly hierarchical structure, as in the visual system~\cite{vanessen1991}. The flow of information within this strongly hierarchical organization is clearly bidirectional. Indeed, it has been shown~\cite{bullier2001} that there are comparable numbers of fibers sending neural signals top-down (i.e., from higher to lower levels of the hierarchy) as they are projecting in the opposite direction. Several theories aim to explain the functional role of cortical connections within hierarchical brain regions processing sensory information. Even though there are relevant interactions between the top-down and bottom-up processes \cite{mcmains2011interactions}, they serve supposedly different functional roles. On the one hand, bottom-up or ``feed-forward'' propagation of information is often interpreted as the integration of sensory input and the extraction of relevant features. On the other hand, the opposite, feedback direction is often considered as a modulatory or contextual signal, which can be involved, for example, in the generation of center-surround interactions \cite{bullier2001role,gilbert2013top}. A comprehensive interpretation of bottom-up and top-down neural flow is provided by the predictive coding (PC) theory \cite{RB99,Huang,Shipp}, which states that each region in the cortical hierarchy generates predictions about what caused its activity. In particular, each region propagates these predictions ''top-down'' to the preceding layer, where a prediction error can be computed as the difference between the prediction and the actual activity. Importantly, this prediction error is then propagated forward to iteratively correct and update the neural representation that will generate the next predictions. Over time (and as long as the sensory input does not change), the system converges into a state where top-down predictions fully explain all bottom-up inputs, and no prediction error is transmitted. 
Even though it is still debated whether PC can provide a unifying framework to describe neural activities and, more broadly, cognition \cite{Aitchison,Bowman}, several theoretical and computational models have been proposed in the light of the PC theory \cite{Shipp,Egner}, which remains very influential in the literature. For example, some experimental studies have recently proposed that oscillatory traveling waves are the neural signatures supporting PC processes in the cortex \cite{alamia,alamiaPlos,schwenk}, revealing the propagation of neural activity as top-down predictions and bottom-up prediction errors. On the other hand, Predictive Coding dynamics proved very fruitful also in machine learning \cite{Chalasani}, increasing the performance of artificial neural networks in different tasks \cite{Lotter,Choksi}, as well as providing evidence on the role of top-down processes in uncertain or adversarial environments \cite{AlamiaFeed}. 
Nevertheless, there is an overall lack of mathematical understanding of emerging dynamic properties that can be generated by (deep) neural networks implementing predictive coding dynamics, and of whether these properties can be informative for biologically relevant models.

In a previous work \cite{FFVR23}, we systematically investigated a linear model where a potentially infinite number of neuronal layers exchanged signals in both directions according to predictive coding principles. We found that the stable propagation of information in such a system can be explored analytically as a function of its initial state, its internal parameters (controlling the strength of inputs, predictions, and error signals) and its connectivity (e.g. convolution kernels). We showed that it is actually possible to determine in which direction, and at what speed neural activity propagates in the system. We also demonstrated that different neural assemblies, which reflect distinct eigenvectors of the underlying connectivity matrices, can simultaneously and independently display different properties in terms of stability, propagation speed, or direction. In some example cases (ring model of orientations), we saw that only a few assemblies support signal propagation, implying that the system acts as a filter on its inputs, and these assemblies propagate
information at different speeds, implementing a coarse-to-fine analysis. Finally, we  derived a time-continuous version of the model which typically allowed us to investigate the inclusion of communication delays
between layers. In this setting, we demonstrated the emergence of genuine oscillatory dynamics and traveling waves in various frequency bands consistent with those observed in the brain. 

However, the model we studied was operating strictly in a linear regime, which is a strong limitation from a biological perspective. Indeed, it is well-known that nonlinearities play a key role in both biological neurons \cite{NonLinearity,MDKF} and artificial network architectures \cite{NonLinearityAI}. In this study, we pursue our mathematical investigations by analytically studying a fully nonlinear model that is continuous in time. In a nutshell, our present analysis generalizes the results found in the linear regime with a systematic characterization of propagation properties within the network as a function of the internal parameters and the shape of the nonlinearity introduced into the model. We precisely determine the conditions under which upward propagation, downward propagation or even propagation failure can occur in both bi-infinite and semi-infinite idealizations of the model. We also study the long-time behavior of the system when a fixed external input is constantly presented at the first layer of the network, and demonstrate the existence of a threshold for the amplitude of the input characterizing whether or not a full propagation within the network can occur. Finally, we also explore the propagation properties and threshold phenomena for so-called flashed external input, an external input consisting of the presentation of constant input with large amplitude for a fixed time window, followed by a reset to a down state of the network for all later times. Our results are largely consistent with PC theories and allow us to identify regions of parameters that could be associated with dysfunctional perceptions.

\section{Derivation of the continuous in time model -- Voltage-like approach}

Motivated by our previous work \cite{FFVR23}, we derive a continuous in time nonlinear model of hierarchical neural networks with predictive coding feedback dynamics, which is inspired by the generic formulation of predictive coding proposed in the context of deep learning models.

Already anticipating the final step of the derivation of the model which will consist in taking a formal continuous in time limit, we introduce $\Delta t>0$ which will represent a time step in our formulation.  We consider the following recurrence equation where $\mathscr{V}_j^n\in\R^d$ represents the network neuronal activity at step $n$ and layer $j$
\bqq
\label{update}
\mathscr{V}_j^{n+1}=(1-\beta \Delta t)\mathscr{V}_j^{n}+ \beta \Delta t  \W^f \mathcal{S}( \mathscr{V}_{j-1}^{n+1})-\alpha \Delta t \mathcal{F}_{j-1}^n -\lambda \Delta t \mathcal{B}_{j}^n, \quad j=1,\cdots,J \,,
\eqq
where $\W^f\in\mathscr{M}_d(\R)$ is a $d\times d$ square matrix representing the weights of feedforward connections which we assume to be the same for each layer such that $\W^f \mathcal{S}( \mathscr{V}_{j-1}^{n+1})$ models a nonlinear instantaneous feedforward drive from layer $j-1$ to layer $j$, controlled by hyper-parameter $\beta \geq 0$. Here the nonlinearity $\mathcal{S}:\R^d\to\R^d$ is defined componentwise as follows:
\bqs
\mathcal{S}(X)=(S (x_1),\cdots,S (x_d))^{\bf t}, \quad X=(x_1,\cdots,x_d)^{\bf t},
\eqs
for some real-valued function $S:\R\to\R$ which shall be fixed later.  The term $\mathcal{F}_{j-1}^n$ encodes a feedforward error correction process, controlled by hyper-parameter $\alpha \geq0$, where the reconstruction error $\mathcal{R}_{j-1}^n$ at layer $j-1$, defined as the square error between the representation $\mathscr{V}_{j-1}^n$ and the predicted neural activity $\W^b\mathcal{S}(\mathscr{V}_j^n)$, that is
\bqs
\mathcal{R}_{j-1}^n:= \frac{1}{2} \|\mathscr{V}_{j-1}^n-\W^b \mathcal{S}(\mathscr{V}_j^n)\|^2,
\eqs
propagates to the layer $j$ to update its activity. Here, $\W^b\in\mathscr{M}_d(\R)$ is a $d\times d$ square matrix representing the weights of feedback connections which we assume to be the same for each layer. The contribution $\mathcal{F}_{j-1}^n$ is then taken to be the gradient of $\mathcal{R}_{j-1}^n$ with respect to $\mathscr{V}_j^n$, that is
\bqs
\mathcal{F}_{j-1}^n = \nabla \mathcal{R}_{j-1}^n =  -D\mathcal{S}(\mathscr{V}_j^n) (\W^b) ^\mathbf{t} \left(\mathscr{V}_{j-1}^n - \W^b \mathcal{S}(\mathscr{V}_j^n)\right),
\eqs
where $D\mathcal{S}:\R^d\to\mathscr{M}_d(\R)$ is defined as
\bqs
D\mathcal{S}(X)=\mathrm{diag}\left(S' (x_1),\cdots,S' (x_d)\right), \quad X=(x_1,\cdots,x_d)^{\bf t}.
\eqs
On the other hand, $\mathcal{B}_{j}^n$ incorporates a top-down prediction to update the neural activity at layer $j$. This term thus reflects a feedback error correction process, controlled by hyper-parameter $\lambda \geq 0$. Similar to the feedforward process, $\mathcal{B}_{j}^n$ is defined as the the gradient of $\mathcal{R}_{j}^n$ with respect to $\mathscr{V}_j^n$, that is
\bqs
\mathcal{B}_{j}^n = \nabla \mathcal{R}_{j}^n =  -\W^b \mathcal{S}(\mathscr{V}_{j+1}^n)+\mathscr{V}_j^n.
\eqs
Thus, the update rule \eqref{update} can be equivalently written as
\bqs
\begin{split}
\mathscr{V}_j^{n+1}-\mathscr{V}_j^{n}&= \beta \Delta t  \left( \W^f \mathcal{S}( \mathscr{V}_{j-1}^{n+1})-\mathscr{V}_j^n\right)+\alpha \Delta t  D\mathcal{S}(\mathscr{V}_j^n) (\W^b) ^\mathbf{t} \left(\mathscr{V}_{j-1}^n - \W^b \mathcal{S}(\mathscr{V}_j^n)\right)\\
&~~~+ \lambda \Delta t \left(\W^b \mathcal{S}(\mathscr{V}_{j+1}^n)-\mathscr{V}_j^n\right)
\end{split}
\eqs
for each $j=1,\cdots,J$ and $n\geq0$. At the incoming layer $j=0$, we impose that 
\bqs
\mathscr{V}_0^n=\mathscr{S}_0^n, \quad n \geq0,
\eqs
where $\mathscr{S}_0^n\in\R^d$ is a given source term, which can be understood as the network's constant visual input. At the final layer $j=J$ of the network, there is no possibility of incoming top-down signal, and thus one naturally imposes that
\bqs
\mathcal{B}_{J}^n=0,
\eqs
such that the dynamic at layer $j=J$ is simply given by
\bqs
\mathscr{V}_J^{n+1}-\mathscr{V}_J^{n}= \beta \Delta t  \left( \W^f \mathcal{S}( \mathscr{V}_{J-1}^{n+1})-\mathscr{V}_J^n\right)+\alpha \Delta t  D\mathcal{S}(\mathscr{V}_J^n) (\W^b) ^\mathbf{t} \left(\mathscr{V}_{J-1}^n - \W^b \mathcal{S}(\mathscr{V}_J^n)\right), \quad n\geq0.
\eqs
Finally, at the initial step $n=0$, the neural network is initialized according to
\bqs
\mathscr{V}_j^0=\mathscr{H}_j, \quad j=1,\cdots,J,
\eqs
for some given initial sequence $\left(\mathscr{H}_j\right)_{j=1,\cdots,J}$.

Proceeding similarly as in \cite{FFVR23}, we now interpret $\mathscr{V}_j^n$ and $\mathscr{S}_0^n$ as the approximations of some smooth functions of time $\mathscr{V}_j(t)$ and $\mathscr{S}_0(t)$ respectively evaluated at $t_n:=n\Delta t$, that is $\mathscr{V}_j^n \sim \mathscr{V}_j(t_n)$ and $\mathscr{S}_0^n\sim\mathscr{S}_0(t_n)$. Under this representation, the above system of equations writes
\bqs
\begin{split}
\frac{\mathscr{V}_j(t_{n+1})-\mathscr{V}_j(t_{n})}{\Delta t}&= \beta  \left( \W^f \mathcal{S}( \mathscr{V}_{j-1}(t_{n+1}))-\mathscr{V}_j(t_{n})\right)+\alpha  D\mathcal{S}(\mathscr{V}_j(t_{n})) (\W^b) ^\mathbf{t} \left(\mathscr{V}_{j-1}(t_{n}) - \W^b \mathcal{S}(\mathscr{V}_j(t_{n}))\right)\\
&~~~+ \lambda \left(\W^b \mathcal{S}(\mathscr{V}_{j+1}(t_{n}))-\mathscr{V}_j(t_{n})\right),
\end{split}
\eqs 
such that by taking the formal limit $\Delta t\rightarrow0$, we obtain the following functional differential equation 
 \bqs
\begin{split}
\mathscr{V}_j'(t)&=\beta \left(\W^f \mathcal{S}(\mathscr{V}_{j-1}(t))-\mathscr{V}_j(t)\right) +\alpha D\mathcal{S}(\mathscr{V}_j(t)) (\W^b) ^\mathbf{t}\ \left( \mathscr{V}_{j-1}(t)-\W^b\mathcal{S}( \mathscr{V}_j(t))\right) \\
&~~~+ \lambda  \left(\W^b \mathcal{S}(\mathscr{V}_{j+1}(t))-\mathscr{V}_j(t)\right),
\end{split}
\eqs
for $t>0$ and $j=1,\cdots,J$. The previous equation describes the temporal evolution of each $\mathscr{V}_j(t)$ now interpreted as smooth function of time. At layer $j=0$, we simply obtain 
\bqs
\mathscr{V}_0(t)=\mathscr{S}_0(t), \quad t>0,
\eqs
while at $j=J$, we naturally get
\bqs
\mathscr{V}_J'(t)=\beta \left(\W^f \mathcal{S}(\mathscr{V}_{J-1}(t))-\mathscr{V}_J(t)\right) +\alpha D\mathcal{S}(\mathscr{V}_J(t)) (\W^b) ^\mathbf{t}\ \left( \mathscr{V}_{J-1}(t)-\W^b\mathcal{S}( \mathscr{V}_J(t))\right), \quad t>0.
\eqs
As a consequence, the full model reads
\bqq\label{model}
\left\{
\begin{split}
\mathscr{V}_1'(t)&=\beta \left(\W^f \mathcal{S}(\mathscr{S}_0(t))-\mathscr{V}_1(t)\right) +\alpha D\mathcal{S}(\mathscr{V}_1(t)) (\W^b) ^\mathbf{t}\ \left( \mathscr{S}_{0}(t)-\W^b\mathcal{S}( \mathscr{V}_1(t))\right) \\
&~~~+ \lambda  \left(\W^b \mathcal{S}(\mathscr{V}_{2}(t))-\mathscr{V}_1(t)\right),\\
\mathscr{V}_j'(t)&=\beta \left(\W^f \mathcal{S}(\mathscr{V}_{j-1}(t))-\mathscr{V}_j(t)\right) +\alpha D\mathcal{S}(\mathscr{V}_j(t)) (\W^b) ^\mathbf{t}\ \left( \mathscr{V}_{j-1}(t)-\W^b\mathcal{S}( \mathscr{V}_j(t))\right) \\
&~~~+ \lambda  \left(\W^b \mathcal{S}(\mathscr{V}_{j+1}(t))-\mathscr{V}_j(t)\right), \quad j=2,\cdots,J-1,\\
\mathscr{V}_J'(t)&=\beta \left(\W^f \mathcal{S}(\mathscr{V}_{J-1}(t))-\mathscr{V}_J(t)\right) +\alpha D\mathcal{S}(\mathscr{V}_J(t)) (\W^b) ^\mathbf{t}\ \left( \mathscr{V}_{J-1}(t)-\W^b\mathcal{S}( \mathscr{V}_J(t))\right),\\
\mathscr{V}_j(0)&=\mathscr{H}_j, \quad j=1,\cdots,J,
\end{split}
\right. \quad t>0.
\eqq

In what follows, we shall work with the following sigmoidal nonlinearity
\bqs
S(x):=\frac{1}{1+e^{-\mu (x-\theta)}},
\eqs
which depends on two parameters $\mu>0$ and $\theta\geq 0$. Here, $\mu$ describes the slope of the nonlinearity while $\theta$ represents a given threshold. Such a nonlinearity is commonly used in neural-mass models as a canonical choice for smooth saturating nonlinear activation function \cite{MDKF,nunez,WC72}.

Throughout, we shall always work under the natural assumption that
\bqs
\alpha+\beta+\lambda>0,
\eqs
indicating that the degenerate situation where all hyper-parameters are set to zero is not considered here since it leads to trivial constant dynamics.

\section{One population model as a mathematical landmark}

As a very first step, we shall work under the key simplifying assumption that each neuron in a layer is only connected to the corresponding neuron in the immediately preceding and following layer, with unit weight in each direction. Mathematically, we thus assume that $\W^f$ and $\W^b$ are both equal to the identity matrix, that is
\bqs
\W^f=\W^b=\mathbf{I}_d.
\eqs
Under such an assumption, the model \eqref{model} reduces to the analysis of the following set of equations
\bqs
v_j'(t)=\beta \left(S(v_{j-1}(t))-v_j(t)\right) +\alpha S'(v_j(t)) \left( v_{j-1}(t)-S( v_j(t))\right)+ \lambda  \left(S(v_{j+1}(t))-v_j(t)\right),
\eqs
for some scalar unknowns $v_j(t)\in\R$ with $j=1,\cdots,J$, such that
\bqs
v_0(t)=s_0(t), \quad \text{ and } \quad v_{J+1}(t)=S^{-1}(v_J(t)), \quad t>0,
\eqs
for a given function $s_0:\R_+\to \R$. The introduction of the extra unknown $v_{J+1}(t)$ is purely artificial, but it allows us to consider the same equation for all layers instead of having a modified equation for layer $J$. Initially, we assume that the network is initialized according to
\bqs
v_j(0)=h_j, \quad j=1,\cdots,J,
\eqs
for some given sequence $(h_j)_{j=1,\cdots,J}$.

In a second step, we shall rescale time through the following change of variable
\bqs
t \leftrightarrow t(\alpha+\beta+\lambda),
\eqs
such that the above equation reads 
\bqs
\begin{split}
v_j'(t)&=\frac{\beta}{\alpha+\beta+\lambda} \left(S(v_{j-1}(t))-v_j(t)\right) +\frac{\alpha}{\alpha+\beta+\lambda} S'(v_j(t)) \left( v_{j-1}(t)-S( v_j(t))\right)\\
&~~~+\frac{\lambda}{\alpha+\beta+\lambda} \left(S(v_{j+1}(t))-v_j(t)\right),
\end{split}
\eqs
for $j=1,\cdots,J$. Thanks to the above change of variable, it is natural to set
\bqs
p:=\frac{\alpha}{\alpha+\beta+\lambda}\in[0,1], \quad \text{ and } \quad q:= \frac{\lambda}{\alpha+\beta+\lambda}\in[0,1],
\eqs
such that
\bqs
\frac{\beta}{\alpha+\beta+\lambda}=1-p-q\in[0,1].
\eqs
The fractions $p$ and $q$ now represent the relative strength of the feedforward error correction process and feedback error correction process respectively, whereas the remaining fraction $1-p-q$ represents the relative strength of the instantaneous feedforward drive.

As a consequence, the model under study is thus
\bqq
\label{cont1d}
v_j'(t)=(1-p-q) \left(S(v_{j-1}(t))-v_j(t)\right) +p S'(v_j(t)) \left( v_{j-1}(t)-S( v_j(t))\right)+ q  \left(S(v_{j+1}(t))-v_j(t)\right),
\eqq
for $t>0$ and $j=1,\cdots,J$, together with
\bqq
\label{BC1d}
v_0(t)=s_0(t), \quad \text{ and } \quad v_{J+1}(t)=S^{-1}(v_J(t)), \quad t>0,
\eqq 
and
\bqq
\label{IC1d}
v_j(0)=h_j, \quad j=1,\cdots,J.
\eqq

\subsection{Wave propagation on an infinite depth network}

As is usual in such a context, it is first useful to consider the model \eqref{cont1d}-\eqref{BC1d}-\eqref{IC1d} in the case where the network is assumed to be bi-infinite (no boundaries). That is we look at \eqref{cont1d} set on the lattice $\Z$, namely 
\bqq\label{cont1dZ}
v_j'(t)=(1-p-q) \left(S(v_{j-1}(t))-v_j(t)\right) +p S'(v_j(t)) \left( v_{j-1}(t)-S( v_j(t))\right)+ q  \left(S(v_{j+1}(t))-v_j(t)\right),
\eqq
for $t>0$ and $j\in\Z$ with initial condition
\bqq
\label{IC1dZ}
v_j(0)=h_j, \quad j\in\Z.
\eqq
We recall that $(p,q)\in[0,1]\times[0,1]$. Our aim is to understand the long-time behavior of the solutions of \eqref{cont1dZ} starting for a class of specific initial conditions. In order to motivate the class of initial conditions that we shall consider, we first look at the intrinsic dynamics of \eqref{cont1dZ}, that is, the dynamics around stationary homogeneous solutions.

\subsubsection{Stationary homogeneous solutions}

Stationary homogeneous solutions of \eqref{cont1dZ} are time independent sequences of the form $v_j(t)=x$ for all $j\in\Z$ for some $x\in\R$ which satisfies the following equation 
\bqs
0 = (-x+S(x))( 1-p-pS'(x)):=F_p(x).
\eqs
Recalling our explicit expression for $S$, solving for $x=S(x)$, we readily obtain an explicit parametrization of the threshold $\theta$ as a function of $x$ and $\mu$ via
\bqs
\theta = x+\frac{1}{\mu} \ln\left( \frac{1-x}{x}\right), \quad x\in(0,1), \quad \mu>0.
\eqs
We will work in the regime where $f:x\mapsto x+\frac{1}{\mu} \ln\left( \frac{1-x}{x}\right)$ has exactly three intervals of monotony on $(0,1)$.  By direct inspection, this is equivalent to imposing that $\mu>4$, which we assume that is holds true throughout the rest of the paper.  

As a consequence we get the existence of $0<x_*(\mu)<x^*(\mu)<1$ such that $f'<0$ on $\left(0,x_*(\mu)\right)\cup(x^*(\mu),1)$ and $f'>0$ on $(x_*(\mu),x*(\mu))$. We can thus uniquely invert $f$ on each interval of monotonicity and get three branches of solutions:
\begin{itemize}
\item a branch $x_d(\theta,\mu)$ for $\theta\in(f(x_*(\mu)),+\infty)$ and $\mu>4$;
\item a branch $x_m(\theta,\mu)$ for $\theta\in(f(x_*(\mu),f(x^*(\mu))$ and $\mu>4$;
\item a branch $x_u(\theta,\mu)$ for $\theta\in[0,f(x^*(\mu)))$ and $\mu>4$.
\end{itemize}
At $\theta =f(x_*(\mu)):=\theta_*(\mu)$, we have $x_d(\theta_*(\mu),\mu)=x_m(\theta_*(\mu),\mu)=x_*(\mu)$, while at $\theta =f(x^*(\mu)):=\theta^*(\mu)$, we have $x_m(\theta^*(\mu),\mu)=x_u(\theta^*(\mu),\mu)=x^*(\mu)$.

\begin{figure}[t!]
\centering
\includegraphics[width=.5\textwidth]{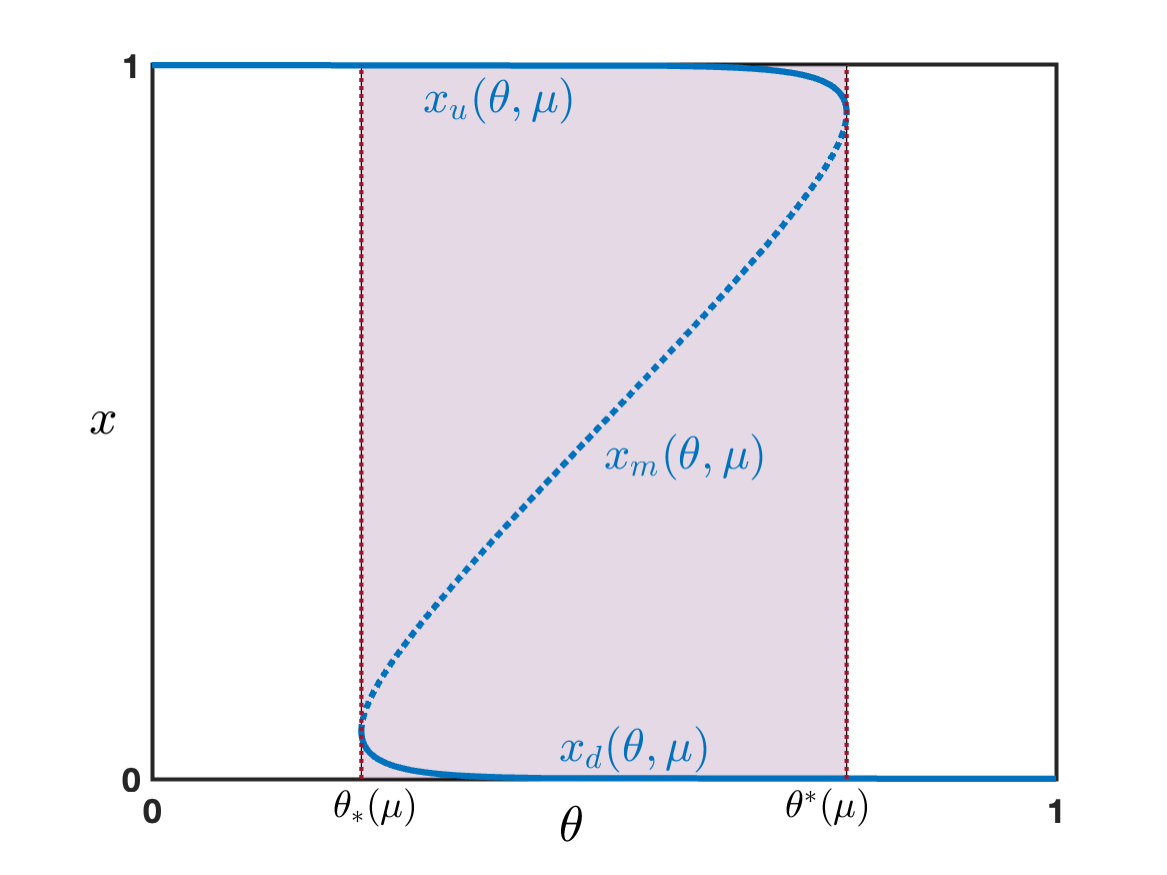}
\caption{Illustration of the three branches of stationary homogeneous solutions $x_{d,m,u}(\theta,\mu)$ solutions of $F_p(x)=0$ as a function of $\theta\in(\theta_*(\mu),\theta^*(\mu))$ when condition \eqref{condpmu} on $\mu$ and $p$ is satisfied. The two branches $x_{d}(\theta,\mu)$ and $x_{u}(\theta,\mu)$ (solid line) are stable while the branch $x_{m}(\theta,\mu)$ (dashed line) is unstable. The so-called bistable regime where the three branches coexist is shaded in violet.}
 \label{fig:statsol}
\end{figure}

On the other hand, solving for $p S'(x)=1-p$ which equivalently reads $S(x)(1-S(x))=\frac{1-p}{p \mu}$, we either obtain two branches of solutions or no solutions depending on the ratio $\frac{1-p}{p \mu}$. More precisely, if $\frac{\mu}{4}<\frac{1-p}{p}$ then there is no such solution, whereas if $\frac{\mu}{4} \geq \frac{1-p}{p}$ there are two branches of solutions. In order to restrict ourselves to a regime where only three stationary homogeneous solutions coexist, we shall work under the assumption that
\bqs
0<\frac{\mu}{4} < \frac{1-p}{p}.
\eqs
Keeping in mind that we have $\mu>4$ and $p\in[0,1]$, the above condition can equivalently be stated as
\bqq
\label{condpmu}
0\leq p < \frac{4}{4+\mu}< \frac{1}{2}, \quad \mu>4.
\eqq
Let us remark that $x_*(\mu)$ and $x^*(\mu)$ can be computed explicitly and are given by
\bqs
x_*(\mu) = \frac{1}{2}-\sqrt{\frac{1}{4}-\frac{1}{\mu}}\in\left(0,\frac{1}{2}\right), \quad x^*(\mu) = \frac{1}{2}+\sqrt{\frac{1}{4}-\frac{1}{\mu}}\in\left(\frac{1}{2},1\right).
\eqs
We refer to Figure~\ref{fig:statsol} for an illustration.

\subsubsection{Linear stability analysis of stationary homogenous solutions}

We start by investigating the linear stability of the stationary homogeneous  solutions. Linearizing around $v_j=x$ with $x\in\left\{x_d(\theta,\mu),x_m(\theta,\mu),x_u(\theta,\mu) \right\}$, the linear equation simply reads
\bqq
u_j'(t)=(1-q) S'(x) u_{j-1}(t)+\left(p S''(x)(x-S(x))-p (S'(x))^2-(1-p) \right)u_j(t)+ q S'(x) u_{j+1}(t),
\label{lincont1d}
\eqq
for each $t>0$ and $j\in\Z$.  Since stationary homogeneous solutions satisfy $x=S(x)$, the above linear equation reduces to 
\bqs
u_j'(t)=(1-q) S'(x) u_{j-1}(t)-\left( p (S'(x))^2+1-p \right)u_j(t)+q S'(x) u_{j+1}(t).
\eqs
To investigate the stability, we look for pure exponential functions of the form $u_j(t)=e^{\nu t + \mbi \varphi j}$ for $\varphi\in[-\pi,\pi]$ for some complex number $\nu$ that satisfies
\bqs
\nu(\varphi)=(1-q) S'(x) e^{-\mbi \varphi}-(1-p)-p (S'(x))^2+q S'(x) e^{\mbi \varphi}, \quad \varphi\in[-\pi,\pi].
\eqs
We note that 
\bqs
\Re(\nu(\varphi))=S'(x) \left(\cos(\varphi)-1\right) +(S'(x)-1)(1-p -pS'(x)), \quad \varphi \in[-\pi,\pi].
\eqs
Since $S'(x)>0$ by definition of the sigmoidal function, we deduce that
\bqs
\underset{\varphi \in[-\pi,\pi]}{\max}\Re(\nu(\varphi))=\Re(\nu(0))=(S'(x)-1)(1-p -pS'(x)).
\eqs
As a consequence, $\underset{\varphi \in[-\pi,\pi]}{\max}\Re(\nu(\varphi))=0$ if and only if $S'(x)=1$ since  $p S'(x_\mu)\neq 1-p$ by assumption \eqref{condpmu} on $p$. Based on the previous analysis on the existence of the branches  $x_{d,m,u}(\theta,\mu)$, we get that $S'(x)=1$ precisely when $\theta=\theta_*(\mu)$ or $\theta=\theta^*(\mu)$. As a consequence, we have
\begin{itemize}
\item for $x=x_d(\theta,\mu)$ then $\underset{\varphi \in[-\pi,\pi]}{\max}\Re(\nu(\varphi))<0$;
\item for $x=x_m(\theta,\mu)$ then $\underset{\varphi \in[-\pi,\pi]}{\max}\Re(\nu(\varphi))>0$;
\item for $x=x_u(\theta,\mu)$ then $\underset{\varphi \in[-\pi,\pi]}{\max}\Re(\nu(\varphi))<0$;
\end{itemize}
 indicating that branches $x_d(\theta,\mu)$ and $x_u(\theta,\mu)$ are linearly stable while the branch $x_m(\theta,\mu)$ is linearly unstable. These stability properties also hold true for the nonlinear dynamics. We shall interpret the stationary solutions along the branch $x_d(\theta,\mu)$ as being \emph{down states}, representing some kind of rest state neural activity of the network, while we interpret the stationary solutions along the branch $x_u(\theta,\mu)$ as being \emph{up states}, representing some kind of high activity of the neural network. Once again, we refer to Figure~\ref{fig:statsol} for an illustration. 

\subsubsection{Bistable traveling waves}

The analysis conducted in the previous section indicates that the regime $\theta\in(\theta_*(\mu),\theta^*(\mu))$, where the three branches of stationary homogeneous solutions $x_{d,m,u}(\theta,\mu)$ coexist, is \emph{bistable} since the two branches $x_d(\theta,\mu)$ and $x_u(\theta,\mu)$ are stable while the branch $x_m(\theta,\mu)$ is unstable. It is thus natural to ask which of the two stable states $x_d(\theta,\mu)$ or $x_u(\theta,\mu)$ is more stable. The notion of stability is of course very important. We shall focus on a very specific class of perturbations. Namely, we aim at considering initial sequences of the form
\bqs
h_j^{u\to d}=\left\{
\begin{array}{cc}
x_u(\theta,\mu), & \text{ if } j\leq0,\\
x_d(\theta,\mu), & \text{ if } j\geq1,
\end{array}
\right. \quad \text{ or } \quad h_j^{d\to u}=\left\{
\begin{array}{cc}
x_d(\theta,\mu), & \text{ if } j\leq0,\\
x_u(\theta,\mu), & \text{ if } j\geq1,
\end{array}
\right.
\eqs
where half of the network is initialized with the up state and the other half with the down state. The question is thus to determine which state will be able to propagate through the network, and if possible, characterize the wave speed of the propagation as a function of the various parameters of the system. We have two parameters for the nonlinearity $\theta$ and $\mu$ which are set as
\bqs
4<\mu, \quad \text{ and } \quad \theta\in(\theta_*(\mu),\theta^*(\mu)),
\eqs
and two hyper-parameters $p$ and $q$ which verify
\bqs
0\leq p < \frac{4}{4+\mu}<\frac{1}{2}, \quad \text{ and } \quad q\in[0,1].
\eqs
We define the parameter set $\mathcal{P}$ as
\bqs
\mathcal{P}:=\left\{\Lambda:= (\theta,\mu,p,q)\in\R_+^4 ~|~ \theta\in(\theta_*(\mu),\theta^*(\mu)), ~ 4<\mu, ~ 0\leq p < \frac{4}{4+\mu},~  q\in[0,1]\right\}.
\eqs

\begin{figure}[t!]
\centering
\subfigure[$c_{u\to d}<0$.]{\includegraphics[width=.32\textwidth]{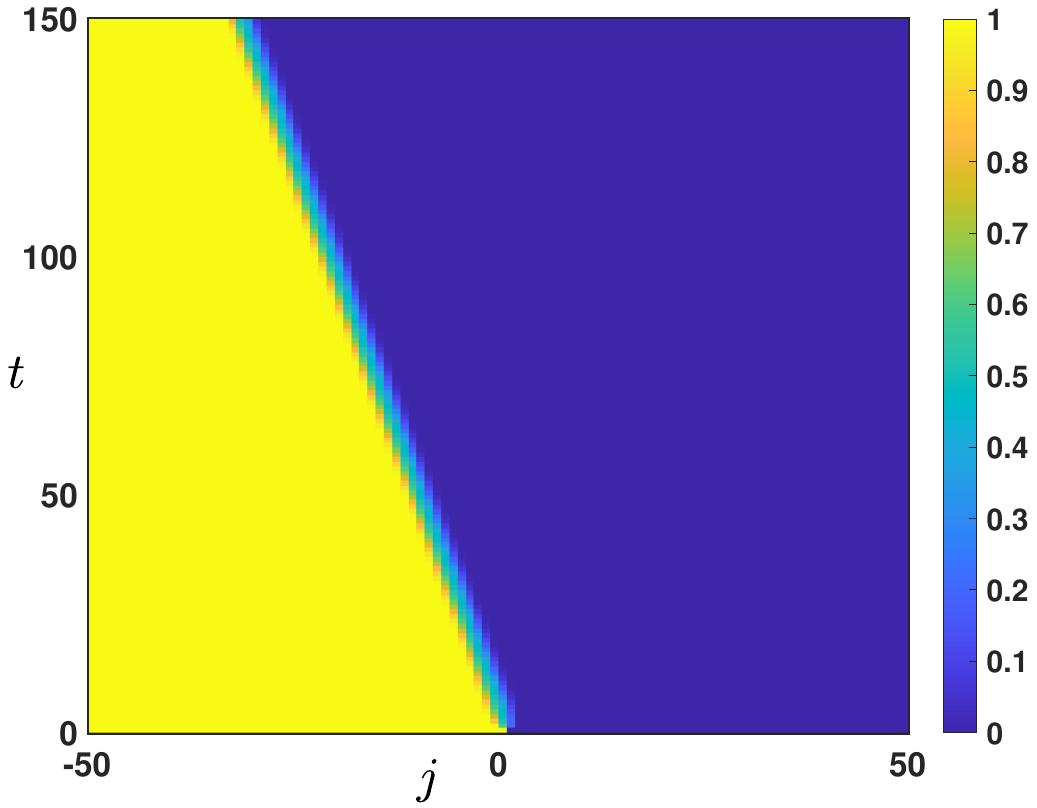}}
\subfigure[Initial condition.]{\includegraphics[width=.3\textwidth]{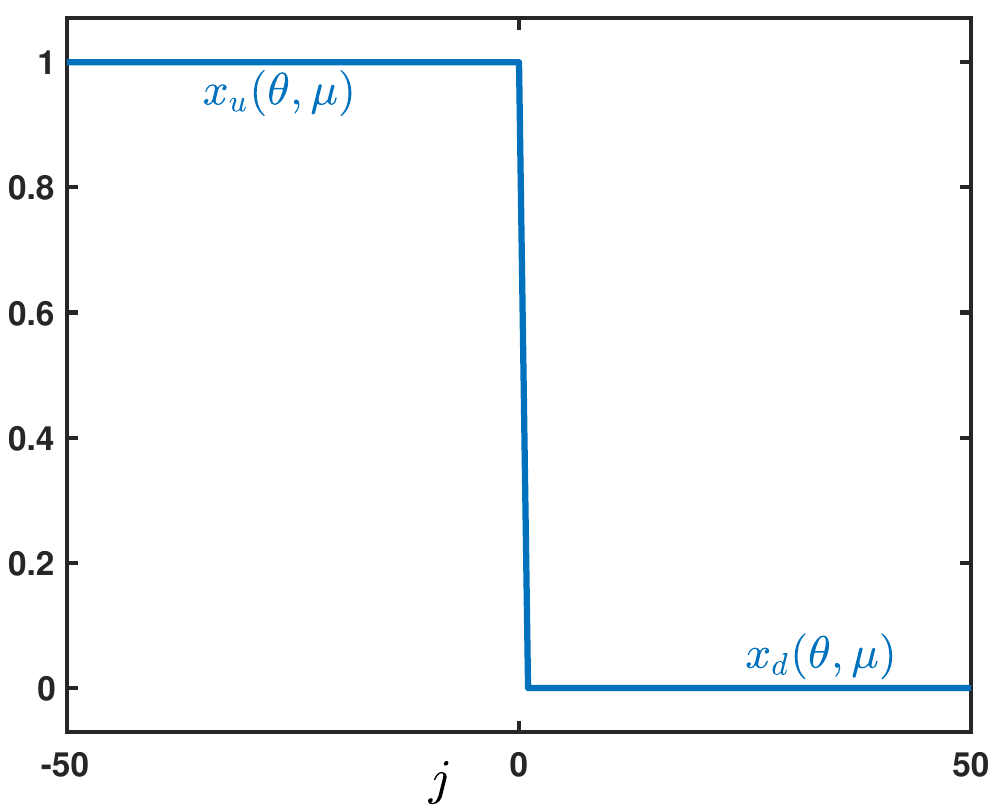}}
\subfigure[$c_{u\to d}>0$.]{\includegraphics[width=.32\textwidth]{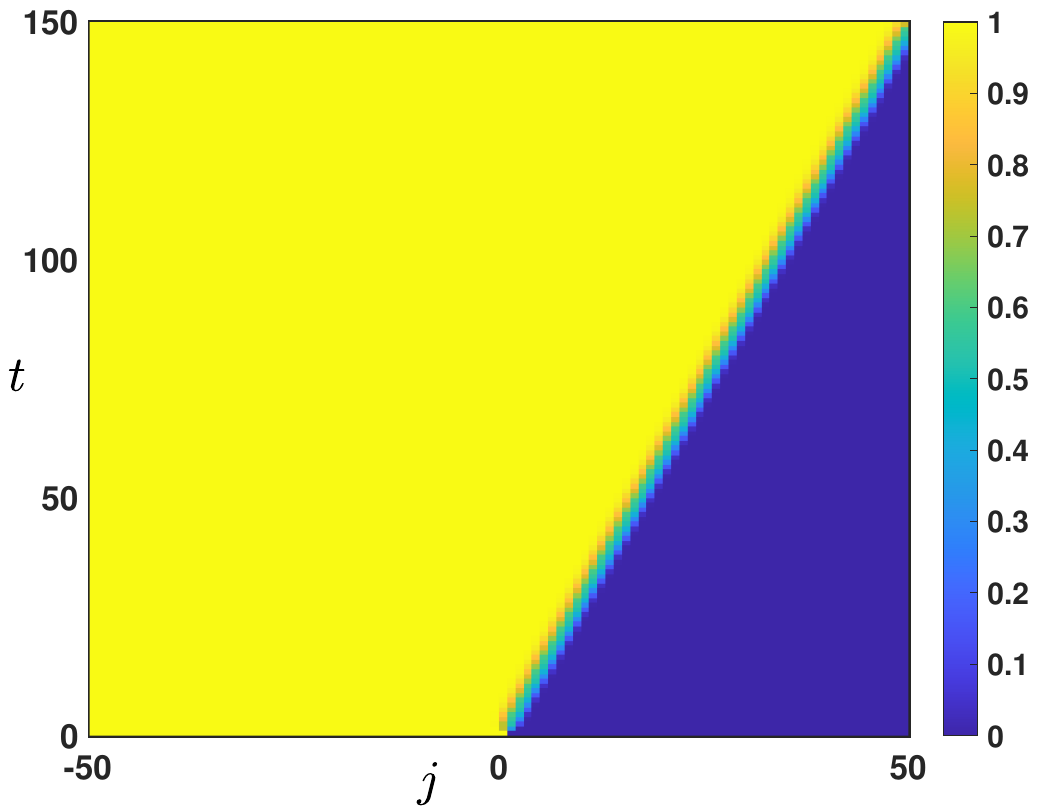}}
\subfigure[$c_{d\to u}<0$.]{\includegraphics[width=.32\textwidth]{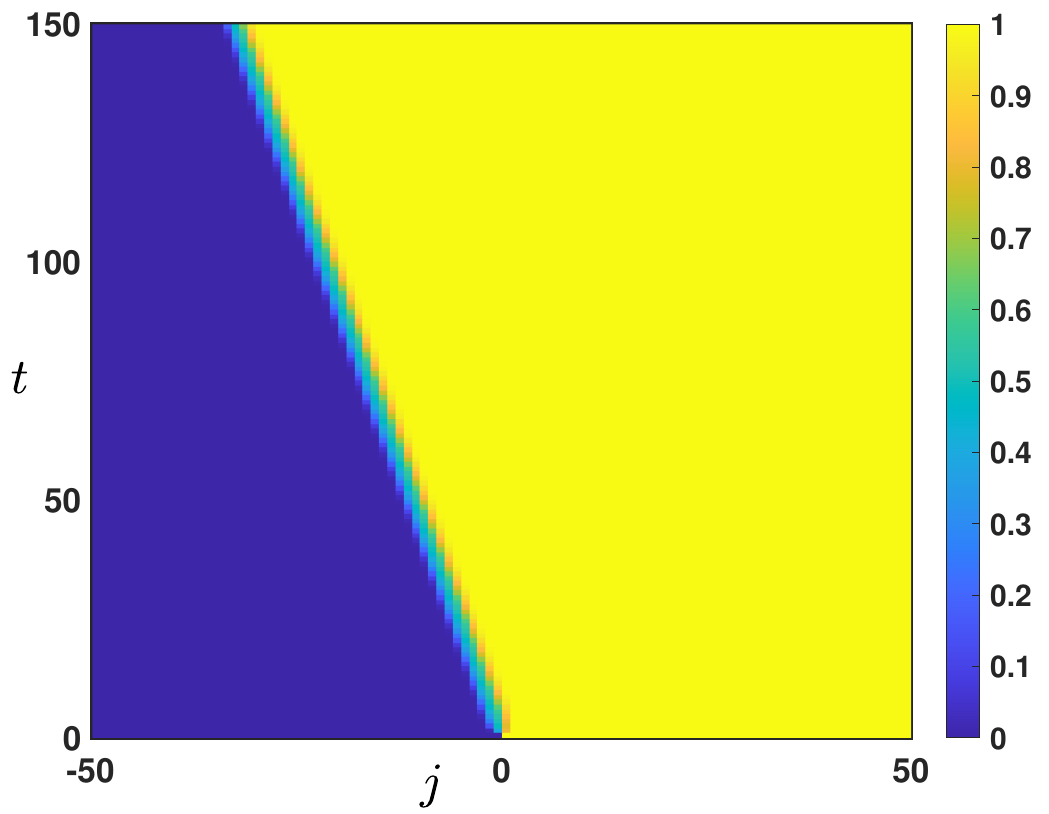}}
\subfigure[Initial condition.]{\includegraphics[width=.3\textwidth]{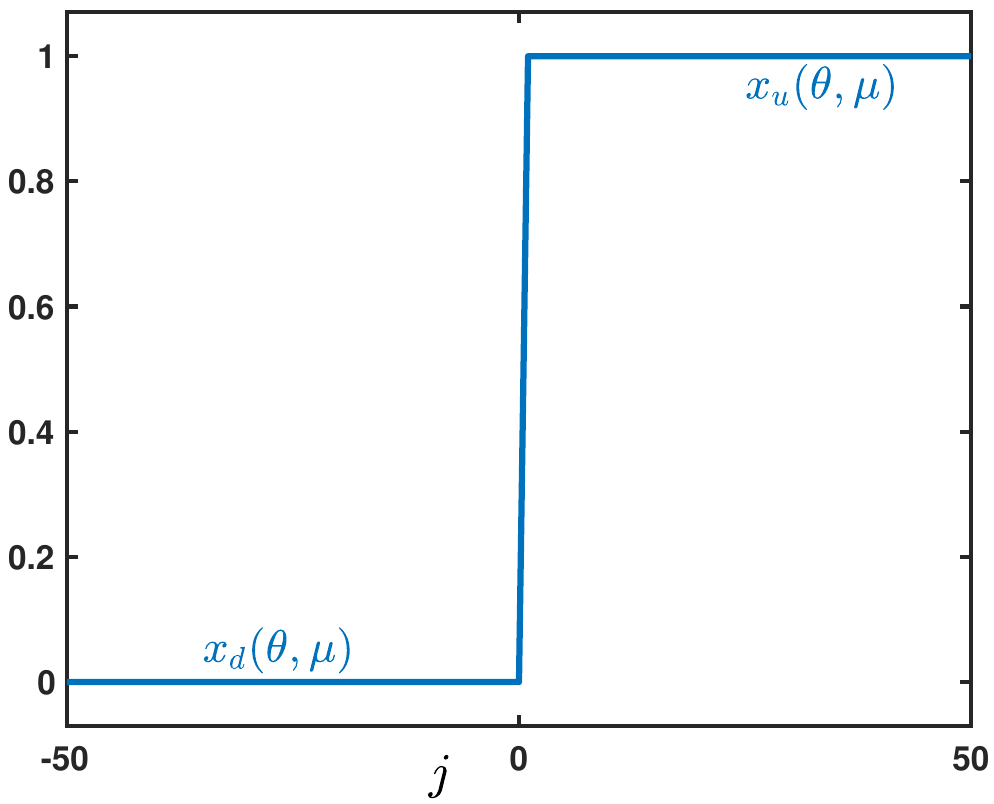}}
\subfigure[$c_{d\to u}>0$.]{\includegraphics[width=.32\textwidth]{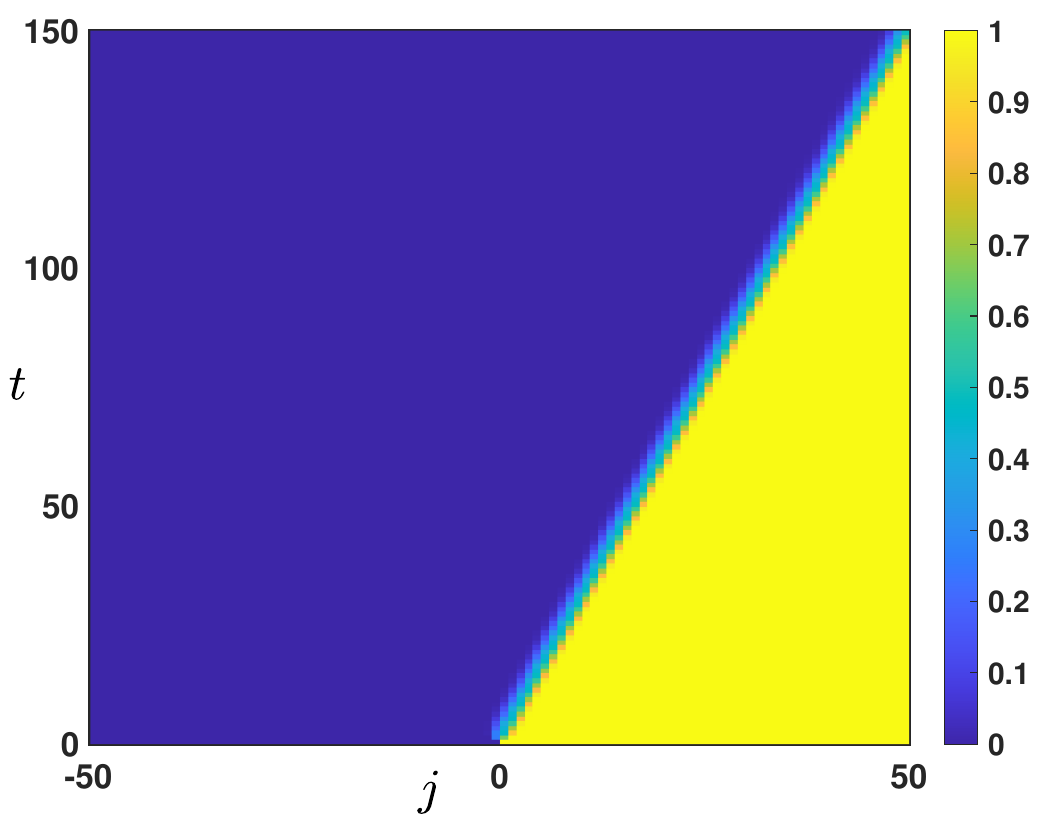}}
\caption{Space-time plots of the solution $(v_j)_{j\in\Z}$ of \eqref{cont1dZ} starting from $(h_j^{u\to d})_{j\in\Z}$ (b) or $(h_j^{d\to u})_{j\in\Z}$ (e) in the case where $c_{u\to d}<0$ (a) and $c_{u\to d}>0$ (c) or $c_{d\to u}<0$ (d) and $c_{d\to u}>0$ (f). Values of the parameters were set to $(\theta,\mu,p)=(0.5,16,0.1)$ with $q=0.6$ in (a)-(d) and $q=0.4$ in (c)-(f).}
 \label{fig:STP}
\end{figure}

We shall first look for special solutions of \eqref{cont1dZ} in the form of traveling wave solutions connecting at infinity the two stable stationary homogeneous solutions. More precisely, we let $v_j(t)=\Phi(j-ct)$ where the profile $\Phi$ and the real $c$ satisfy
\bqq\label{TW}
\left\{
\begin{split}
-c\Phi'(\xi)&=\mathscr{N}(\Phi(\xi-1),\Phi(\xi),\Phi(\xi+1),\Lambda), \quad \forall\xi\in\R,\\
\Phi(-\infty)&=x_u(\theta,\mu), \quad \Phi(+\infty)=x_d(\theta,\mu),
\end{split}
\right.
\eqq
and similarly we let $v_j(t)=\Psi(j-ct)$ where the profile $\Psi$ and the real $c$ satisfy
\bqq\label{TWbis}
\left\{
\begin{split}
-c\Psi'(\xi)&=\mathscr{N}(\Psi(\xi-1),\Psi(\xi),\Psi(\xi+1),\Lambda), \quad \forall\xi\in\R,\\
\Psi(-\infty)&=x_d(\theta,\mu), \quad \Psi(+\infty)=x_u(\theta,\mu),
\end{split}
\right.
\eqq
where we have set
\bqs
\mathscr{N}(u,v,w,\Lambda):=(1-p-q) \left(S(u)-v\right) +p S'(v) \left( u-S(v)\right)+q\left(S(w)-v\right).
\eqs
Here, the map $\mathscr{N}$ is a short-end notation to denote the nonlinear terms appearing in the right-hand side of \eqref{cont1dZ}.  Let us remark that this map $\mathscr{N}: [0,1]^3\times \mathcal{P}\rightarrow\R$ is smooth in all its arguments and that for any $v\in[0,1]$ and $\Lambda\in\mathcal{P}$ one has that $\mathscr{N}(v,v,v,\Lambda)=(-v+S(v))(1-p-pS'(v))=F_p(v)$ is of bistable type since we work in the regime $\theta\in(\theta_*(\mu),\theta^*(\mu))$ with $\mu>4$ and $p<\frac{4}{4+\mu}<\frac{1}{2}$. Furthermore, we also have
\bqs
\forall (u,v,w,\Lambda)\in [0,1]^3\times \mathcal{P}, \quad \partial_u \mathscr{N}(u,v,w,\Lambda)=(1-p-q)S'(u)+pS'(v)>0, \quad \partial_w \mathscr{N}(u,v,w,\Lambda)=qS'(w)>0.
\eqs

\begin{figure}[t!]
\centering
\subfigure[$q=0.35$.]{\includegraphics[width=.325\textwidth]{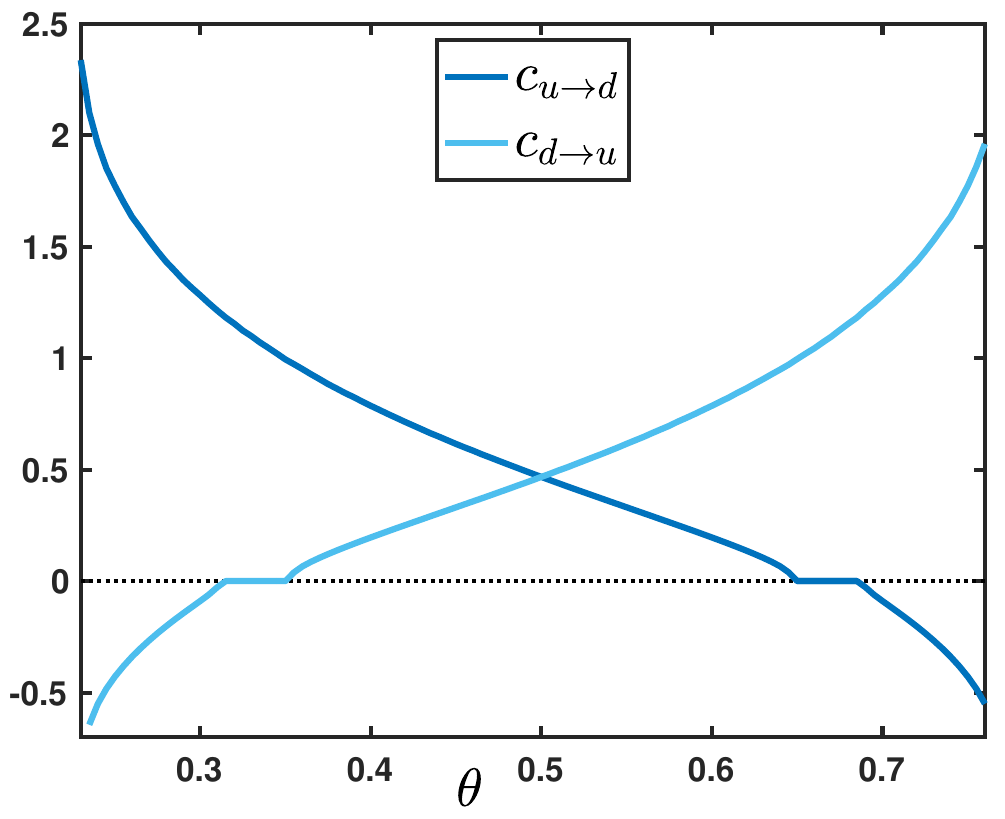}}
\subfigure[$q=0.5$.]{\includegraphics[width=.32\textwidth]{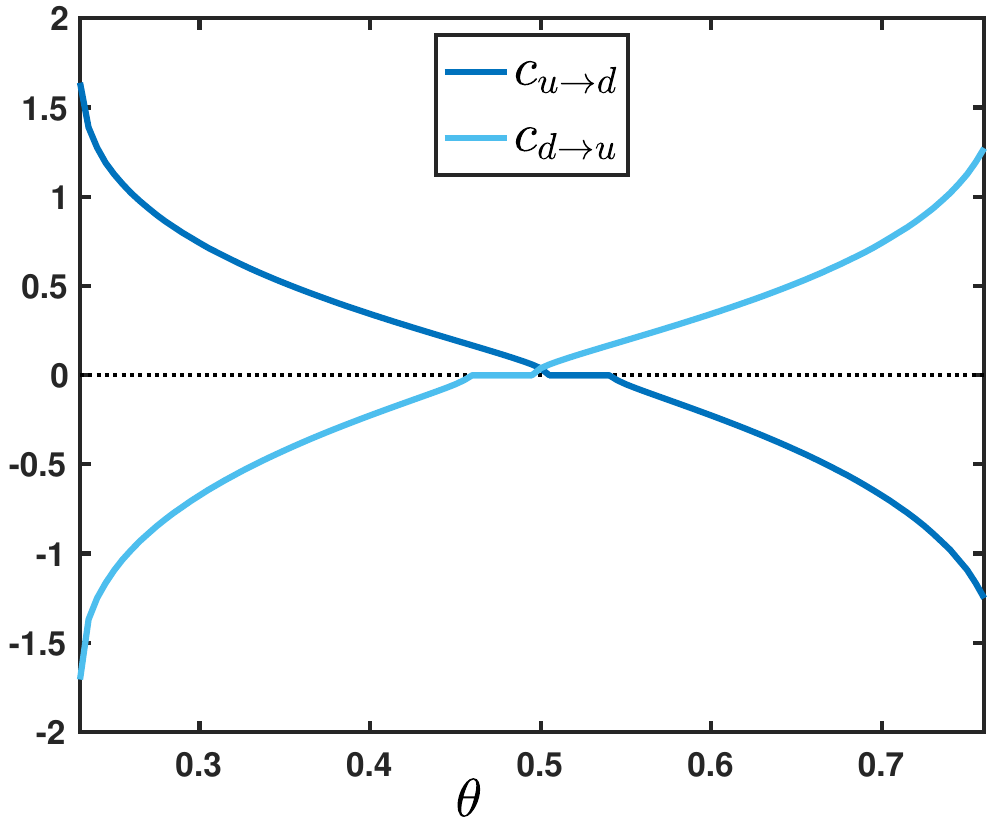}}
\subfigure[$q=0.65$.]{\includegraphics[width=.325\textwidth]{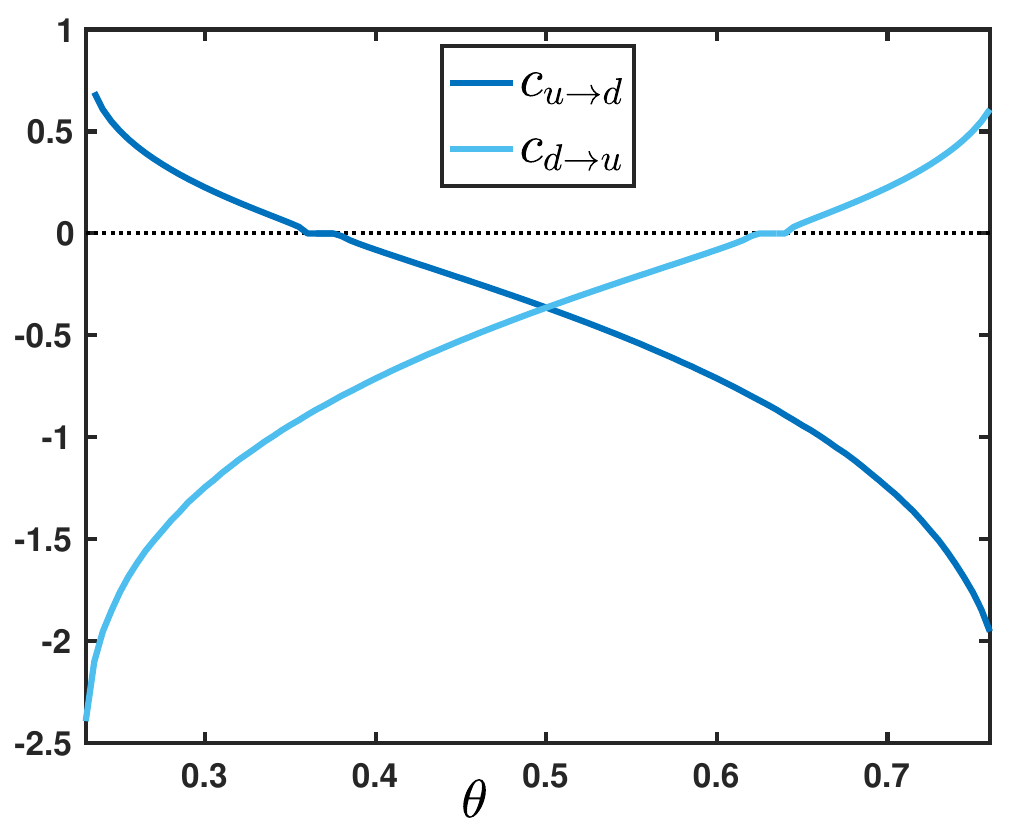}}
\subfigure[$\theta=0.35$.]{\includegraphics[width=.32\textwidth]{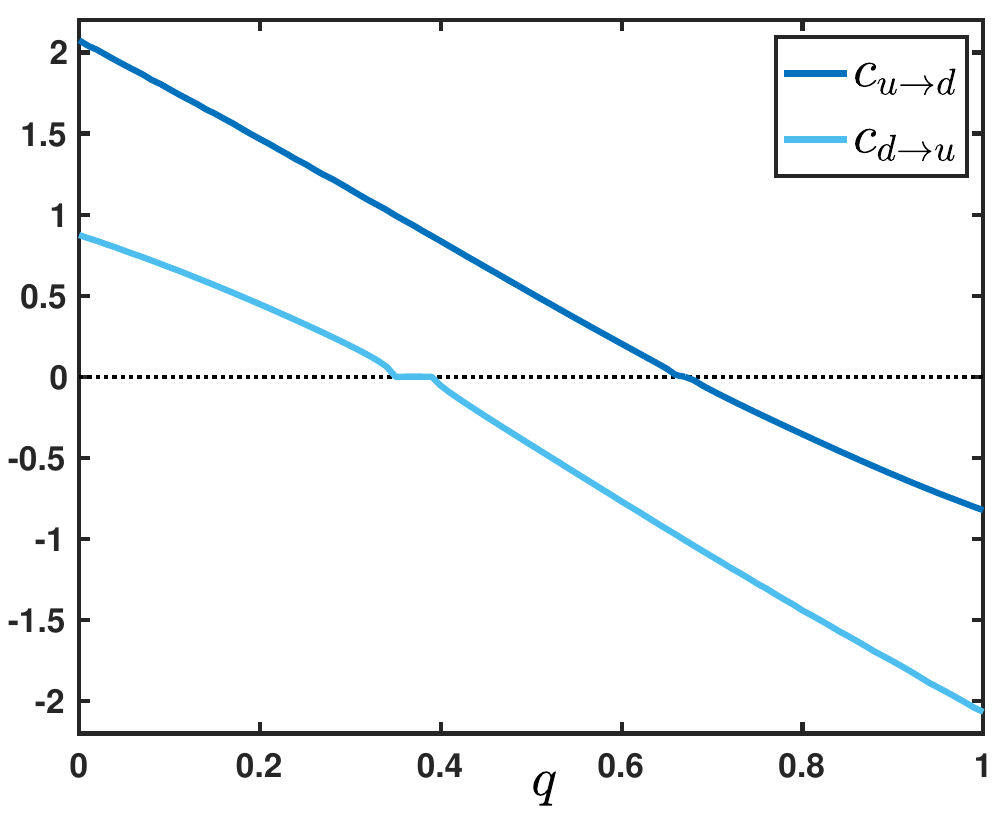}}
\subfigure[$\theta=0.5$.]{\includegraphics[width=.32\textwidth]{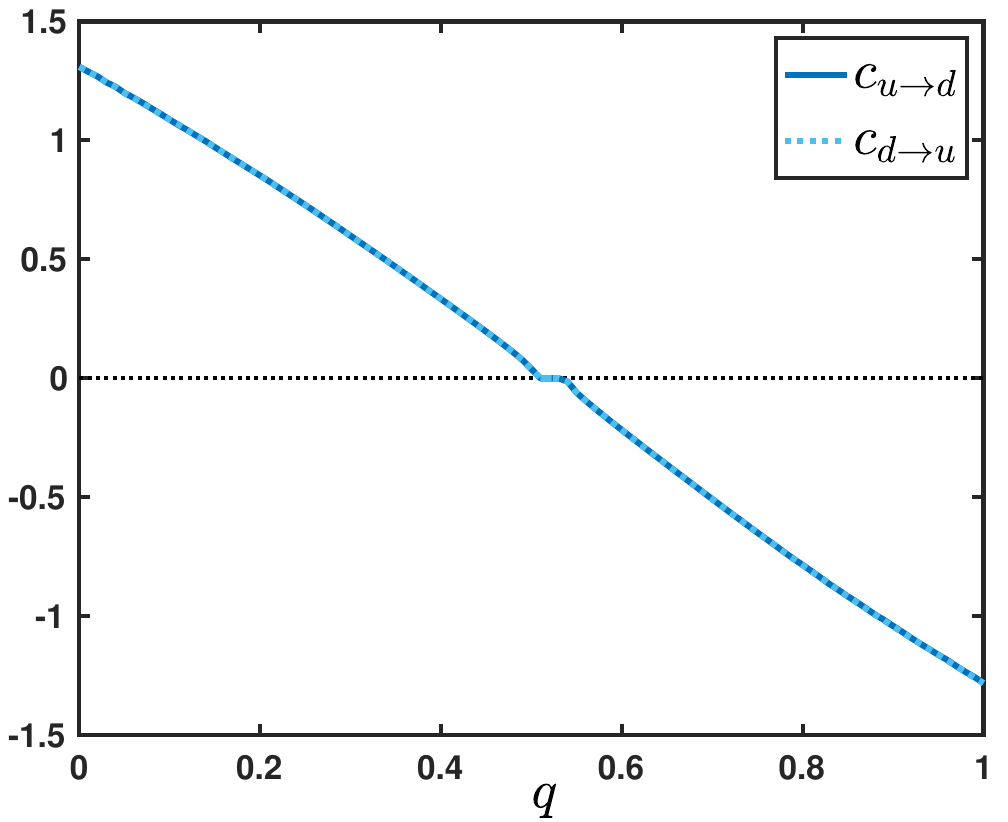}}
\subfigure[$\theta=0.65$.]{\includegraphics[width=.32\textwidth]{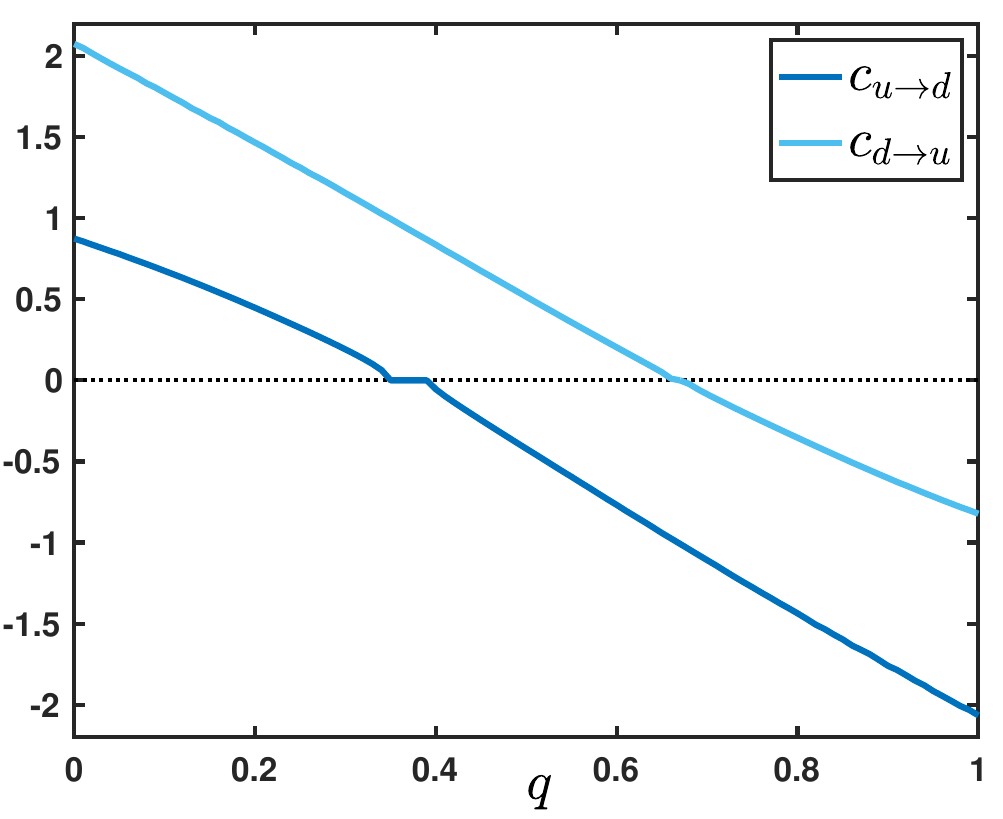}}
\caption{(a)-(b)-(c) Plots of the wave speeds $c_{u\to d}$ (blue) and $c_{d\to u}$ (light blue) as a function of $\theta\in(\theta_*(\mu),\theta^*(\mu))$ for several values of the parameter $q\in\left\{0.35,0.5,0.65\right\}$. (d)-(e)-(f) Plots of the wave speeds $c_{u\to d}$ (blue) and $c_{d\to u}$ (light blue) as a function of $q\in[0,1]$ for several values of the parameter $\theta\in\left\{0.35,0.5,0.65\right\}$. We observe that at $\theta=\frac{1}{2}$ both wave speeds are equal $c_{u\to d}=c_{d\to u}$ while for $\theta\lessgtr\frac{1}{2}$ one has $c_{d\to u}\lessgtr c_{u\to d}$ for all $q\in[0,1]$. The other parameters $\mu$ and $p$ are set to $\mu=16$ and $p=0.1$.}
 \label{fig:speed}
\end{figure}

As a consequence, we can use a result of Mallet-Paret~\cite{MP99} which asserts that for every $\Lambda\in\mathcal{P}$, there exists $c_{u\to d}\in\R$ (resp. $c_{d\to u}\in\R$) and a monotone decreasing (resp. montone increasing) solution $\Phi_{u\to d}$ of \eqref{TW} (resp. $\Psi_{d\to u}$ of \eqref{TWbis}). This wave speed $c_{u\to d}=c_{u\to d}(\Lambda)$ (resp. $c_{d\to u}=c_{d\to u}(\Lambda)$) is unique and depends continuously on $\Lambda\in\mathcal{P}$ and smoothly on $\Lambda\in\mathcal{P}$ when $c_{u\to d}(\Lambda)\neq0$ (resp. $c_{d\to u}(\Lambda)\neq0$). If $c_{u\to d}(\Lambda)\neq0$ (resp. $c_{d\to u}(\Lambda)\neq0$), then the profile $\Phi_{u\to d}$ (resp. $\Psi_{d\to u}$) is unique up to translations and satisfies $\Phi'_{u\to d}(\xi)<0$ (resp. $\Psi_{d\to u}'(\xi)>0$) for all $\xi\in\R$. Furthermore, when the wave speed is non zero, solutions $(v_j)_{j\in\Z}$ of \eqref{cont1dZ} starting from $(h_j^{u\to d})_{j\in\Z}$ or $(h_j^{d\to u})_{j\in\Z}$ asymptotically exponentially converge towards a shift of the traveling waves $\Phi_{u\to d}$ or $\Psi_{d\to u}$ as proved in \cite{CGW08,Z91}. More precisely, for each $\Lambda\in\mathcal{P}$ such that $c=c_{u\to d}(\Lambda)\neq0$, there exists a constant $\nu>0$ such that for all initial condition $(h_j^{u\to d})_{j\in\Z}$ the solution $(v_j)_{j\in\Z}$ of \eqref{cont1dZ} starting from $(h_j^{u\to d})_{j\in\Z}$ is globally defined, and there exists $\xi_0\in\R$ and $C>0$ such that
\bqs
\sup_{j\in\Z}\left| v_j(t)-\Phi_{u\to d}(j-ct+\xi_0) \right| \leq C e^{-\nu t}, \quad t\geq0.
\eqs
A similar result holds when $c=c_{d\to u}(\Lambda)\neq0$ for the solution $(v_j)_{j\in\Z}$ of \eqref{cont1dZ} starting from $(h_j^{d\to u})_{j\in\Z}$ with global asymptotic exponential convergence towards $\Psi_{d\to u}$. As a consequence, the sign of the above wave speed is determinant in understanding the long time dynamics of the network initialized with the sequences $(h_j^{u\to d})_{j\in\Z}$ or $(h_j^{d\to u})_{j\in\Z}$ as can be shown in Figure~\ref{fig:STP} where we report space-time plots of the solutions to \eqref{cont1dZ} starting from either of these two initial conditions. In simpler words, Figure~\ref{fig:STP} shows that considering any given point in the network (e.g., $j=0$ in the x-axis of all panels), the up-and down-state activity can either propagate bottom-up or top-down, depending on the network's parameters.

We illustrate in Figure~\ref{fig:speed} the dependence of the wave speeds $c_{u\to d}$ (blue curve) and $c_{d\to u}$ (light blue curve) when a parameter is varied while the other parameters are kept fixed. In Figures~\ref{fig:speed}(a)-(b)-(c), we vary continuously the threshold $\theta\in(\theta_*(\mu),\theta^*(\mu))$ and let $q\in\left\{0.35,0.5,0.65\right\}$ with $(\mu,p)=(16,0.1)$ fixed. We observe that the curve $\theta\mapsto c_{u\to d}(\theta)$ is monotone decreasing, while $\theta\mapsto c_{d\to u}(\theta)$ is monotone increasing and that the two curves $c_{u\to d}$ and $c_{d\to u}$ always intersect at $\theta=\frac{1}{2}$. We also note that both wave speeds identically vanish on some non empty interval. In the literature, this is referred to as \emph{propagation failure} or \emph{pinning} \cite{MP99,fath,keener87,EN93}. For values of the parameter in such intervals, one can typically show that there exists infinitely many stationary solutions to \eqref{cont1dZ} which asymptotically converge towards $x_{u/d}(\theta,\mu)$ and $x_{d/u}(\theta,\mu)$ at $\pm\infty$. Finally, we remark that when the fraction $q<\frac{1}{2}$ both wave speeds are positive for large intervals of the threshold whereas when $q>\frac{1}{2}$, we observe the opposite that is both wave speeds are positive for small intervals of the parameter $\theta$ and thus negative on large intervals. When $\theta$ is fixed and now $q\in[0,1]$ is varied, we observe that at $\theta=\frac{1}{2}$ both wave speeds are equal $c_{u\to d}(q)=c_{d\to u}(q)$ while for $\theta\lessgtr\frac{1}{2}$ one has $c_{d\to u}(q)\lessgtr c_{u\to d}(q)$ for all $q\in[0,1]$, see Figures~\ref{fig:speed}(d)-(e)-(f). Once again, as $q$ is varied, there also exists \emph{propagation failure} and non empty intervals where the wave speeds identically vanish.

\begin{figure}[t!]
\centering
\subfigure[Sign of $c_{u\to d}$.]{\includegraphics[width=.32\textwidth]{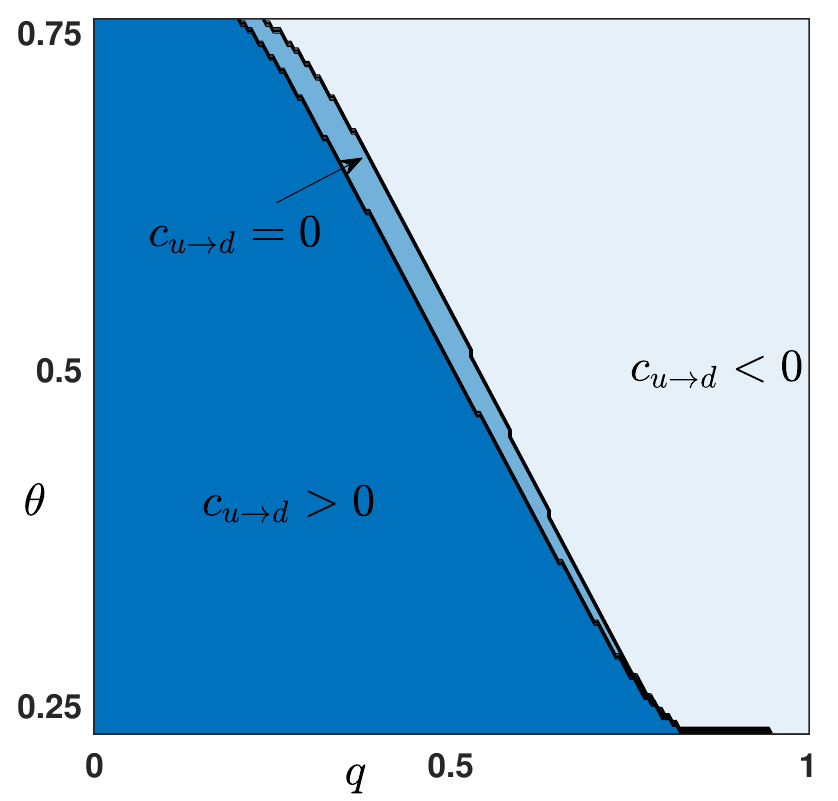}}
\subfigure[Sign of $c_{d\to u}$.]{\includegraphics[width=.32\textwidth]{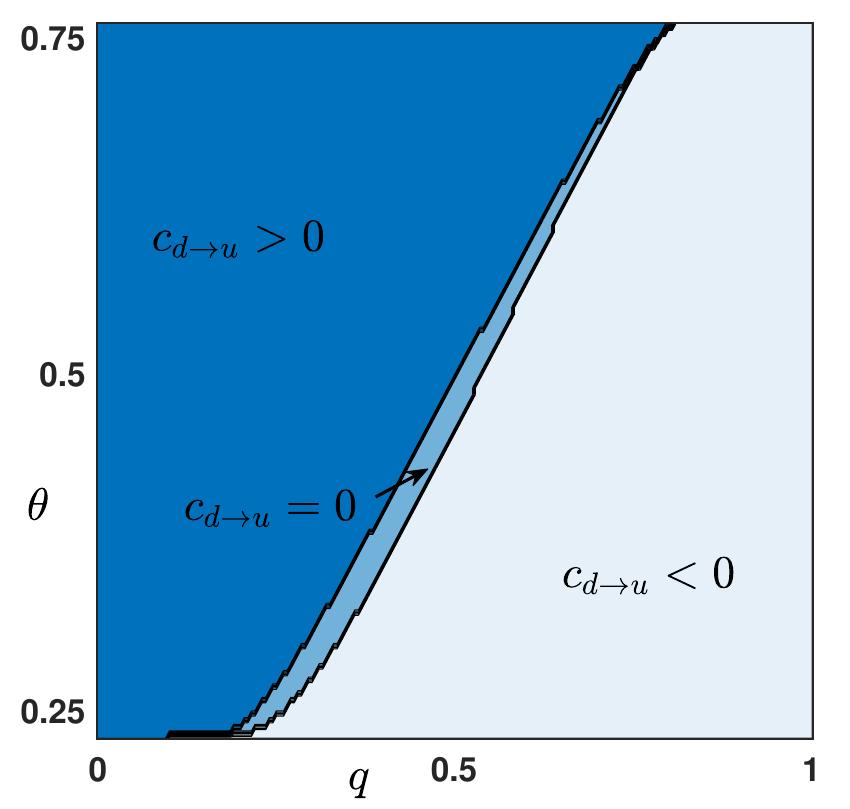}}
\subfigure[Sign of $(c_{u\to d},c_{d\to u})$.]{\includegraphics[width=.32\textwidth]{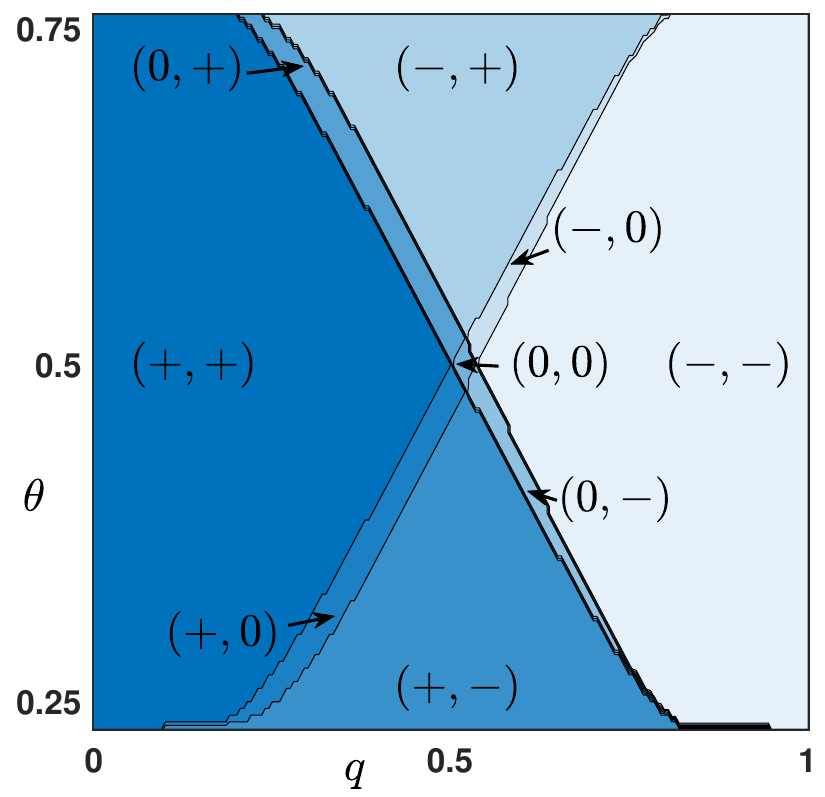}}
\caption{(a)-(b) Sign of the wave speeds $c_{u\to d}$ and $c_{d\to u}$  as a function of $q\in[0,1]$ and $\theta\in(\theta_*(\mu),\theta^*(\mu))$. Dark blue region corresponds to positive wave speeds while light blue region corresponds to negative speeds. Region of zero wave speeds is depicted in sky blue. (c) Combined sign of $(c_{u\to d},c_{d\to u})$ as a function of $q\in[0,1]$ and $\theta\in(\theta_*(\mu),\theta^*(\mu))$. The other parameters $\mu$ and $p$ are set to $(\mu,p)=(16,0.1)$.}
 \label{fig:signspeedtq}
\end{figure}

As we already emphasized, the sign of the wave speed is the key information in order to determine the long time behavior of the solutions of \eqref{cont1dZ} starting from the initial conditions $(h_j^{u\to d})_{j\in\Z}$ or $(h_j^{d\to u})_{j\in\Z}$. In Figure~\ref{fig:signspeedtq}(a)-(b), we show the numerically computed sign of the wave speeds $c_{u\to d}$ and $c_{d\to u}$  as a function of $(q,\theta)\in[0,1]\times(\theta_*(\mu),\theta^*(\mu))$. We observe three distinct regions where $c_{u\to d}$ (respectively $c_{d\to u}$) is positive (dark blue region), vanishes (sky blue region) or is negative (light blue region). It is also possible to visualize the joint sign $(c_{u\to d},c_{d\to u})$ as a function of $(q,\theta)\in[0,1]\times(\theta_*(\mu),\theta^*(\mu))$, see Figure~\ref{fig:signspeedtq}(c). All nine possible combinations of pair of signs $\left\{+,0,-\right\}$ correspond to nine disjoint regions in parameter space. We notably notice that for small values of $0\leq q\lesssim 0.2$, both wave speeds are positive independently of the value of the threshold, while for large values of $0.8 \lesssim q\leq 1$ both wave speeds are negative. This makes sense if one remembers that parameter $q$ reflects the strength of the feedback error correction in the predictive coding framework. Interestingly, there exists a small region near $(q,\theta)\sim(0.5,0.5)$ where $c_{u\to d}=c_{d\to u}=0$ and no propagation is possible.

\begin{figure}[t!]
\centering
\subfigure[$\theta=0.35$.]{\includegraphics[width=.32\textwidth]{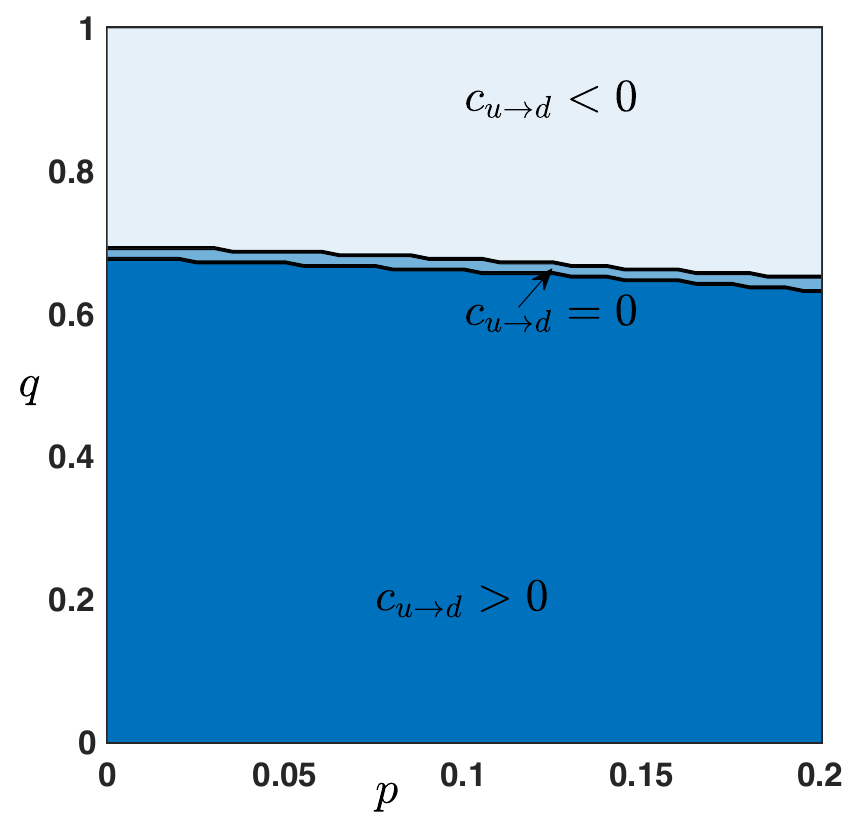}}
\subfigure[$\theta=0.5$.]{\includegraphics[width=.32\textwidth]{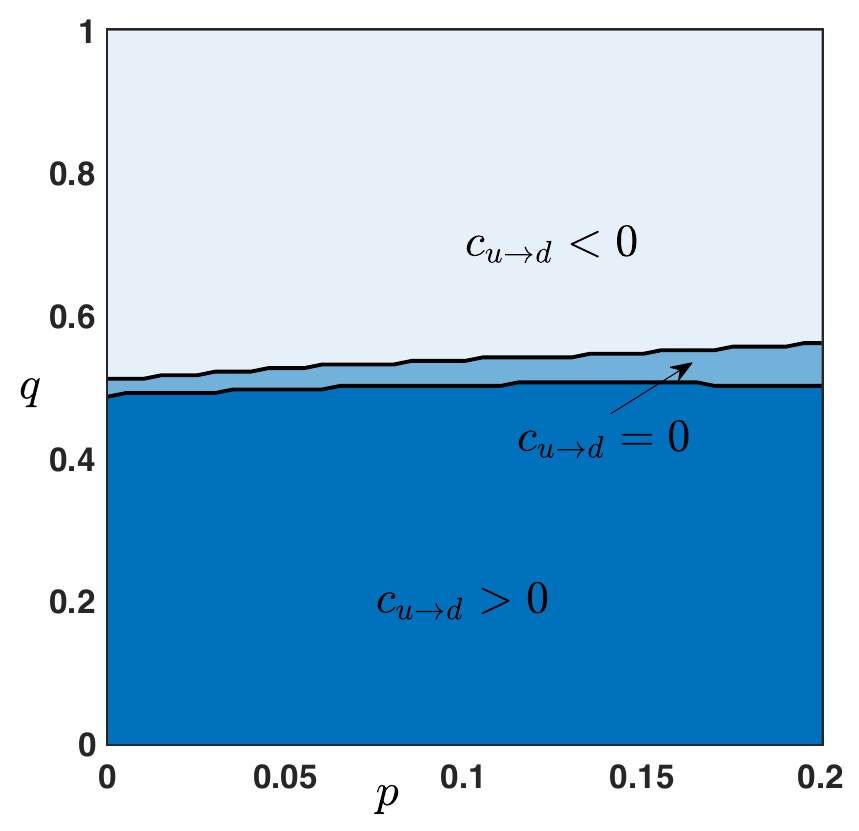}}
\subfigure[$\theta=0.65$.]{\includegraphics[width=.32\textwidth]{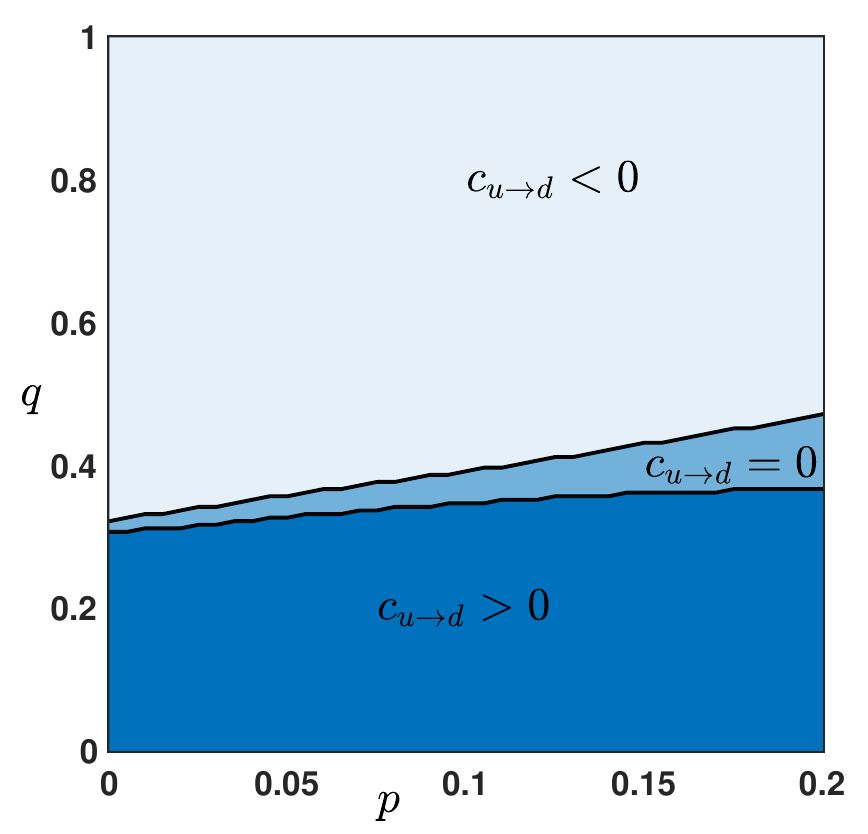}}
\subfigure[$\theta=0.35$.]{\includegraphics[width=.32\textwidth]{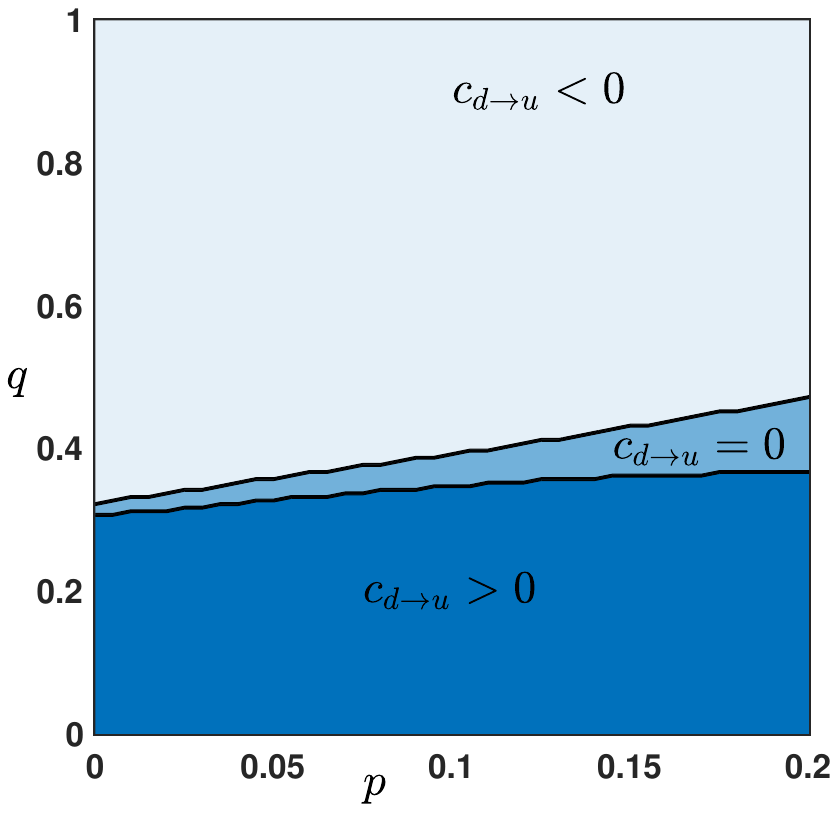}}
\subfigure[$\theta=0.5$.]{\includegraphics[width=.32\textwidth]{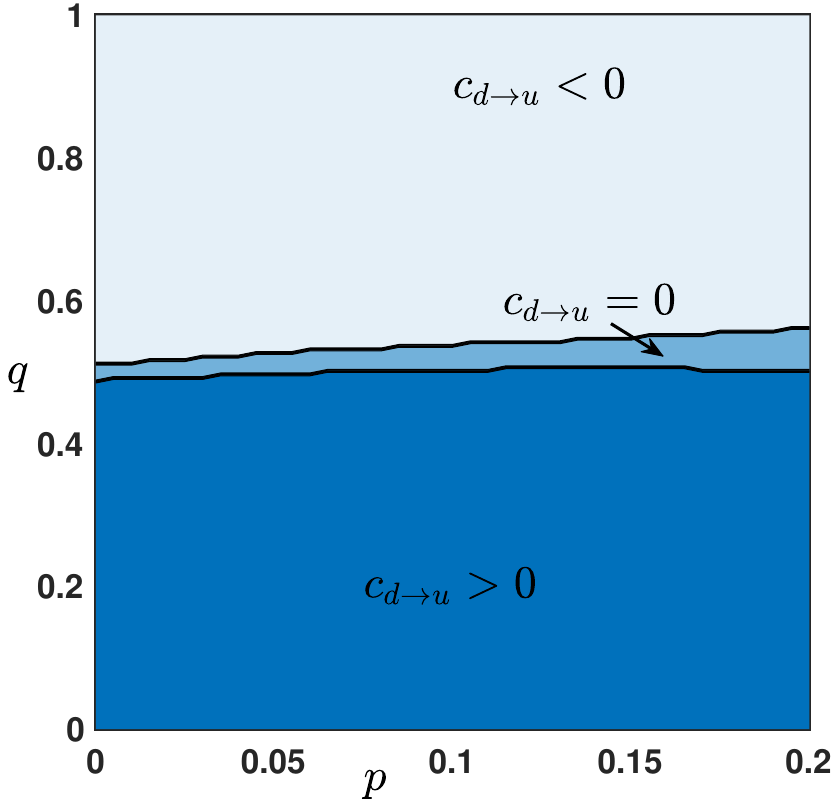}}
\subfigure[$\theta=0.65$.]{\includegraphics[width=.32\textwidth]{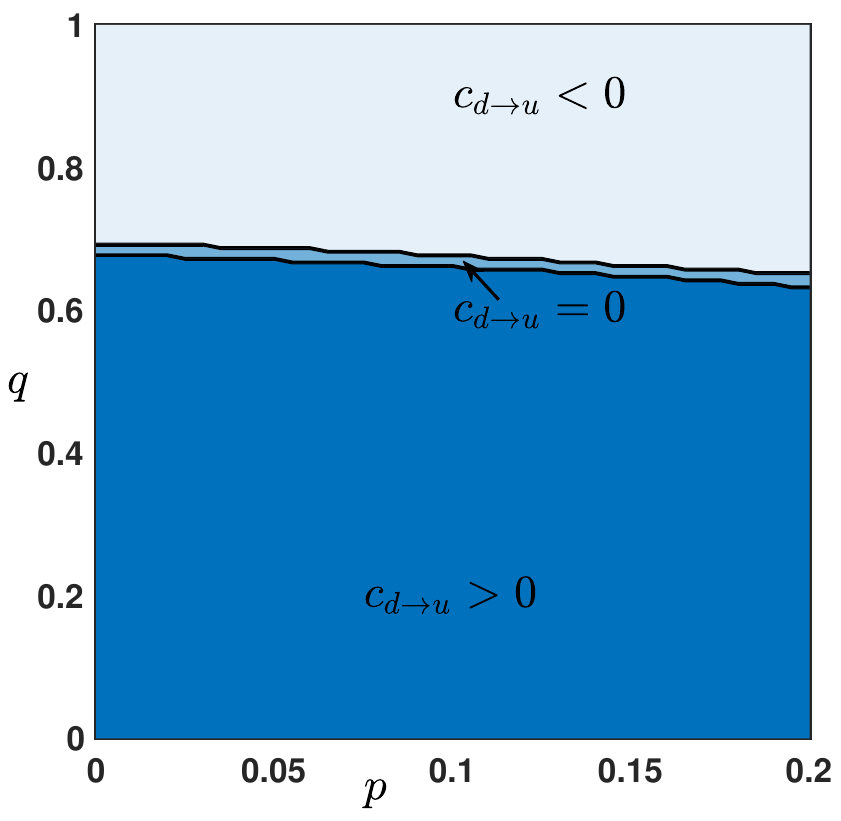}}
\caption{Sign of the wave speed $c_{u\to d}$ and $c_{d\to u}$ as a function of $(p,q)\in\left[0,\frac{4}{4+\mu}\right]\times[0,1]$ for several values of the parameter $\theta\in(\theta_*(\mu),\theta^*(\mu))$. Dark blue region corresponds to positive wave speeds while light blue region corresponds to negative speeds. Region of zero wave speeds is depicted in sky blue. The parameter $\mu$ is set to $\mu=16$.}
 \label{fig:signspeedpq}
\end{figure}

We have also numerically computed the sign of the wave speeds when this time $p$ and $q$ are varied in $\left[0,\frac{4}{4+\mu}\right]\times[0,1]$ with $\theta$ and $\mu$ being kept fixed (see Figure~\ref{fig:signspeedpq}). We most notably remark that the sign of the wave speeds, and thus the corresponding regions in the parameter space, is almost independent of $p$. For small values of $\theta$ wave speeds are positive in larger regions of parameter space whereas for large values of $\theta$ wave speeds are positive in smaller regions. It is also interesting to remark that the larger $\theta$ is, the larger is the region of parameter space where wave speeds vanish. Finally, we notice the following symmetry
\bqs
\mathrm{sgn}~ c_{u\to d}(\theta,\mu,p,q)= \mathrm{sgn}~ c_{d\to u}(1-\theta,\mu,p,q), \quad \forall (\theta,\mu,p,q)\in\mathcal{P},
\eqs
and actually one has the stronger fact that
\bqs
c_{u\to d}(\theta,\mu,p,q)=  c_{d\to u}(1-\theta,\mu,p,q), \quad \forall (\theta,\mu,p,q)\in\mathcal{P}.
\eqs

\subsection{Wave initiation on a semi-infinite depth network -- The bottom-up case}

\begin{figure}[t!]
\centering
\subfigure[Stagnation: $s_0<s_0^*$.]{\includegraphics[width=.35\textwidth]{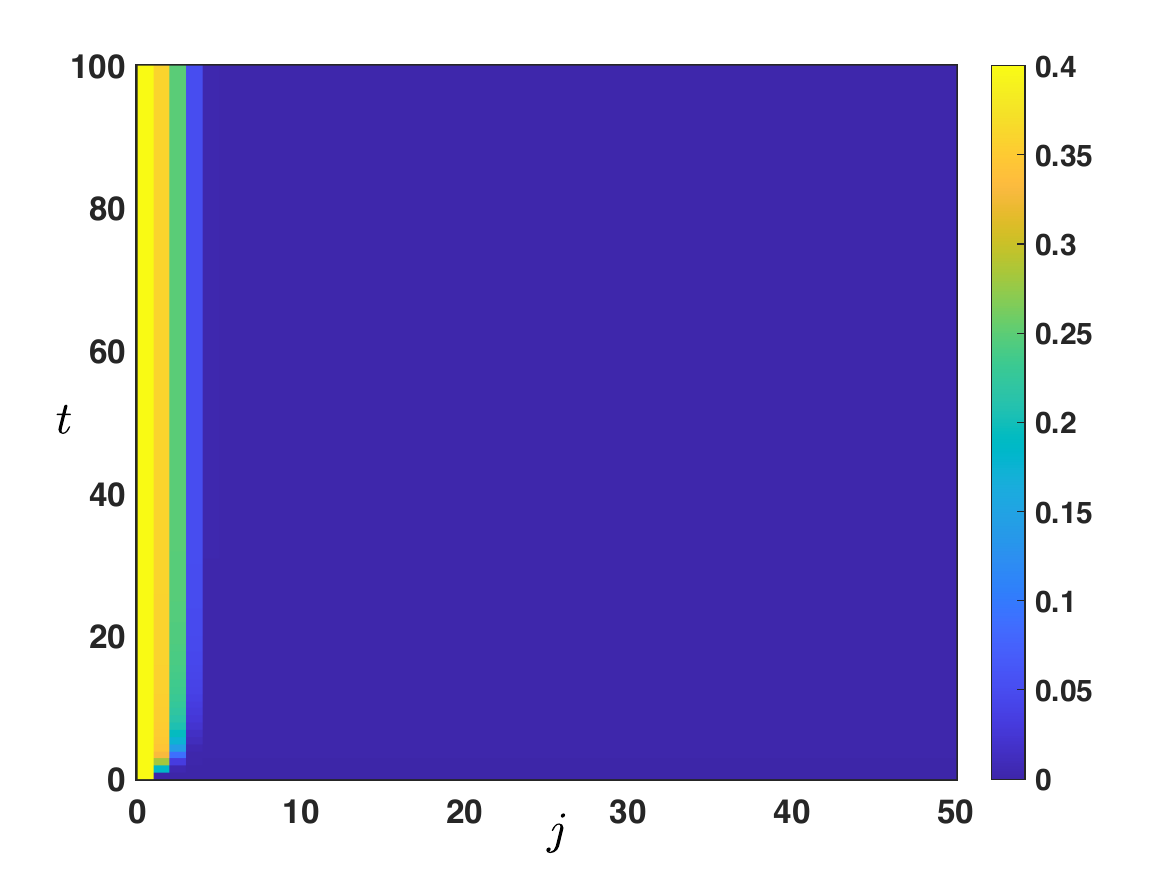}}\hspace{1cm}
\subfigure[Propagation: $s_0>s_0^*$.]{\includegraphics[width=.35\textwidth]{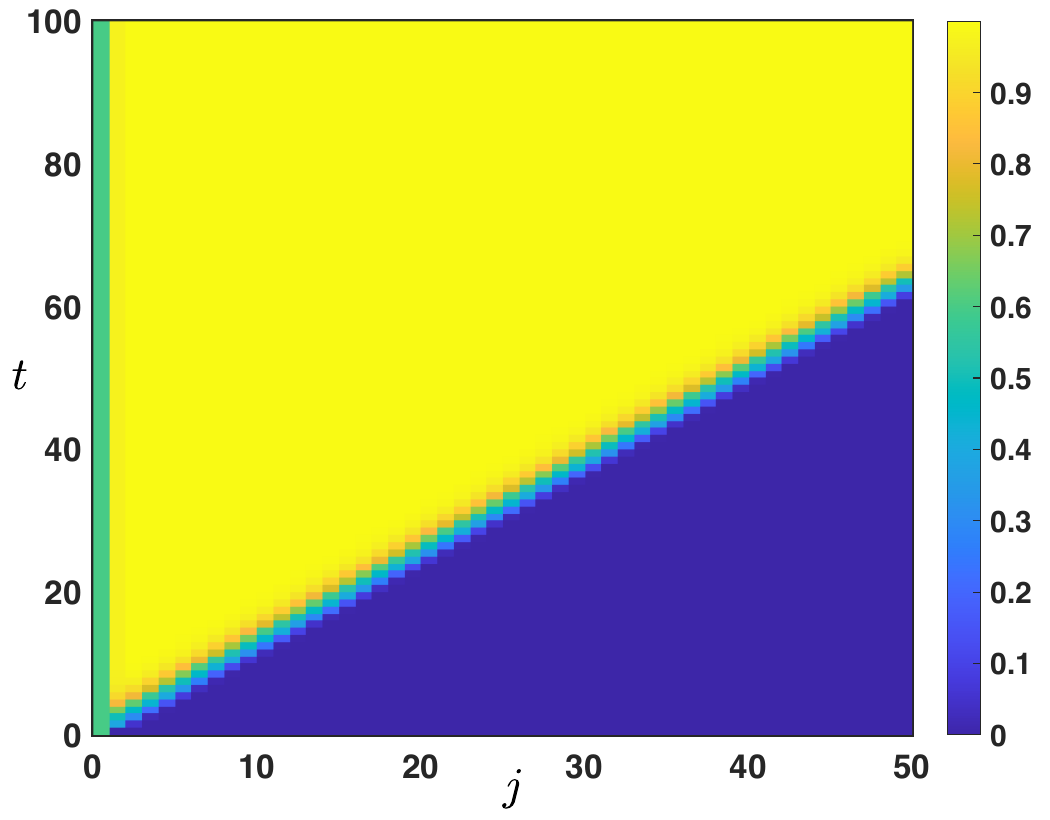}}
\caption{Space-time plot of the solutions of \eqref{cont1dN}-\eqref{IC1dN} when $s_0(t)=s_0$ for all $t>0$. Depending on the amplitude $s_0 \lessgtr s_0^*$ for some threshold $s_0^*$, the solution can either uniformly converge towards the stationary solution $\boldsymbol{x}^d(s_0)$ (stagnation) or locally uniformly converge towards the stationary solution $\boldsymbol{x}^u(s_0)$ (propagation). }
 \label{fig:STPbord}
\end{figure}

We now relax the assumption that the network under consideration is bi-infinite and work instead with a semi-infinite network, that is we consider
\bqq\label{cont1dN}
\forall t>0, \quad \left\{\begin{split}
v_j'(t) & =\mathscr{N}\left(v_{j-1}(t),v_j(t),v_{j+1}(t)\right), \quad j \geq 1, \\
v_0(t)& = s_0(t).
\end{split}
\right.
\eqq
Our primary interest will be to understand the influence of the time varying external input $s_0$ on the long time dynamics of the network. Throughout this section, we will always assume that the network is initialized to be in the down state of the system, that is, we always assume that the initial condition is set to
\bqq
\label{IC1dN}
v_j(0)=u_d(\theta,\mu), \quad j\geq 1,
\eqq
and the parameters $(\theta,\mu,p,q)\in\mathcal{P}$. In such a setting, the external input can be interpreted as an external perturbation (i.e., a sensory stimulus) provided during a resting state, and we would like to somehow quantify the \emph{strength} of the perturbation which may generate a propagation upward in the network. We will essentially consider two forms of external input. In a first step, we will work under the assumption that the external input is constant $s_0(t)=s_0\in\R$, meaning that it is constantly fixed at a given value for all positive times and we will always assume that $s_0\geq x_d(\theta,\mu)$. In a second step, we will rather consider the case where the external input is presented constantly at a fixed amplitude $s_0$ during an initial time window $[0,\tau]$ for $\tau>0$ and then set to the down state $x_d(\theta,\mu)$ for later times, that is  $s_0(t)=s_0\mathds{1}_{[0,\tau]}(t)+x_d(\theta,\mu)\mathds{1}_{(\tau,+\infty)}(t)$ for all $t>0$. In this second case, we shall always assume that $s_0>x_d(\theta,\mu)$.

\subsubsection{Threshold phenomena for constant external input}

We study the long time behavior of the solutions of \eqref{cont1dN}-\eqref{IC1dN} with constant external input, that is throughout we set
\bqs
\forall t>0, \quad s_0(t)=s_0.
\eqs
Regarding the amplitude $s_0$ of the external input, we already make the intuitive remark that it needs to be large enough in order for the corresponding solution to be activated and pushed to the up state of the network, and without loss of generality we will always assume that it is larger than the down state, that is $s_0\geq x_d(\theta,\mu)$. Based on comparison principle techniques (see Lemma~\ref{lemCPN}), we have that if $s_0\in[x_d(\theta,\mu),x_m(\theta,\mu)]$, then $v_j(t)\in[x_d(\theta,\mu),x_m(\theta,\mu)]$ for all $t>0$ and $j\geq1$ such that one needs $s_0>x_m(\theta,\mu)$ in order for the solution to eventually converge to the up state. When $s_0\in[x_d(\theta,\mu),x_m(\theta,\mu)]$, one typically observes (see Figure~\ref{fig:STPbord}(a)) that the solution uniformly converges towards a stationary solution of \eqref{cont1dN}, that is a sequence $(x_j)_{j\in\N}$ solution of the following form
\bqq
\left\{
\begin{split}
0&= \mathscr{N}\left(x_{j-1},x_j,x_{j+1}\right), \quad j \geq 1, \\
x_0&=s_0,\\
x_j& \underset{j\rightarrow+\infty}{\longrightarrow} x_d(\theta,\mu).
\end{split}
\right.
\label{statsolNd}
\eqq
Let us also note that in the special case $s_0=x_{k}(\theta,\mu)$ with $k\in\left\{d,m,u\right\}$, the constant solution $(x_{k}(\theta,\mu))_{j\in\N}$ is a stationary homogeneous solution of \eqref{cont1dN}.

\begin{figure}[t!]
\centering
\subfigure[$\theta=0.35$.]{\includegraphics[width=.32\textwidth]{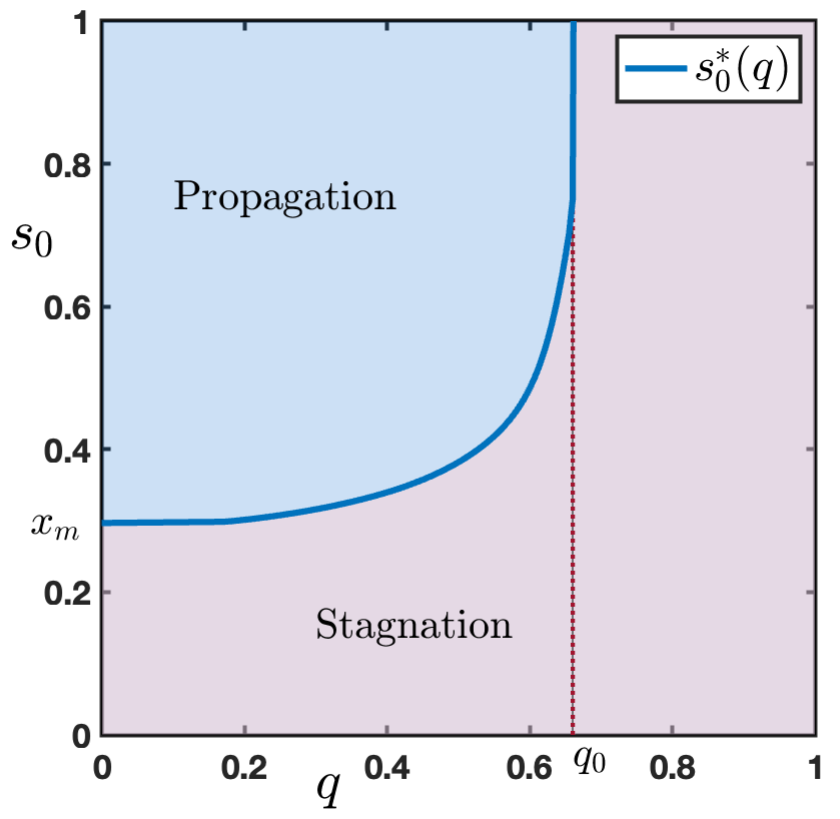}}
\subfigure[$\theta=0.5$.]{\includegraphics[width=.32\textwidth]{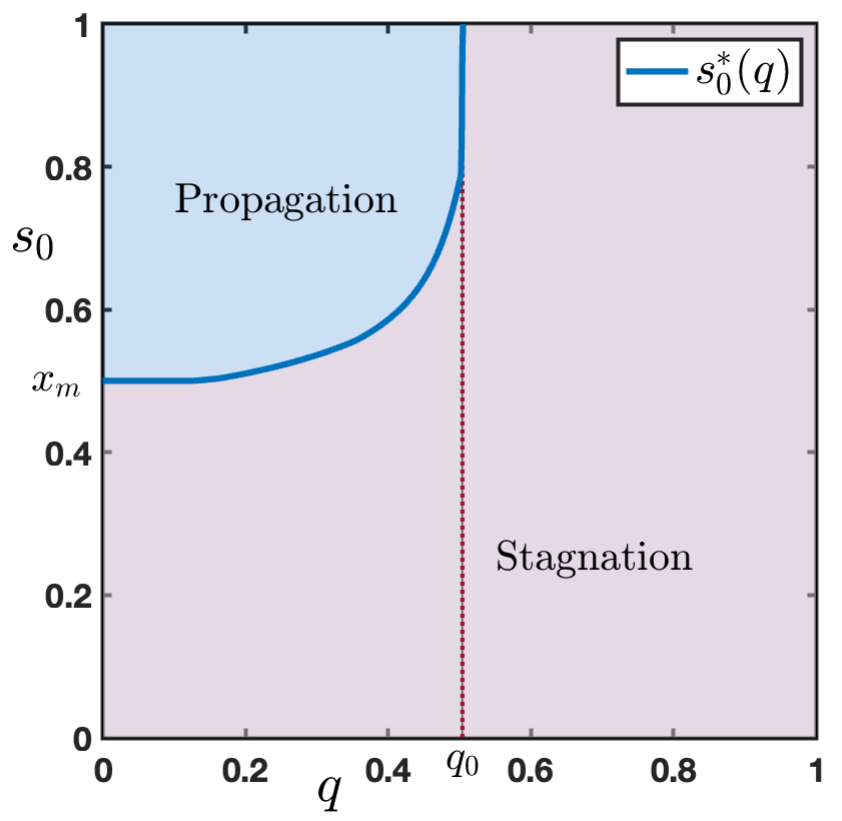}}
\subfigure[$\theta=0.65$.]{\includegraphics[width=.32\textwidth]{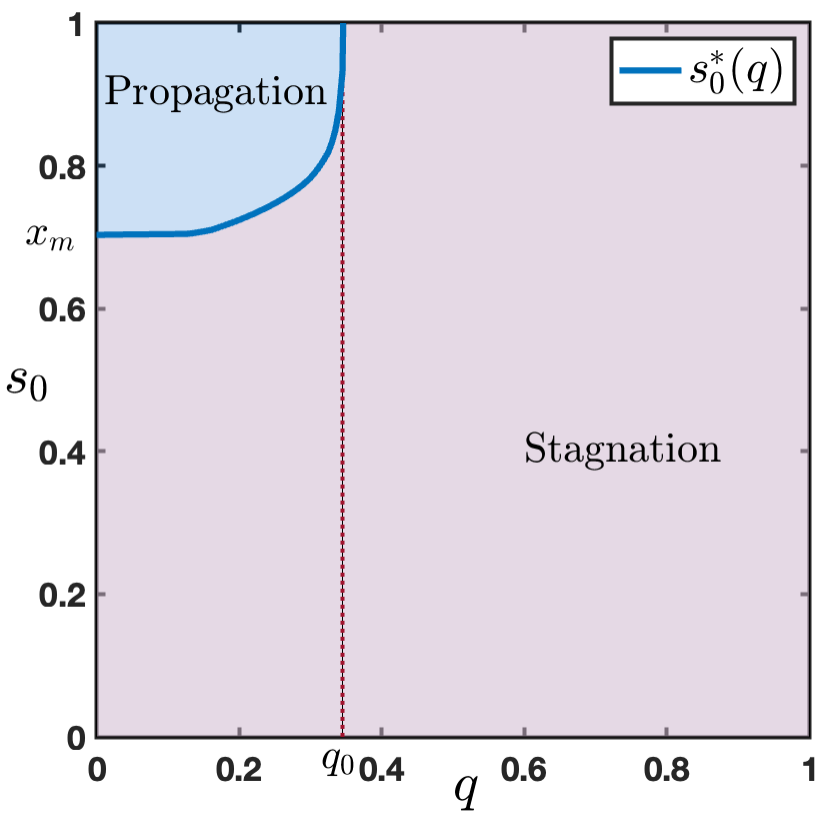}}
\caption{Plot of the sharp threshold $s_0^*(q)$ for $q\in[0,q_0)$ for several values of $\theta\in\left\{0.35,0.5,0.65\right\}$. The curve $\left\{ s_0^*(q)~|~q\in[0,q_0) \right\}$ delimits two regions in the rectangle $[0,1]\times[x_d(\theta,\mu),+\infty)$. In the region $\left\{(q,s_0) ~|~ s_0> s_0^*(q)\text{ and }q\in[0,q_0) \right\}$ (blue) there is propagation while in the region $\left\{(q,s_0) ~|~ s_0< s_0^*(q) \text{ and }q\in[0,q_0) \right\}\cup \left\{(q,s_0) ~|~ q_0\leq q \leq 1 \text{ and } s_0 \geq x_d(\theta,\mu) \right\}$ (violet) there is stagnation. Other values of the parameters are set to $(\mu,p)=(16,0.1)$.}
 \label{fig:s0qvary}
\end{figure}

For future reference, we shall denote by $\boldsymbol{x}^d(s_0)=(x_j^d)_{j\in\N}$ a stationary solution of \eqref{statsolNd} where we make explicit the dependence on $s_0$. Another intuitive remark can also be made at this stage. From the preceding section, in the bi-infinite case, we have seen that upward propagation of the up state to the down state can only occur with a positive wave speed $c_{u\to d}>0$. As a consequence, for the up state to propagate upward in the network in the semi-infinite case, one needs for $\Lambda\in\mathcal{P}$ to be in a region of parameter space where $c_{u\to d}(\Lambda)>0$. Let us note that in this case (see Figure~\ref{fig:STPbord}(b)) the solution locally uniformly converges towards a stationary solution of \eqref{cont1dN} of the form 
\bqq
\left\{
\begin{split}
0&= \mathscr{N}\left(x_{j-1},x_j,x_{j+1}\right), \quad j \geq 1, \\
x_0&=s_0,\\
x_j& \underset{j\rightarrow+\infty}{\longrightarrow} x_u(\theta,\mu).
\end{split}
\right.
\label{statsolNu}
\eqq
Compared to the stationary solution of \eqref{statsolNd} the key difference is that now the asymptotic state at $j=+\infty$ is given by the up state $x_u(\theta,\mu)$. We shall also denote by $\boldsymbol{x}^u(s_0)=(x_j^u)_{j\in\N}$ a stationary solution of \eqref{statsolNu}.

\begin{figure}[t!]
\centering
\subfigure[$q=0.35$.]{\includegraphics[width=.32\textwidth]{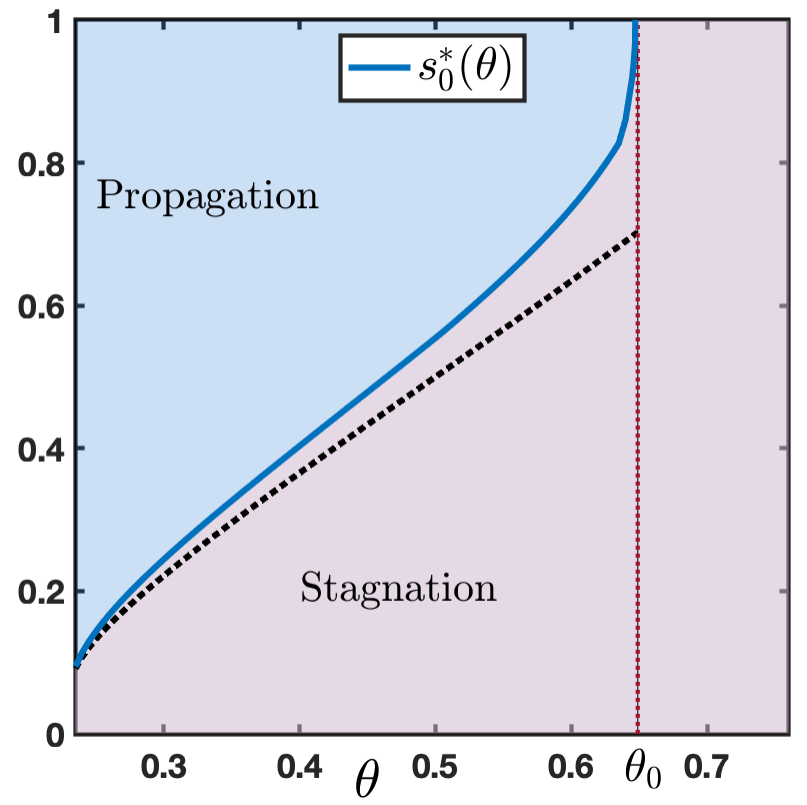}}
\subfigure[$q=0.5$.]{\includegraphics[width=.32\textwidth]{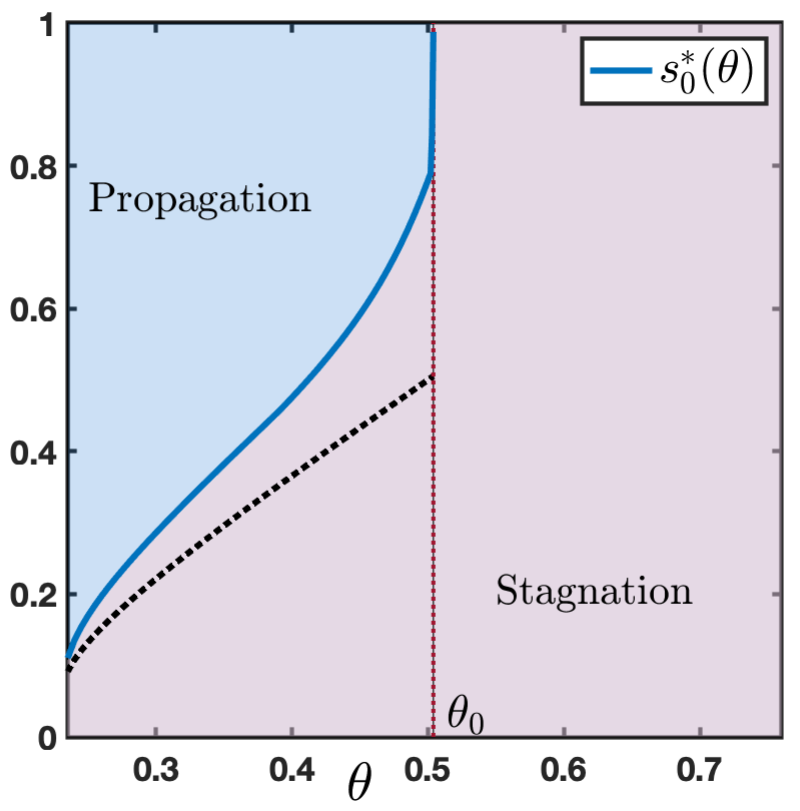}}
\subfigure[$q=0.65$.]{\includegraphics[width=.32\textwidth]{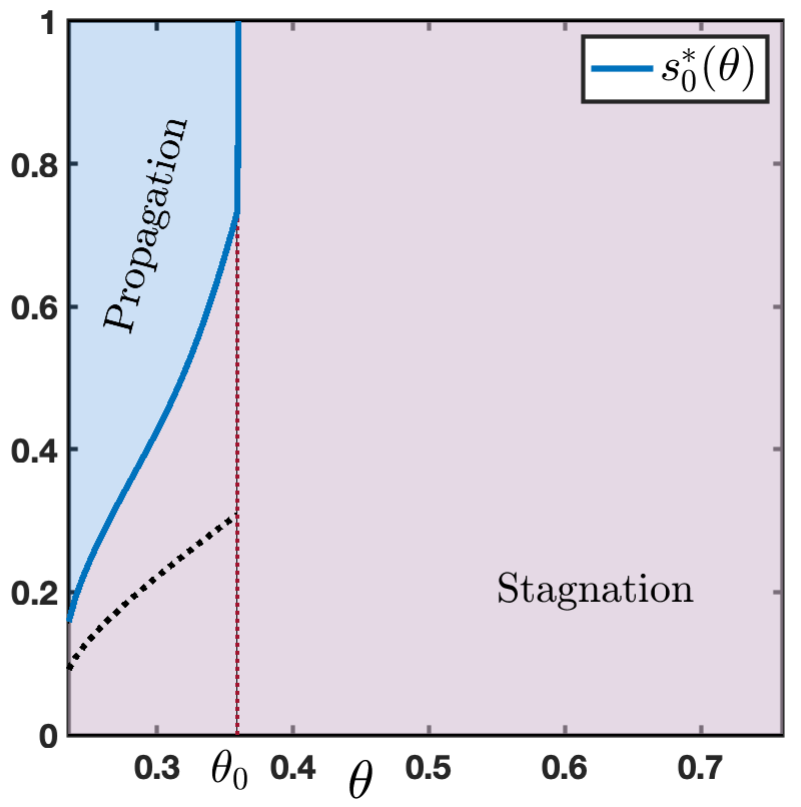}}
\caption{Plot of the sharp threshold $s_0^*(\theta)$ for $\theta\in(\theta_*(\mu),\theta_0)$ for several values of $q\in\left\{0.35,0.5,0.65\right\}$. The curve $\left\{ s_0^*(\theta)~|~\theta\in(\theta_*(\mu),\theta_0) \right\}$ delimits two regions in the rectangle $(\theta_*(\mu),\theta^*(\mu))\times[x_d(\theta,\mu),+\infty)$. In the region $\left\{(\theta,s_0) ~|~ s_0> s_0^*(\theta)\text{ and }\theta\in(\theta_*(\mu),\theta_0) \right\}$ (blue) there is propagation while in the region $\left\{(\theta,s_0) ~|~ s_0< s_0^*(\theta) \text{ and }\theta\in(\theta_*(\mu),\theta_0) \right\}\cup \left\{(\theta,s_0) ~|~ \theta_0\leq \theta < \theta^*(\mu) \text{ and } s_0 \geq x_d(\theta,\mu) \right\}$ (violet) there is stagnation. The dashed black line presents the unstable stationary state $x_m(\theta,\mu)$ showing that $s_0^*(\theta)\geq x_m(\theta,\mu)$ for $\theta\in(\theta_*(\mu),\theta_0)$. Other values of the parameters are set to $(\mu,p)=(16,0.1)$.}
 \label{fig:s0Thetavary}
\end{figure}

As a consequence, for a given $s_0\geq x_d(\theta,\mu)$, we say that there is \emph{propagation} if $\boldsymbol{v}=(v_j)_{j\in\N}$ the solution of \eqref{statsolNd} satisfies 
\bqs
\boldsymbol{v}(t) \underset{t\rightarrow+\infty}{\longrightarrow} \boldsymbol{x}^u(s_0) \text{ locally uniformly on } \N\footnote{By local uniform convergence on $\N$, we mean that for each $0\leq k < \ell <+\infty$, one has $\underset{{k\leq j \leq \ell}}{\sup}\left| v_j(t)- x_j ^u\right|\underset{t\rightarrow+\infty}{\longrightarrow}0$.},
\eqs
while there is \emph{stagnation} if the corresponding solution satisfies
\bqs
\boldsymbol{v}(t) \underset{t\rightarrow+\infty}{\longrightarrow} \boldsymbol{x}^d(s_0) \text{ uniformly on } \N\footnote{By uniform convergence on $\N$, we mean that $\underset{j\in\N}{\sup}\left| v_j(t)- x_j ^d\right|\underset{t\rightarrow+\infty}{\longrightarrow}0$.}.
\eqs

In plain English, \emph{stagnation} corresponds to a state in which the activity is strictly related to the instantaneous presence of the input and doesn't propagate through the network. On the other hand, \emph{propagation} occurs when the activity, initiated by the input, spreads bottom-up across the nodes of the network and keeps traveling through the nodes after the input ceases. 

As we have already seen, for all $s_0\in[x_d(\theta,\mu),x_m(\theta,\mu)]$ stagnation occurs. Our numerical investigations (see Figures~\ref{fig:s0qvary} and~\ref{fig:s0Thetavary}) show the existence of a sharp threshold between stagnation and propagation as $s_0$ is further increased. More precisely, there exists $s_0^*>x_m(\theta,\mu)$ such that the following dichotomy holds:
\begin{itemize}
\item for all $s_0\in[x_d(\theta,\mu),s_0^*)$ there is stagnation;
\item for all $s_0>s_0^*$ there is propagation.
\end{itemize}
In Figure~\ref{fig:s0qvary} we plot the value of the sharp threshold $s_0^*$ as a function of $q$ for several values of $\theta\in\left\{0.35,0.5,0.65\right\}$ while the other values of the parameters are being fixed. The sharp threshold is well defined for all values of $q\in[0,q_0)$ where we denote by $q_0\in(0,1)$ the smallest value of $q\in[0,1]$ for which $c_{u\to d}=0$, that is $q_0=\min\left\{ q\in[0,1]~|~ c_{u\to d}=0\right\}$. We remark, as expected, that $s_0^*(q)\geq x_m(\theta,\mu)$ for all $q\in[0,q_0)$. We also note that $s_0^*(q)\rightarrow+\infty$ as $q\rightarrow q_0^-$. We report in Figure~\ref{fig:s0Thetavary} similar behaviors when this time the sharp threshold $s_0^*$ is considered as a function of $\theta$ for several values of $q\in\left\{0.35,0.5,0.65\right\}$ while the other values of the parameters are being fixed. The sharp threshold is well defined for all values of $\theta\in(\theta_*(\mu),\theta_0)$ where we denote by $\theta_0\in(\theta_*(\mu),\theta^*(\mu))$ the smallest value of $\theta$ for which $c_{u\to d}=0$, that is $\theta_0=\min\left\{ \theta\in(\theta_*(\mu),\theta^*(\mu))~|~ c_{u\to d}=0\right\}$. Once again, we remark that $s_0^*(\theta)\geq x_m(\theta,\mu)$ for all $\theta\in(\theta_*(\mu),\theta_0)$ and $s_0^*(\theta)\rightarrow+\infty$ as $\theta\rightarrow \theta_0^-$.

\subsubsection{Threshold phenomena for flashed external input}

\begin{figure}[t!]
\centering
\subfigure[Propagation failure.]{\includegraphics[width=.32\textwidth]{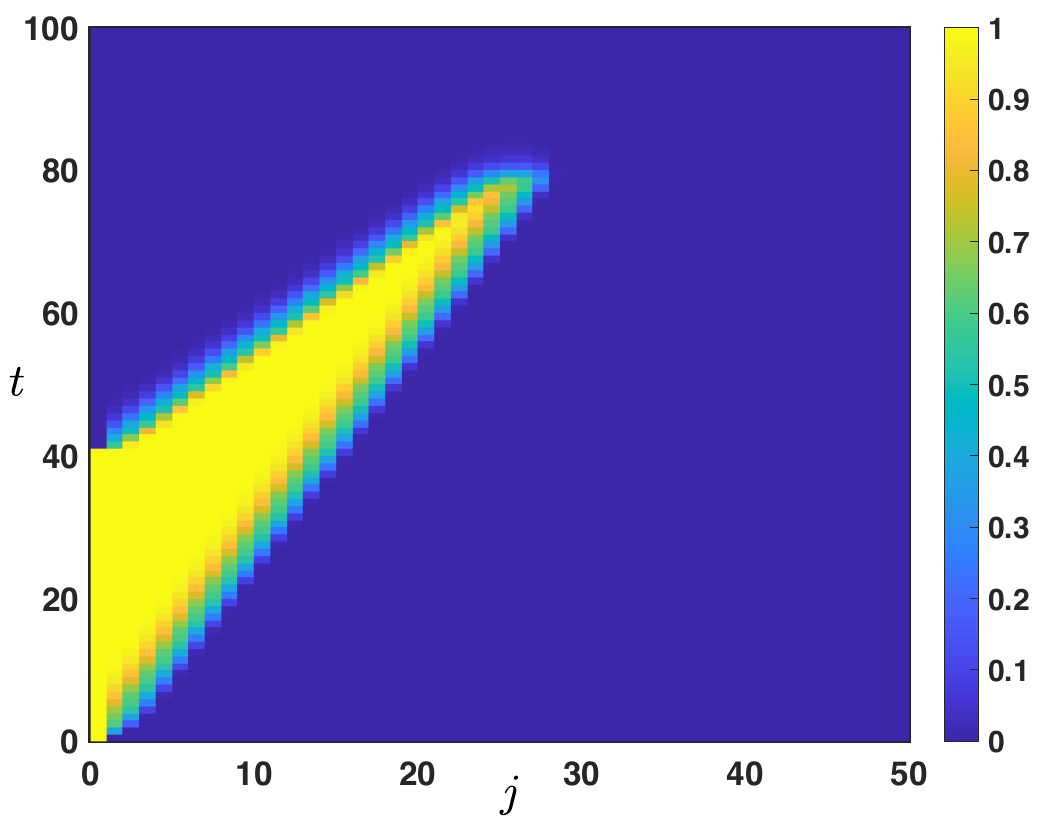}}
\subfigure[Stacked propagation.]{\includegraphics[width=.32\textwidth]{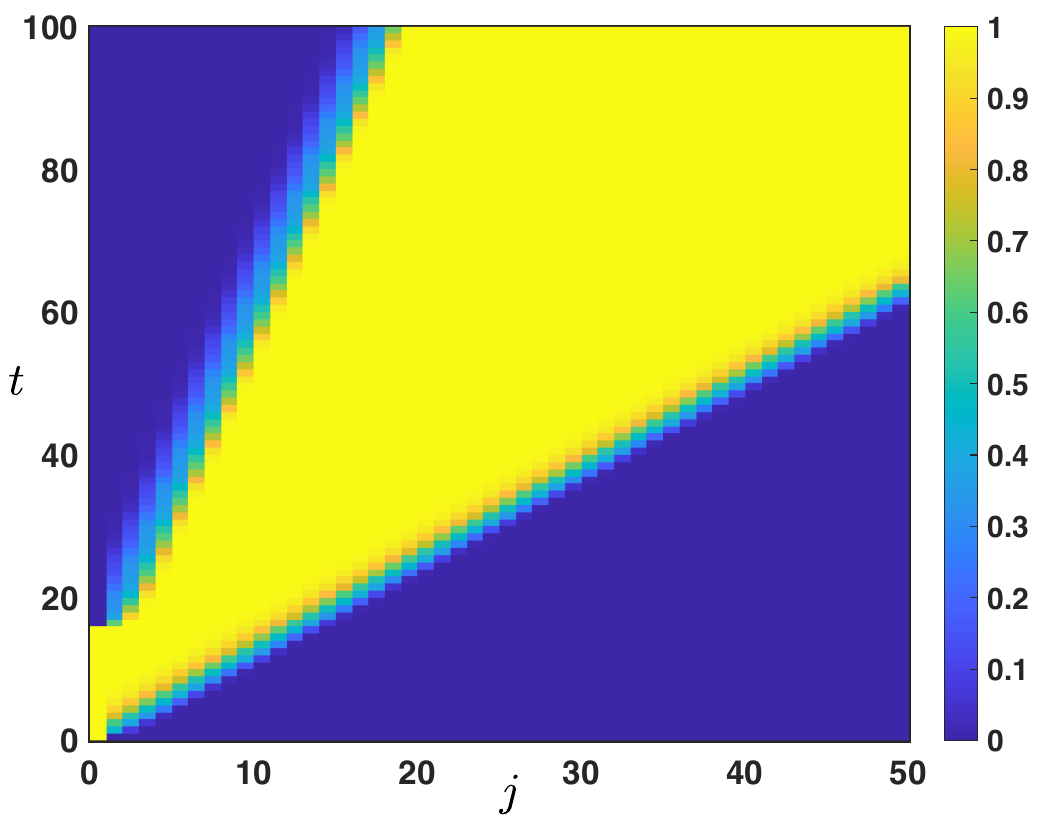}}
\subfigure[Front propagation.]{\includegraphics[width=.32\textwidth]{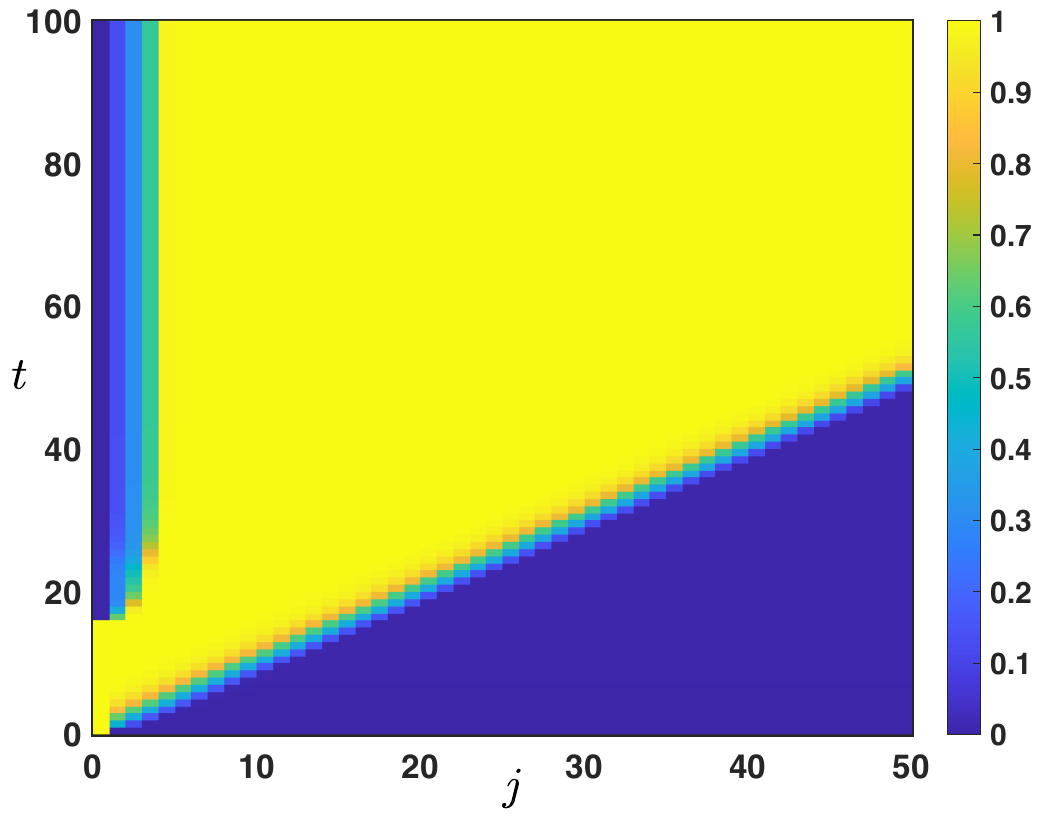}}
\caption{Space-time plot of the solutions of \eqref{cont1dN}-\eqref{IC1dN} when $s_0(t)=x_u(\theta,\mu)\mathds{1}_{[0,\tau]}(t)+x_d(\theta,\mu)\mathds{1}_{(\tau,+\infty)}(t)$ for all $t>0$. Depending on $\tau$ and the parameters several scenarios are possible. (a) When parameters are set such that $0<c_{u\to d}<c_{d\to u}$, we observe a propagation failure and uniform convergence towards the down state. (b) When parameters are set such that $0 < c_{d\to u} <  c_{u\to d}$ , we observe a propagation through the network  of the up state in the form of a stacked interface (down state / up state / down state)  for large enough $\tau$. (c) When parameters are set such that $c_{d\to u}\leq 0 < c_{u\to d}$, we observe a propagation of the up state in the form of a front through the network for large enough $\tau$.}
 \label{fig:STPbordtau}
\end{figure}

We now turn our attention to the case of an external input which consists in the presentation of a constant input $s_0$ during an initial time window $[0,\tau]$ for $\tau>0$ and then a reset to the down state for later times, that is throughout we set
\bqs
\forall t>0, \quad s_0(t)=s_0\mathds{1}_{[0,\tau]}(t)+x_d(\theta,\mu)\mathds{1}_{(\tau,+\infty)}(t),
\eqs
with $s_0>x_d(\theta,\mu)$.  We now have two parameters $(s_0,\tau)$ which determine the external input. To fix ideas we will only consider the case when the amplitude $s_0$ is set to the up state, that is $s_0=x_u(\theta,\mu)$. Once again, we would like to understand the long time behavior of the solution of \eqref{cont1dN} starting from the initial condition \eqref{IC1dN}. Based on the study of the previous section, we need to consider values of the parameters for which at least $s_0=x_u(\theta,\mu)>s_0^*$, since for $s_0=x_u(\theta,\mu)<s_0^*$ we will necessarily have uniform convergence towards the constant solution $(x_d(\theta,\mu))_{j\in\N}$, and thus no propagation of the up state upward in the network. We have reported in Figure~\ref{fig:STPbordtau} several illustrative space-time plots of the solution for large enough $\tau$. The fate of the solution is directly linked to the sign of the wave speeds $c_{u\to d}$ and $c_{d\to u}$ studied previously. As we have already seen in the case of a constant external input, in order for the up state to have a chance to propagate through the network, one needs to work in the regime where $c_{u\to d}>0$, which is already ensured since we have assumed that $x_u(\theta,\mu)>s_0^*$ implicitly implying that $c_{u\to d}>0$. As a consequence, we need to distinguish between three cases:

\begin{itemize}
\item[(i)] Case $0<c_{u\to d}<c_{d\to u}$. We typically observe a propagation failure, see Figure~\ref{fig:STPbordtau}(a), and the solution uniformly converges towards the down state: that is
\bqs
\underset{j\in\N}{\sup}\left| v_j(t)- x_d(\theta,\mu)\right|\underset{t\rightarrow+\infty}{\longrightarrow}0.
\eqs
Indeed, the interface between the down state and the up state propagates at a larger speed than the interface between up state and the down state. 
\item[(ii)] Case $0 < c_{d\to u} <  c_{u\to d}$. We typically observe the propagation of the up state upward in the network in the form of a \emph{stacked interface} which consists in the concatenation of the two interfaces, one between the down state and the up state and one between the up state and the down state, each such interface propagating at speed $c_{d\to u} $ and  $c_{u\to d}$ respectively. More precisely, for any small $\epsilon>0$, one has
\bqs
\underset{(c_{d\to u}-\epsilon)t \leq j \leq  (c_{u\to d}-\epsilon)t}{\sup}\left| v_j(t)- x_u(\theta,\mu)\right| \underset{t\rightarrow+\infty}{\longrightarrow}0,
\eqs
while
\bqs
\underset{0\leq j \leq (c_{d\to u}-\epsilon)t}{\sup}\left| v_j(t)- x_d(\theta,\mu)\right| \underset{t\rightarrow+\infty}{\longrightarrow}0, \quad \text{ and } \quad \underset{ (c_{u\to d}-\epsilon)t\leq j}{\sup}\left| v_j(t)- x_d(\theta,\mu)\right| \underset{t\rightarrow+\infty}{\longrightarrow}0.
\eqs
This is illustrated in Figure~\ref{fig:STPbordtau}(b).
\item[(iii)] $ c_{d\to u} \leq 0<  c_{u\to d}$. In this case, we observe a propagation of the up state in the form of a front through the network, and the solution locally uniformly converges towards $\boldsymbol{x}^u(x_d(\theta,\mu))$, stationary solution of \eqref{statsolNu} with $s_0=x_d(\theta,\mu)$. We refer to Figure~\ref{fig:STPbordtau}(c) for a visualization of this scenario.
\end{itemize}

\begin{figure}[t!]
\centering
\subfigure[Small $\tau$.]{\includegraphics[width=.32\textwidth]{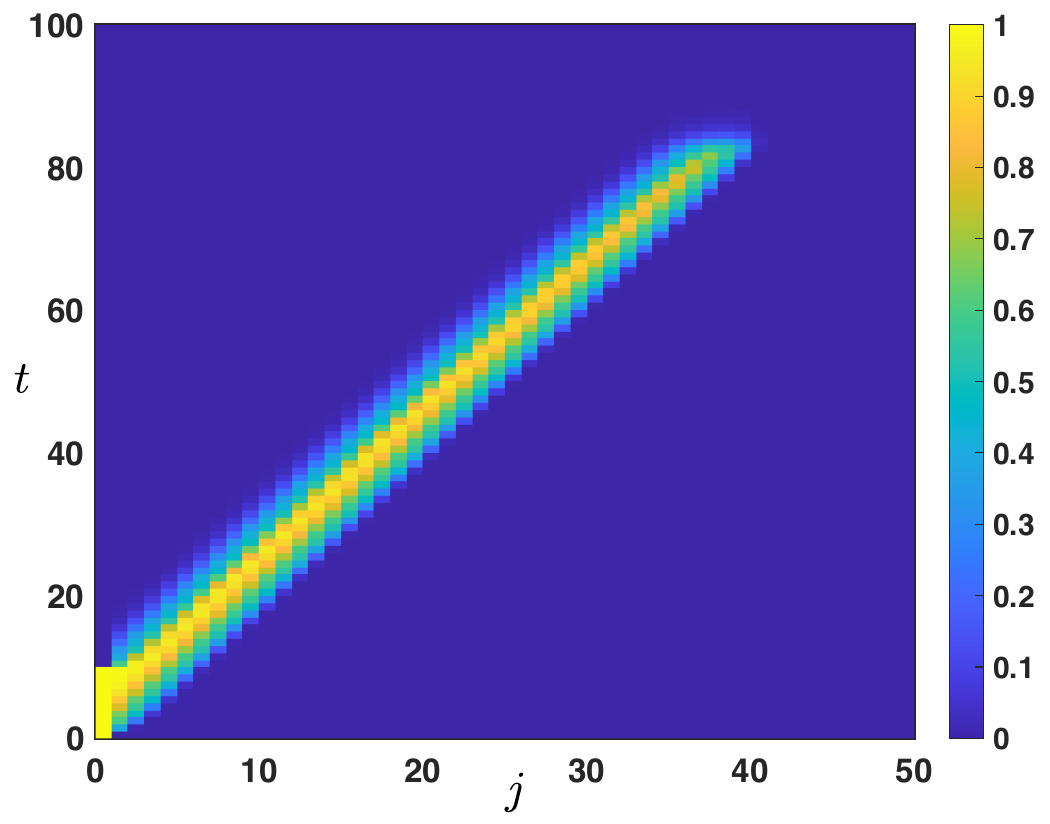}}
\subfigure[Medium $\tau$.]{\includegraphics[width=.32\textwidth]{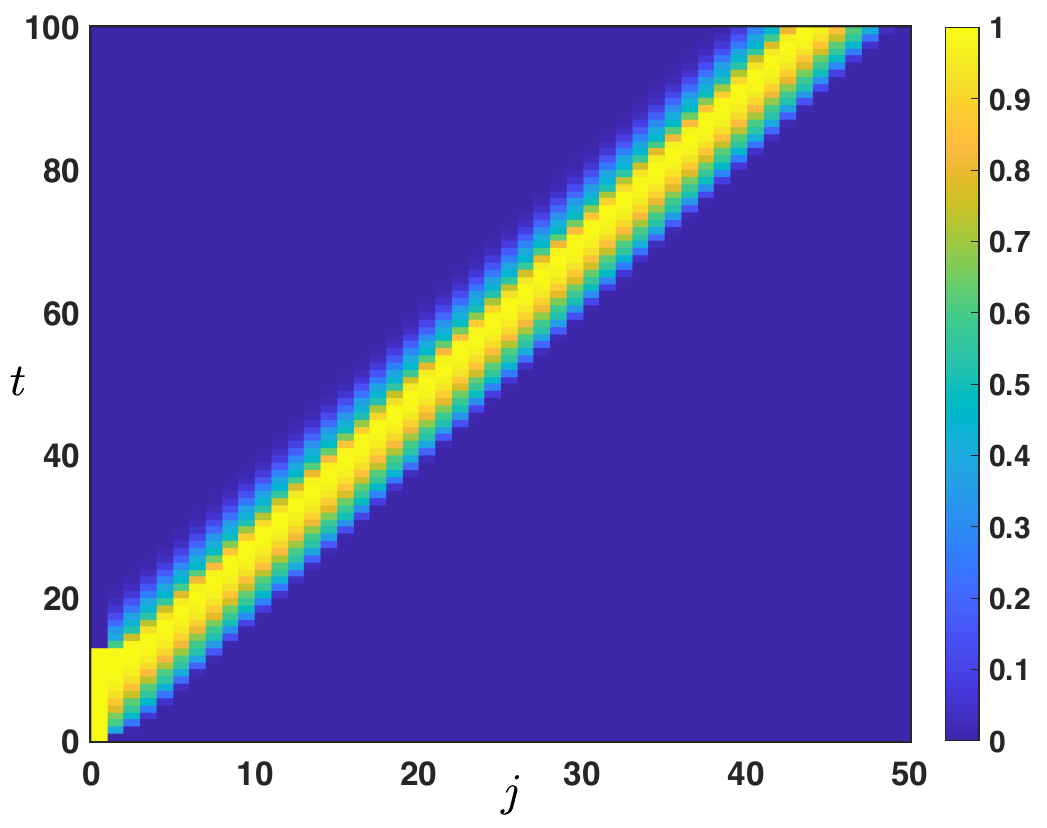}}
\subfigure[Large $\tau$.]{\includegraphics[width=.32\textwidth]{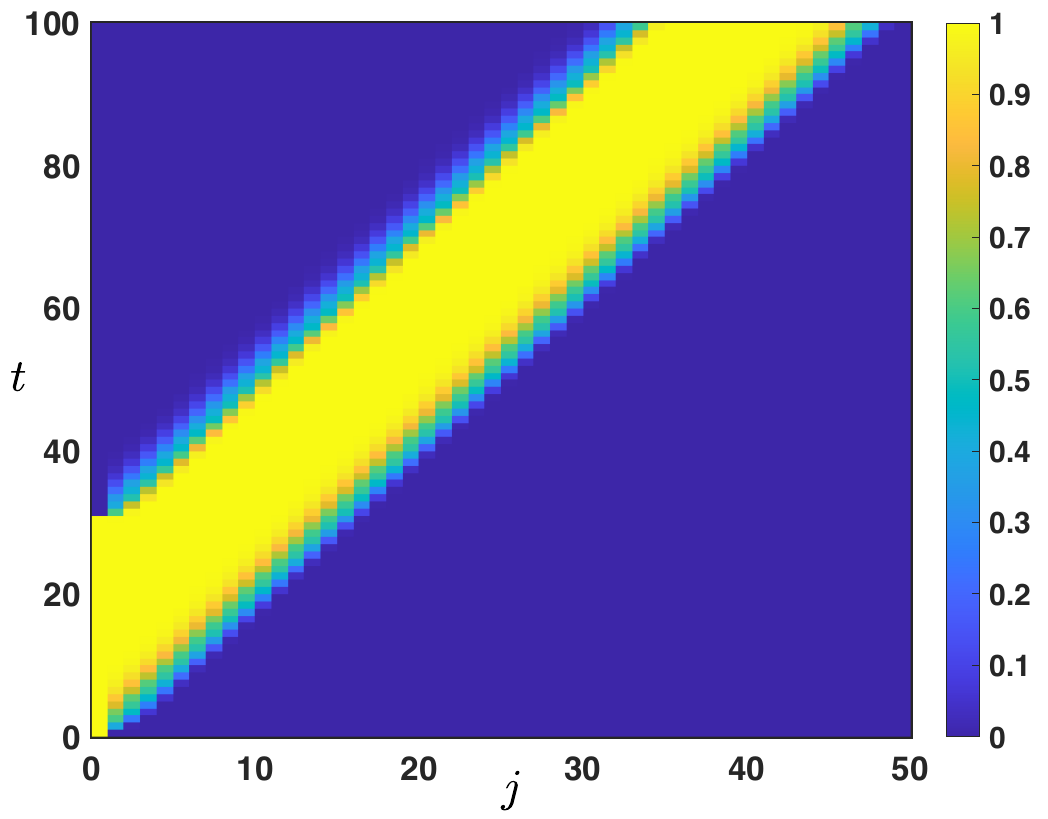}}
\caption{Space-time plot of the solutions of \eqref{cont1dN}-\eqref{IC1dN} when $s_0(t)=x_u(\theta,\mu)\mathds{1}_{[0,\tau]}(t)+x_d(\theta,\mu)\mathds{1}_{(\tau,+\infty)}(t)$ for all $t>0$ when parameters are set such that $0<c_{d\to u}=c_{u\to d}$. For $\tau$ larger than a critical value, we observe a propagation through the network  of the up state in the form of traveling pulse whose width depends monotonically on $\tau$.}
 \label{fig:STPbordtaupulse}
\end{figure}

There exists in fact a fourth situation which does not precisely fall into our above classification. It lies at the boundary between cases (i) and (ii). It is the case where both wave speed $c_{d\to u} $ and  $c_{u\to d}$ are equal and positive, that is $0<c_{d\to u}=c_{u\to d}$. This case happens when $\theta=\frac{1}{2}$ and $q\in[0,1]$ is sufficiently small, see Figure~\ref{fig:speed}. In that case, we typically observe a propagation through the network  of the up state in the form of \emph{traveling pulse}, see Figure~\ref{fig:STPbordtaupulse}(b)-(c), where the width of this traveling pulse depends on $\tau$. Indeed, the larger $\tau$ and the larger is the width, indicating that there is a whole family of such traveling pulses (indexed by the width, and thus by $\tau$). The proper mathematical study of such special solution is beyond the scope of the present study and is left for future work. 

It is possible to further interpret the four states shown in  Figure~\ref{fig:STPbordtau} and Figure~\ref{fig:STPbordtaupulse} from a biological perspective.  Specifically, we may relate the propagation failure to some types of unconscious perception: when the stimulus is weak or below the visual threshold, it may generate activity that propagates through the early stages of the cortex yet without reaching the hierarchical higher areas \cite{dehIgn}. 
On the other hand, when the down state propagates more slowly than the up-state (panel b of Figure~\ref{fig:STPbordtau}) or with the same speed (Figure~\ref{fig:STPbordtaupulse}), we observe the case in which the activity travels through the network followed by the down state propagation, which eventually resolves the activity to the down state. 
To some extent, panel 'c' of Figure~\ref{fig:STPbordtau} represents the only case in which the activity stably propagates through the network in both directions: from lower levels forward due to the presence of the stimulus, and from higher nodes backward, as the top-down influence contributes to maintaining the activity at each node over time (i.e., there is no propagation of the down state). As a final word, we note that our modeling approach has intrinsic limitations (e.g. lack of adaptation and stabilization of the activity over time) and thus, the above biological interpretations need to be taken with caution and in light of our simplifying modeling assumptions.

\begin{figure}[t!]
\centering
\subfigure[$q=0.35$.]{\includegraphics[width=.32\textwidth]{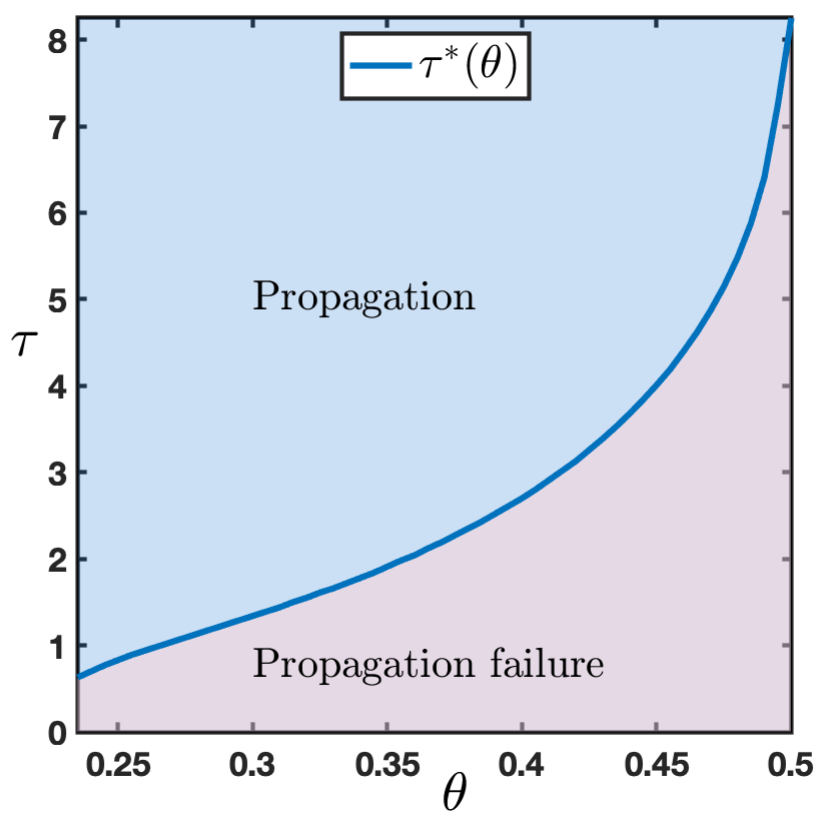}}
\subfigure[$q=0.5$.]{\includegraphics[width=.32\textwidth]{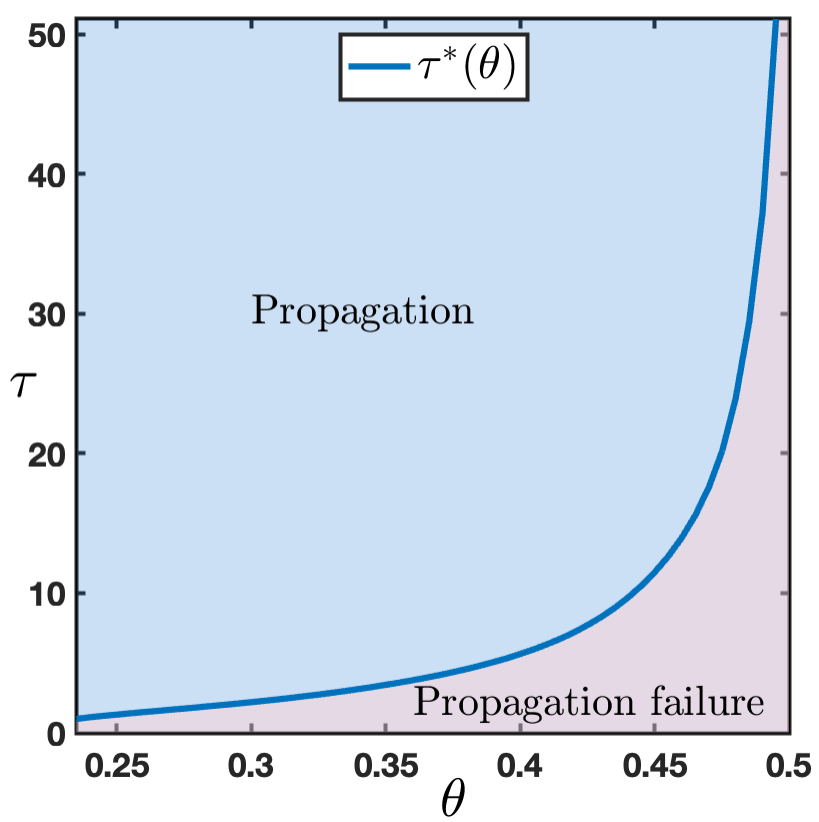}}
\subfigure[$q=0.65$.]{\includegraphics[width=.32\textwidth]{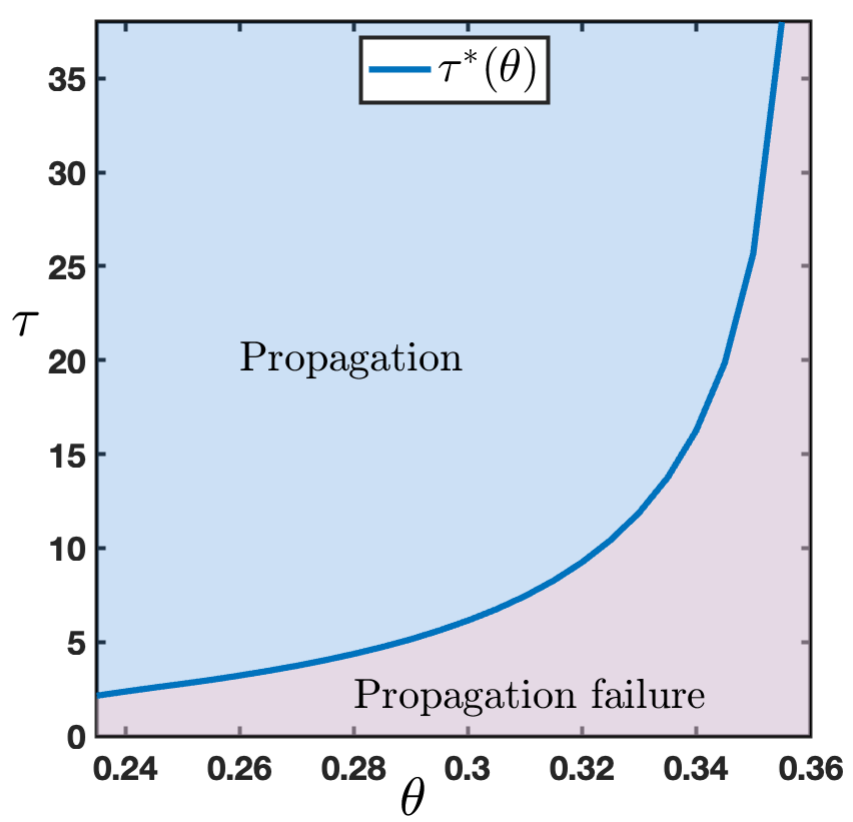}}
\caption{Sharp threshold $\tau^*$ as a function of $\theta$ for several values of $q\in\left\{0.35,0.5,0.65\right\}$ while other values of the parameters are set to $(\mu,p)=(16,0.1)$. The range of $\theta$ is given by those values where $c_{u\to d}>0$ and $c_{d\to u}<c_{u\to d}$, for $q=0.35$ or $q=0.5$ the range is $\theta\in(\theta_*(\mu),1/2)$, while for $q=0.65$ it is given by $\theta\in(\theta_*(\mu),\theta_0)$, here $\theta_0\sim0.355$.}
 \label{fig:tau}
\end{figure}

Now that we have identified the different possible scenarios regarding the long time behavior of the solutions of \eqref{cont1dN}-\eqref{IC1dN} under the assumption that $\tau$ was chosen large enough, we would like to quantify, if it exists, an eventual threshold value $\tau^*$ which is the onset between propagation failure (uniform convergence to the down state) and propagation (either in the form of a stacked interface or a front interface). We numerically computed such a threshold, as illustrated in Figure~\ref{fig:tau}, for several representative values of the parameters. This threshold value $\tau^*$ can be interpreted as the minimal time of presentation of an external input (here constantly fixed to the up state) which is needed to generate a full propagation of neural activity  throughout the network (either as a stacked propagation or front propagation).

\subsection{Wave initiation on a semi-infinite depth network -- The top-down case}

As a complement to the study conducted in the previous section, we also work with a semi-infinite network, but this time we consider
\bqq\label{cont1dNneg}
\forall t>0, \quad \left\{\begin{split}
v_j'(t) & =\mathscr{N}\left(v_{j-1}(t),v_j(t),v_{j+1}(t)\right), \quad j \leq -1, \\
v_0(t)& = s_0(t).
\end{split}
\right.
\eqq
This means that we are interested in top-down propagation. Once again, the network is initially in the down state, that is
\bqq
\label{IC1dNneg}
v_j(0)=u_d(\theta,\mu), \quad j\leq -1,
\eqq
with parameters $(\theta,\mu,p,q)\in\mathcal{P}$. Reproducing the analysis of the bottom-up case of the previous section, we shall only focus on two specific external inputs:
\begin{itemize}
\item constant input of the form $s_0(t)=s_0\geq x_d(\theta,\mu)$ for all $t>0$;
\item flashed input of the form $s_0(t)=s_0\mathds{1}_{[0,\tau]}(t)+x_d(\theta,\mu)\mathds{1}_{(\tau,+\infty)}(t)$ for all $t>0$ with $s_0>x_d(\theta,\mu)$ and $\tau>0$.
\end{itemize}

\subsubsection{Threshold phenomena for constant external input}

We first study the long time behavior of the solutions of \eqref{cont1dNneg}-\eqref{IC1dNneg} with a constant external input given by $s_0(t)=s_0\geq x_d(\theta,\mu)$ for all $t>0$. Once again, our aim is to characterize the strength of the amplitude $s_0$ which can give rise to a full activation of the network to the up state. Similarly as in the bottom-up case, if $s_0\in[x_d(\theta,\mu),x_m(\theta,\mu)]$, then the corresponding solution satisfies $v_j(t)\in[x_d(\theta,\mu),x_m(\theta,\mu)]$ for all $t>0$ and $j\leq-1$, which means that $s_0$ needs to be sufficiently strong, that is $s_0>x_m(\theta,\mu)$ in order for the solution to eventually converge to the up state. 

\begin{figure}[t!]
\centering
\subfigure[Stagnation: $s_0<s_0^*$.]{\includegraphics[width=.35\textwidth]{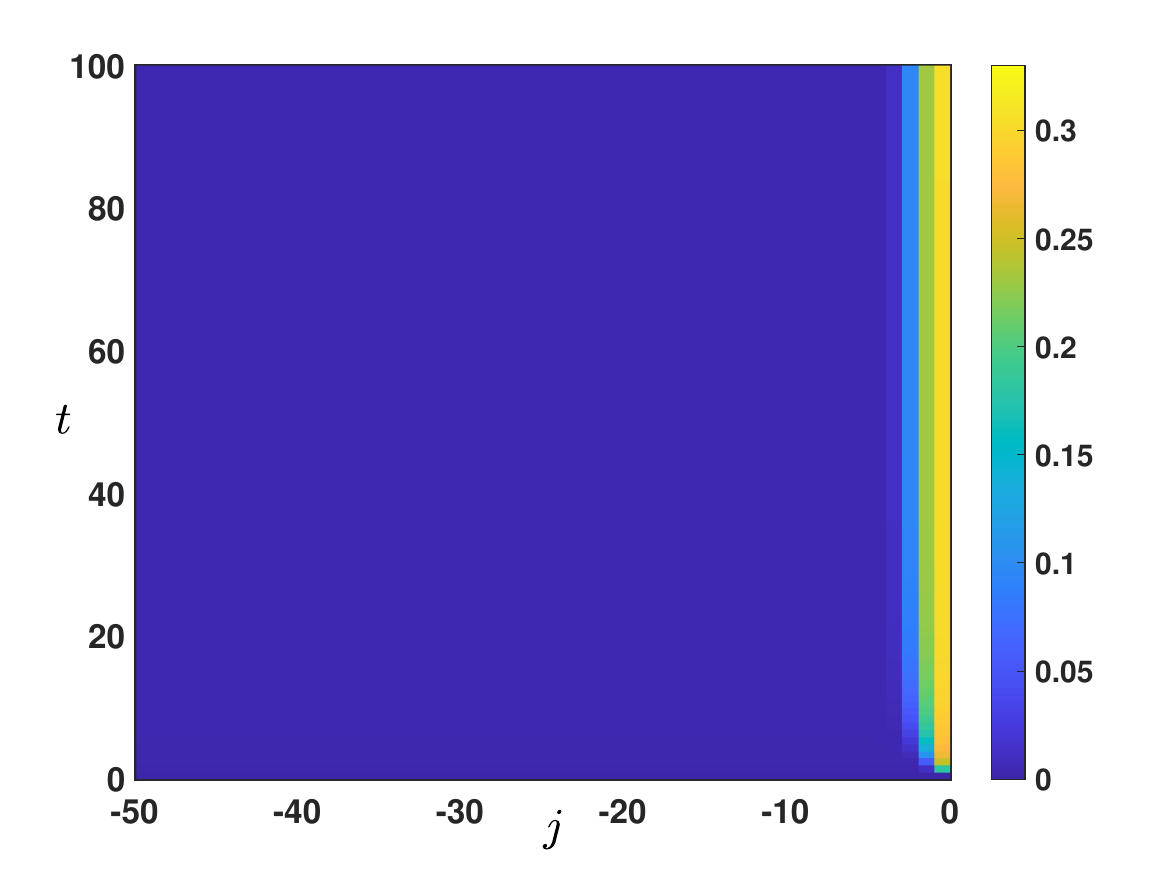}}\hspace{1cm}
\subfigure[Propagation: $s_0>s_0^*$.]{\includegraphics[width=.35\textwidth]{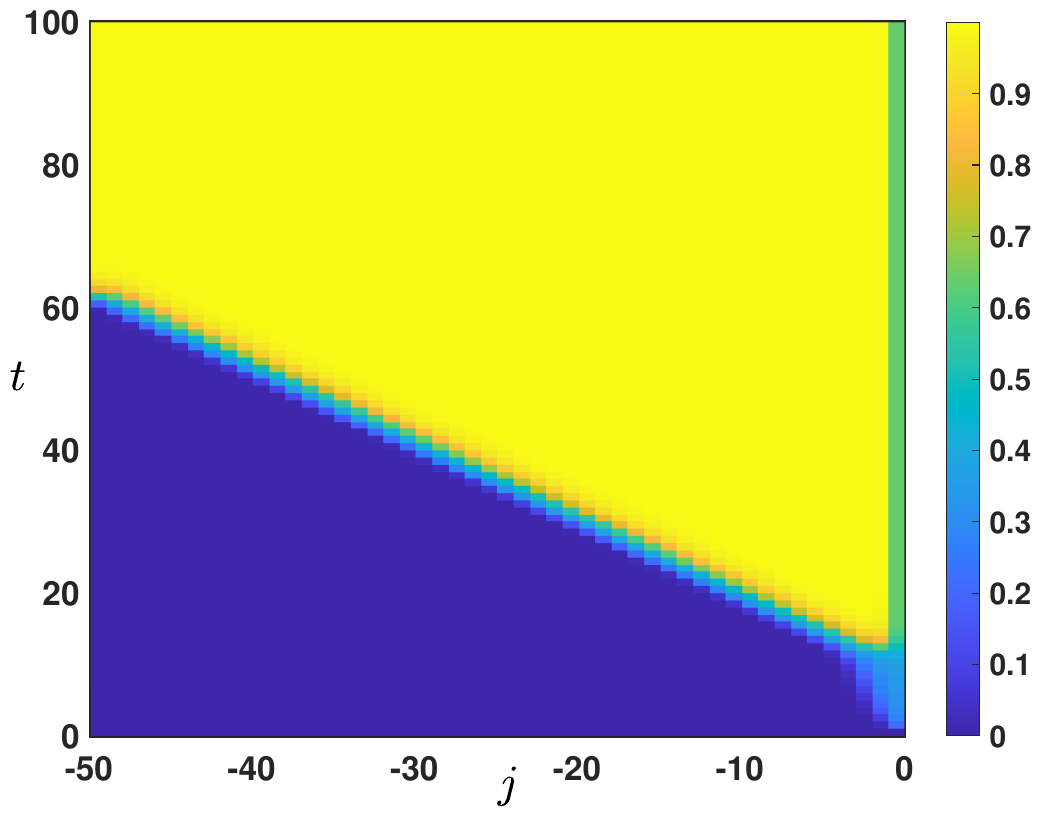}}
\caption{Space-time plot of the solutions of \eqref{cont1dNneg}-\eqref{IC1dNneg} when $s_0(t)=s_0$ for all $t>0$. Depending on the amplitude $s_0 \lessgtr s_0^*$ for some threshold $s_0^*$, the solution can either uniformly converge towards the stationary solution $\boldsymbol{x}^d(s_0)$ (stagnation) or locally uniformly converge towards the stationary solution $\boldsymbol{x}^u(s_0)$ (top-down propagation). }
 \label{fig:STPbordNeg}
\end{figure}

Another important remark to make is the following. From the analysis of the bi-infinite case, we know that downward propagation of the up state to the down state
can only occur with a negative wave speed $c_{d\to u}<0$. As a consequence, for the up state to propagate
downward in the network in the negative semi-infinite case, one needs for $\Lambda\in\mathcal{P}$ to be in a region of parameter space where $c_{d\to u}(\Lambda)<0$

\begin{figure}[t!]
\centering
\subfigure[$\theta=0.35$.]{\includegraphics[width=.32\textwidth]{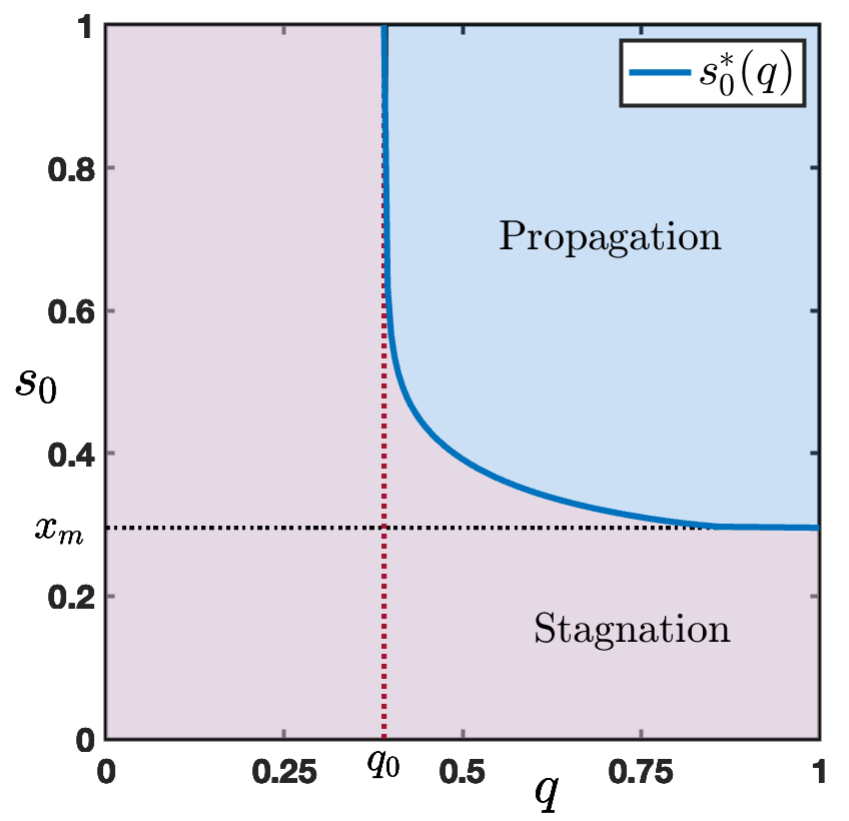}}
\subfigure[$\theta=0.5$.]{\includegraphics[width=.32\textwidth]{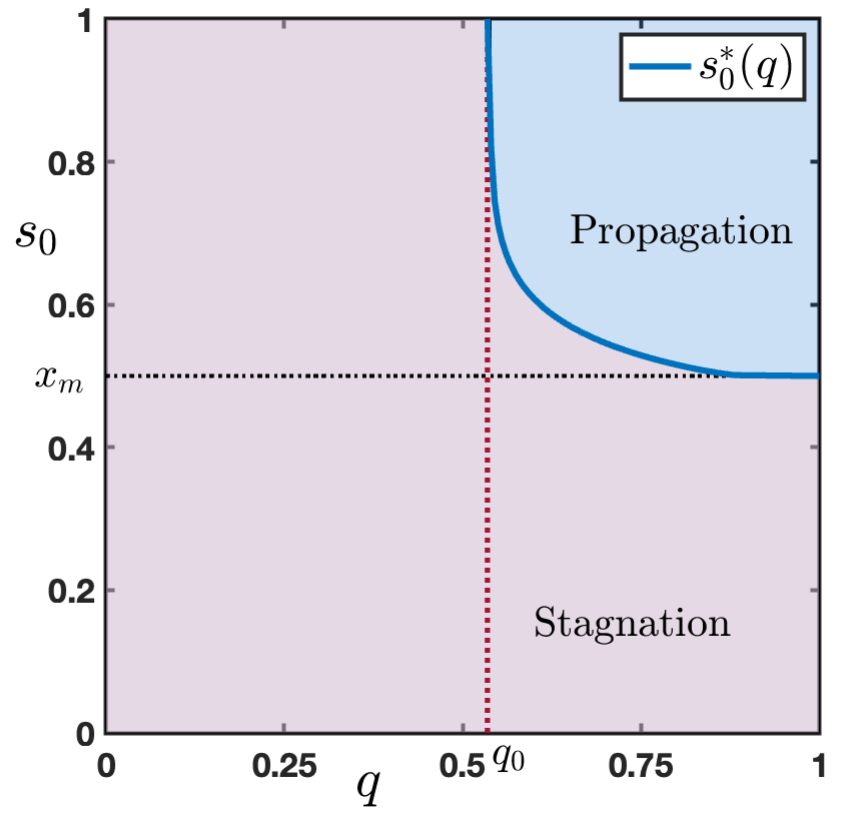}}
\subfigure[$\theta=0.65$.]{\includegraphics[width=.32\textwidth]{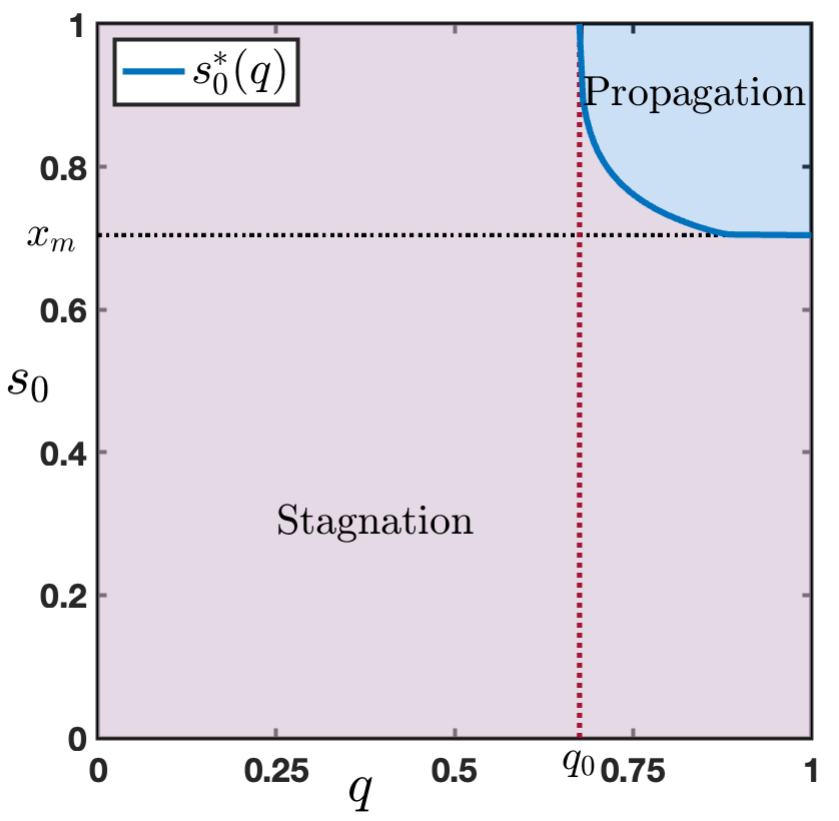}}
\subfigure[$q=0.35$.]{\includegraphics[width=.32\textwidth]{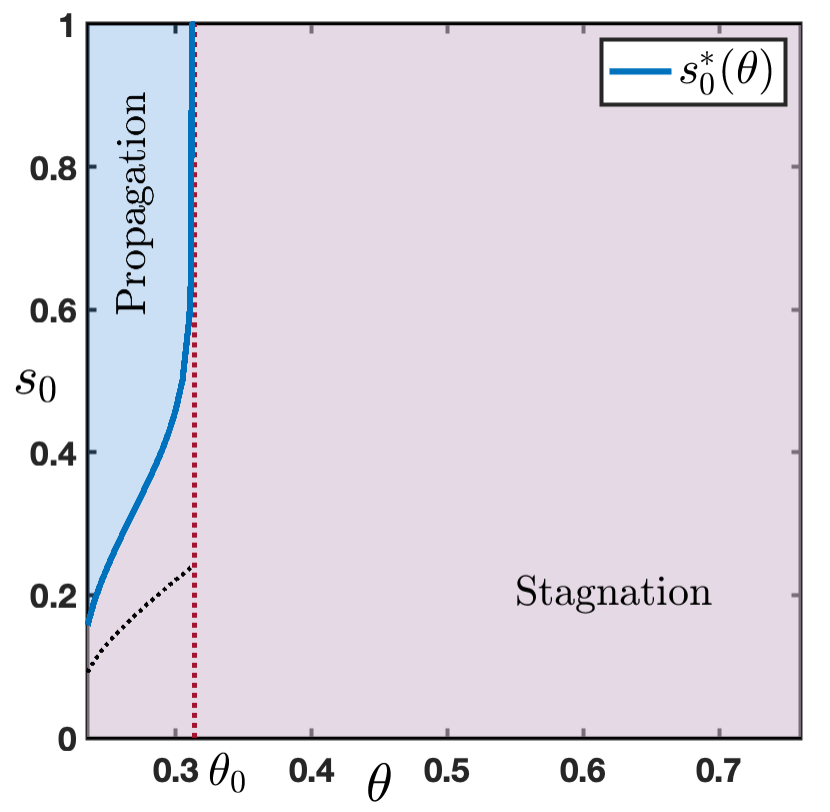}}
\subfigure[$q=0.5$.]{\includegraphics[width=.32\textwidth]{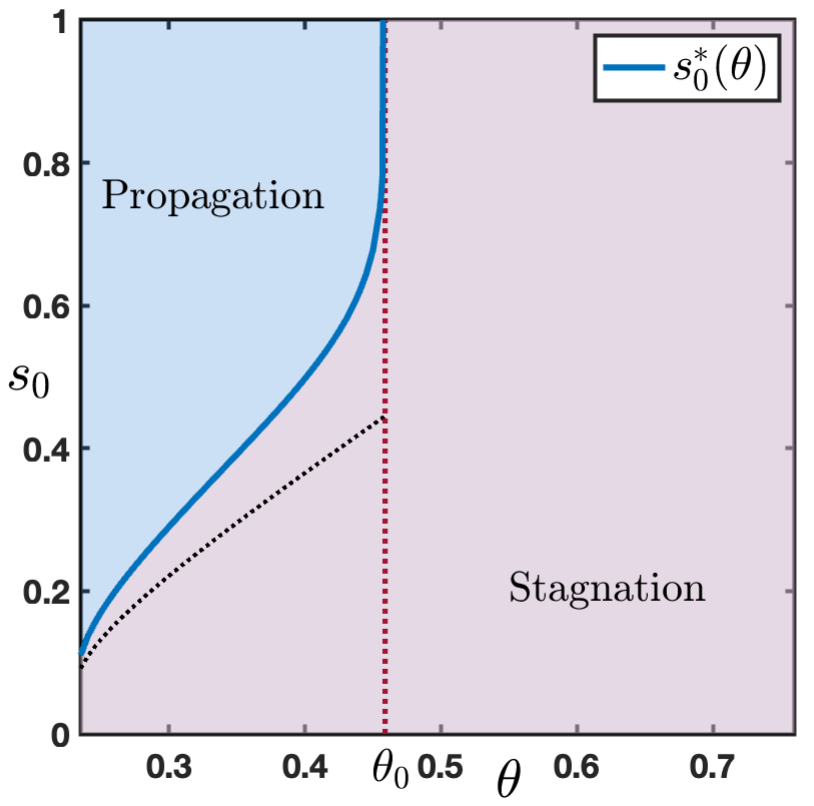}}
\subfigure[$q=0.65$.]{\includegraphics[width=.32\textwidth]{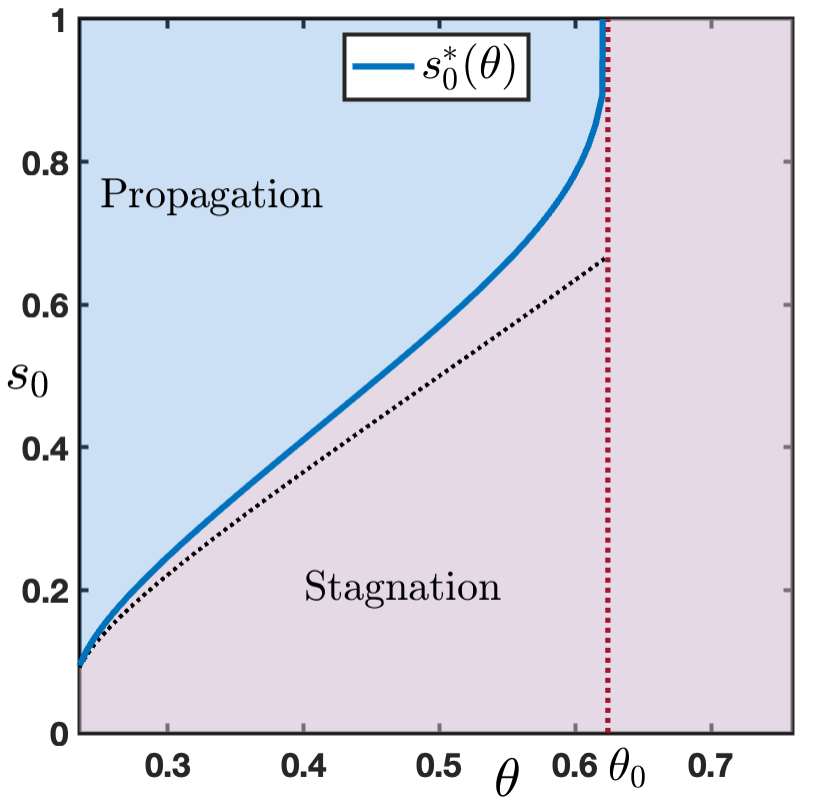}}
\caption{(a)-(b)-(c) Plot of the sharp threshold $s_0^*$ as a function of $q$ for several values of $\theta\in\left\{0.35,0.5,0.65\right\}$. (d)-(e)-(f) Plot of the sharp threshold $s_0^*$ as a function of $\theta$ for several values of $q\in\left\{0.35,0.5,0.65\right\}$. In the blue region there is propagation while in the violet region there is stagnation. Other values of the parameters are set to $(\mu,p)=(16,0.1)$.}
 \label{fig:s0Negqvary}
\end{figure}

To define the equivalent notion of stagnation and propagation in the top-down case, we introduce special stationary solutions of \eqref{cont1dNneg}. First, we denote by $\boldsymbol{y}^d(x_0)=(y_j^d)_{j\in\N_-}$ a sequence solution of
\bqq
\left\{
\begin{split}
0&= \mathscr{N}\left(y_{j-1},y_j,y_{j+1}\right), \quad j \leq -1, \\
y_0&=s_0,\\
y_j& \underset{j\rightarrow+\infty}{\longrightarrow} x_d(\theta,\mu),
\end{split}
\right.
\label{statsolNegd}
\eqq
and $\boldsymbol{y}^u(x_0)=(y_j^u)_{j\in\N_-}$ a sequence solution of  
\bqq
\left\{
\begin{split}
0&= \mathscr{N}\left(y_{j-1},y_j,y_{j+1}\right), \quad j \leq -1, \\
y_0&=s_0,\\
y_j& \underset{j\rightarrow+\infty}{\longrightarrow} x_u(\theta,\mu).
\end{split}
\right.
\label{statsolNegu}
\eqq
For a given $s_0\geq x_d(\theta,\mu)$, we say that there is \emph{propagation} if $\boldsymbol{v}=(v_j)_{j\in\N_-}$ the solution of \eqref{statsolNegd} satisfies 
\bqs
\boldsymbol{v}(t) \underset{t\rightarrow+\infty}{\longrightarrow} \boldsymbol{y}^u(s_0) \text{ locally uniformly on } \N_-,
\eqs
while there is \emph{stagnation} if the corresponding solution satisfies
\bqs
\boldsymbol{v}(t) \underset{t\rightarrow+\infty}{\longrightarrow} \boldsymbol{y}^d(s_0) \text{ uniformly on } \N_-.
\eqs
As we have already discussed, for all $s_0\in[x_d(\theta,\mu),x_m(\theta,\mu)]$ stagnation occurs. Our numerical investigations, see Figure~\ref{fig:STPbordNeg} and  Figure~\ref{fig:s0Negqvary}, show the existence of a sharp threshold between stagnation and propagation as $s_0$ is further increased. More precisely, there exists $s_0^*>x_m(\theta,\mu)$ such that the following dichotomy holds.  For all $s_0\in[x_d(\theta,\mu),s_0^*)$ there is stagnation, while for all $s_0>s_0^*$ there is propagation.

In Figure~\ref{fig:s0Negqvary}, we plot the value of the sharp threshold $s_0^*$ as a function of $q$ for several values of $\theta\in\left\{0.35,0.5,0.65\right\}$ (see panels (a)-(b)-(c)). The sharp threshold is well defined for all values of $q\in(q_0,1]$ where we denote by $q_0\in(0,1)$ the largest value of $q\in[0,1]$ for which $c_{d\to u}=0$, that is $q_0=\max\left\{q\in[0,1]~|~c_{d\to u}=0\right\}$. As anticipated, we note that $s_0^*(q)\geq x_m(\theta,\mu)$ for all $q\in(q_0,1]$ and $s_0^*(q)\rightarrow+\infty$ as $q\to q_0^+$.  When the sharp threshold $s_0^*$ is a function of $\theta$, we observe similar behaviors as reported in panels (d)-(e)-(f) of Figure~\ref{fig:s0Negqvary} corresponding to several values of $q\in\left\{0.35,0.5,0.65\right\}$. Here, the sharp threshold is well defined for all values of $\theta\in(\theta_*(\mu),\theta_0)$ where we denote by $\theta_0\in(\theta_*(\mu),\theta^*(\mu))$ the smallest value of $\theta$ for which $c_{d\to u}=0$, that is $\theta_0=\min\left\{ \theta\in(\theta_*(\mu),\theta^*(\mu))~|~ c_{d\to u}=0\right\}$. Once again, we remark that $s_0^*(\theta)\geq x_m(\theta,\mu)$ for all $\theta\in(\theta_*(\mu),\theta_0)$ and $s_0^*(\theta)\rightarrow+\infty$ as $\theta\rightarrow \theta_0^-$.

\begin{figure}[t!]
\centering
\includegraphics[width=.42\textwidth]{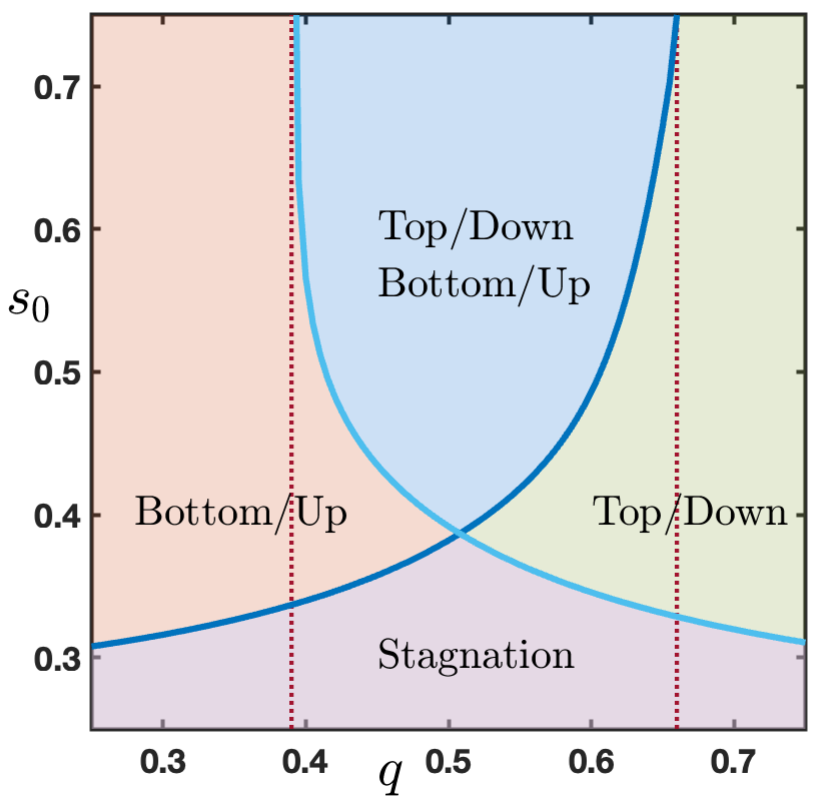}
\caption{Merging of Figure~\ref{fig:s0qvary}(a) associated to system \eqref{cont1dN} and Figure~\ref{fig:s0Negqvary}(a) associated to \eqref{cont1dNneg} for $\theta=0.35$. The violet region corresponds to stagnation for both systems. In the blue region, there is both bottom-up propagation for system \eqref{cont1dN} and top-down propagation for system \eqref{cont1dNneg}. In the orange (resp. green) region, there is bottom-up (resp. top-down) propagation for system \eqref{cont1dN} (resp. \eqref{cont1dNneg}) while stagnation for system \eqref{cont1dNneg} (resp. \eqref{cont1dN}).}
\label{fig:combineds0}
\end{figure}

Comparing Figure~\ref{fig:s0Negqvary} panels (a)-(b)-(c) with the bottom-up case in Figure~\ref{fig:s0qvary}, we observe that for values of $\theta$ slightly below $1/2$, there exists a non empty interval of values of $q$ for which we have at the same time both bottom-up propagation for system \eqref{cont1dN} and top-down propagation for system \eqref{cont1dNneg} with constant external input $s_0$. This is summarized in Figure~\ref{fig:combineds0}. On the other hand, for $\theta$ larger than $1/2$, there is no common range for $q$ where bottom-up propagation for system \eqref{cont1dN} and top-down propagation for system \eqref{cont1dNneg} can occur at the same time. We postulate that a \emph{normal} neural working regime is when both systems admit propagation (bottom-up for \eqref{cont1dN} and top-down for \eqref{cont1dNneg}). This allows to identify ranges of values for the parameters $q$ and $\theta$ where this regime exists. Roughly, it corresponds to $\theta\in(\theta_*(\mu),1/2)$ and $q\in(q_0^{u\to d},q_0^{d\to u})$ where we have denoted by $q_0^{u\to d}$ and $q_0^{d\to u}$ the two threshold values $q_0$ of $q$ for which $s_0^*(q)$ is well defined (see Figure~\ref{fig:s0qvary}(a) and Figure~\ref{fig:s0Negqvary}(a)).

\subsubsection{Threshold phenomena for flashed external input}

Finally, we turn our attention to the case of an external input which consists in the presentation of the up state $x_u(\theta,\mu)$ during an initial time window $[0,\tau]$ for $\tau>0$ and then a reset to the down state $x_d(\theta,\mu)$ for later times. That is, we work with an external input of the form
\bqs
\forall t>0, \quad s_0(t)=x_u(\theta,\mu)\mathds{1}_{[0,\tau]}(t)+x_d(\theta,\mu)\mathds{1}_{(\tau,+\infty)}(t).
\eqs
Once again, we would like to understand the long time behavior of the solution of \eqref{cont1dNneg} starting from the initial condition \eqref{IC1dNneg}. Based on the study of the previous section, we need to consider values of the parameters for which at least $x_u(\theta,\mu)>s_0^*$, since for $x_u(\theta,\mu)<s_0^*$ we will necessarily have uniform convergence towards the constant solution $(x_d(\theta,\mu))_{j\in\N_-}$, and thus no propagation of the up state downward in the network. 

\begin{figure}[t!]
\centering
\subfigure[Propagation failure.]{\includegraphics[width=.32\textwidth]{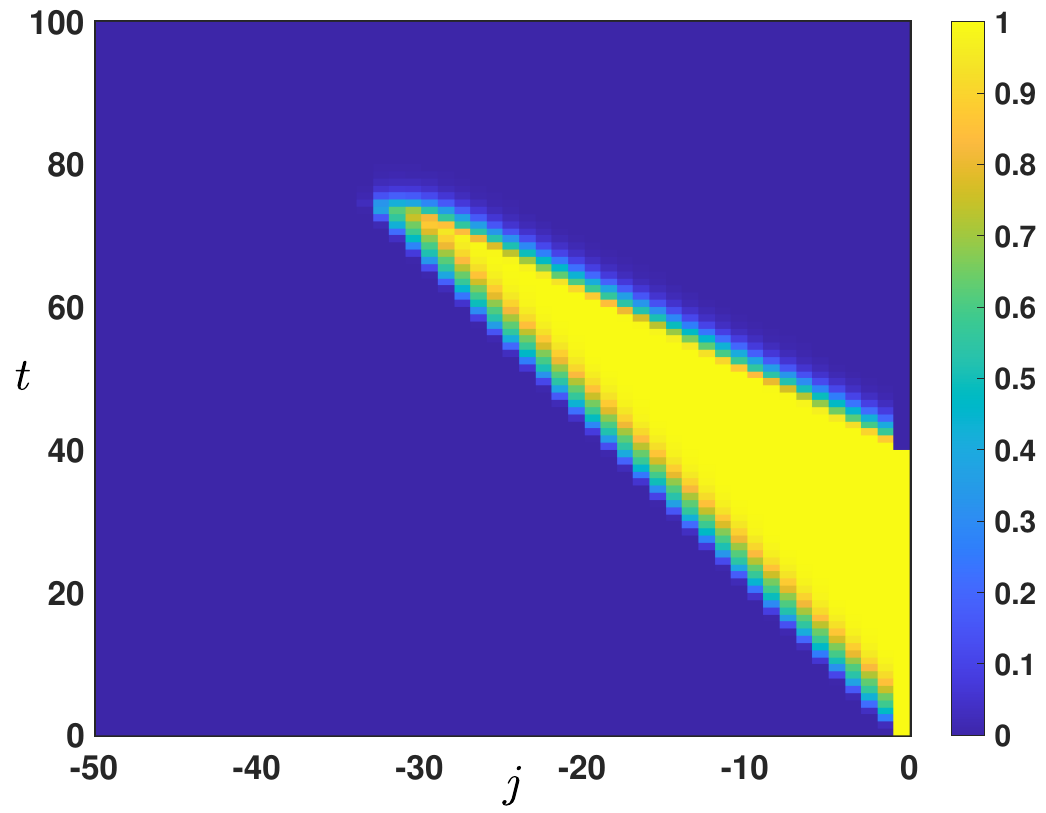}}
\subfigure[Stacked propagation.]{\includegraphics[width=.32\textwidth]{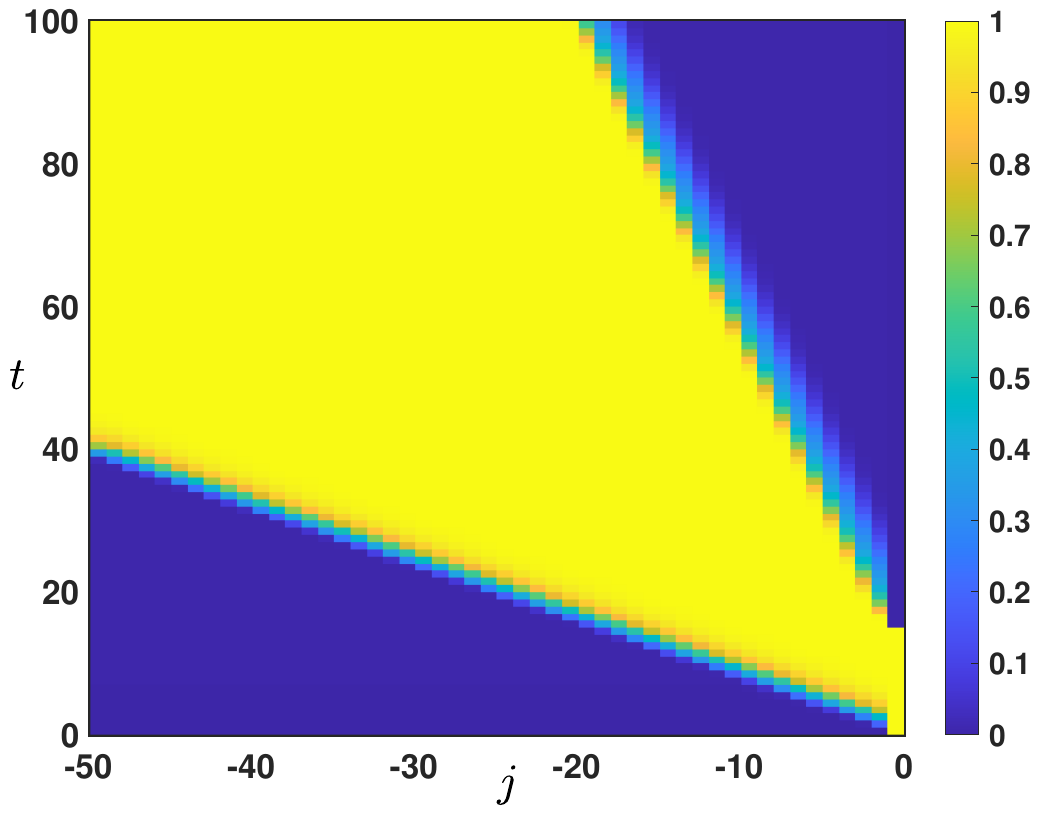}}
\subfigure[Front propagation.]{\includegraphics[width=.32\textwidth]{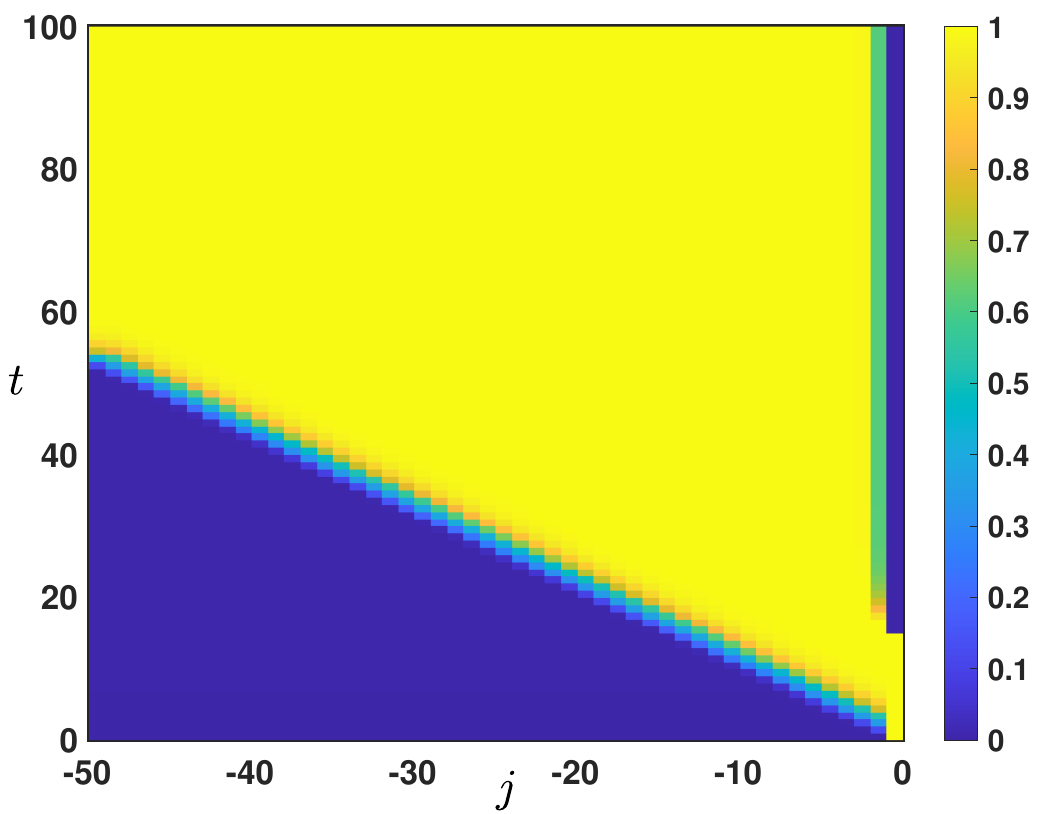}}
\caption{Space-time plot of the solutions of \eqref{cont1dNneg}-\eqref{IC1dNneg} when $s_0(t)=x_u(\theta,\mu)\mathds{1}_{[0,\tau]}(t)+x_d(\theta,\mu)\mathds{1}_{(\tau,+\infty)}(t)$ for all $t>0$. Depending on $\tau$ and the parameters several scenarios are possible. (a) When parameters are set such that $c_{u\to d}<c_{d\to u}<0$, we observe a propagation failure and uniform convergence towards the down state. (b) When parameters are set such that $c_{d\to u} <  c_{u\to d}<0$ , we observe a propagation through the network  of the up state in the form of a stacked interface (down state / up state / down state)  for large enough $\tau$. (c) When parameters are set such that $c_{d\to u}< 0 \leq c_{u\to d}$, we observe a propagation of the up state in the form of a front through the network for large enough $\tau$.}
 \label{fig:STPbordtauNeg}
\end{figure}

We have reported in Figure~\ref{fig:STPbordtauNeg} several illustrative space-time plots of the solution for large enough $\tau$. The fate of the solution is directly linked to the sign of the wave speeds $c_{d\to u}$ and $c_{u\to d}$. As we have already seen in the case of a constant external input, in order for the up state to have a chance to propagate through the network, one needs to work in the regime where $c_{d\to u}<0$, which is already ensured since we have assumed that $x_u(\theta,\mu)>s_0^*$ implicitly implying that $c_{d\to u}<0$. As a consequence, we need to distinguish between three cases:
\begin{itemize}
\item[(i)] Case $c_{u\to d}<c_{d\to u}<0$. We typically observe a propagation failure, see Figure~\ref{fig:STPbordtauNeg}(a), and the solution uniformly converges towards the down state: that is
\bqs
\underset{j\in\N_-}{\sup}\left| v_j(t)- x_d(\theta,\mu)\right|\underset{t\rightarrow+\infty}{\longrightarrow}0.
\eqs
Indeed, the interface between the up state and the down state propagates at a larger speed than the interface between down state and the up state. 
\item[(ii)] Case $ c_{d\to u} <  c_{u\to d}<0$. We typically observe the propagation of the up state downward in the network in the form of a \emph{stacked interface} which consists in the concatenation of the two interfaces, one between the up state and the down state and one between the down state and the up state, each such interface propagating at speed $c_{u\to d} $ and  $c_{d\to u}$ respectively. More precisely, for any small $\epsilon>0$, one has
\bqs
\underset{(c_{d\to u}-\epsilon)t \leq j \leq  (c_{u\to d}-\epsilon)t}{\sup}\left| v_j(t)- x_u(\theta,\mu)\right| \underset{t\rightarrow+\infty}{\longrightarrow}0,
\eqs
while
\bqs
\underset{ j \leq (c_{d\to u}-\epsilon)t}{\sup}\left| v_j(t)- x_d(\theta,\mu)\right| \underset{t\rightarrow+\infty}{\longrightarrow}0, \quad \text{ and } \quad \underset{ (c_{u\to d}-\epsilon)t\leq j \leq 0}{\sup}\left| v_j(t)- x_d(\theta,\mu)\right| \underset{t\rightarrow+\infty}{\longrightarrow}0.
\eqs
This is illustrated in Figure~\ref{fig:STPbordtauNeg}(b).
\item[(iii)] $ c_{d\to u} < 0\leq  c_{u\to d}$. In this case, observe a propagation of the up state in the form of a front through the network, and the solution locally uniformly converges towards $\boldsymbol{y}^u(x_d(\theta,\mu))$, stationary solution of \eqref{statsolNegu} with $s_0=x_d(\theta,\mu)$. We refer to Figure~\ref{fig:STPbordtauNeg}(c) for a visualization of this scenario.
\end{itemize}
Once again, there is the critical case where $c_{d\to u} =  c_{u\to d}<0$, which only occurs whenever $\theta=\frac{1}{2}$ and $q\in[0,1]$ is sufficiently large, see Figure~\ref{fig:speed}(e) and Figure~\ref{fig:signspeedtq}(c). In that case, we also typically observe a propagation through the network of the up state in the form of a traveling pulse, as for the bottom-up case, where the width of this traveling pulse depends on the presentation time $\tau$ with similar properties. 

\begin{figure}[t!]
\centering
\subfigure[$q=0.35$.]{\includegraphics[width=.32\textwidth]{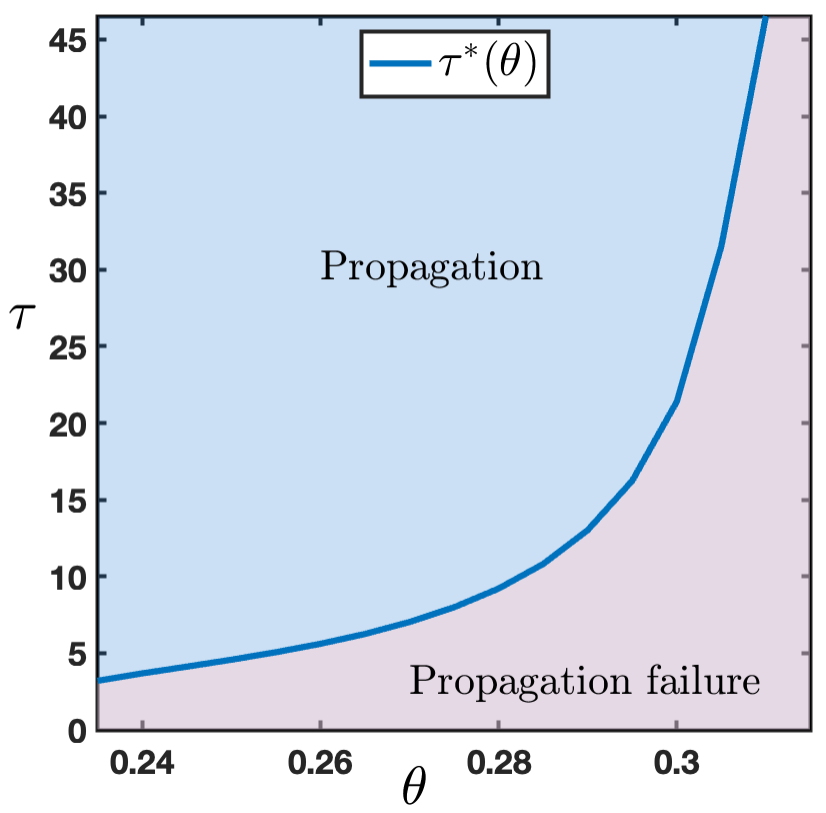}}
\subfigure[$q=0.5$.]{\includegraphics[width=.32\textwidth]{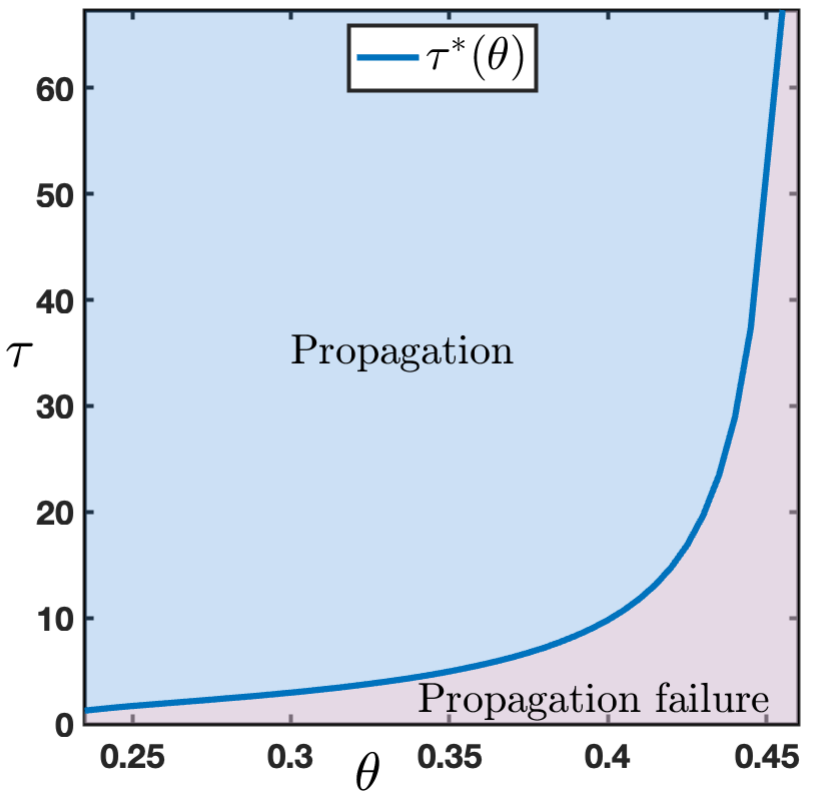}}
\subfigure[$q=0.65$.]{\includegraphics[width=.32\textwidth]{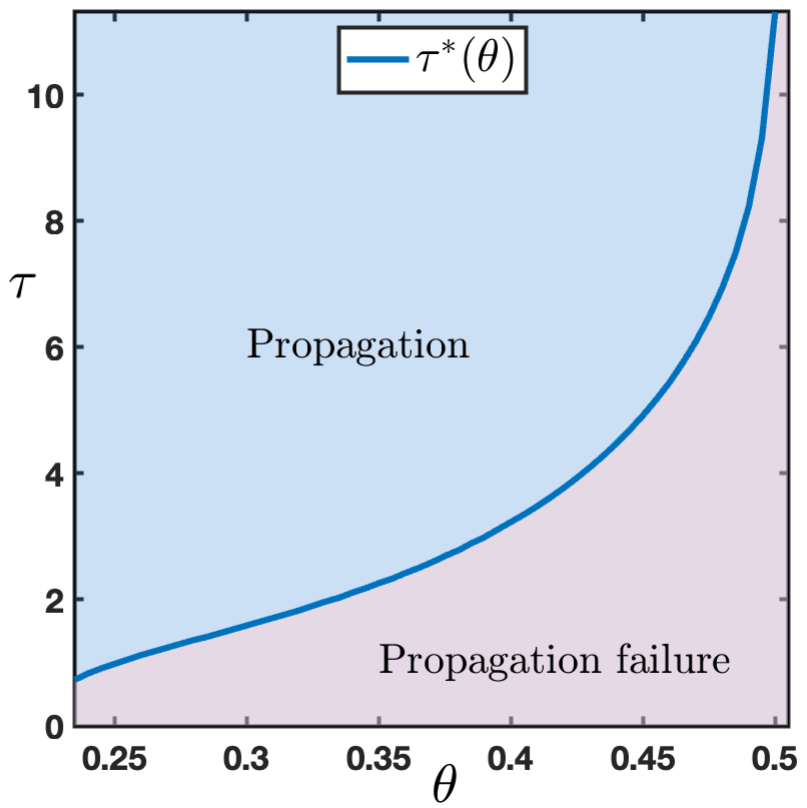}}
\caption{Sharp threshold $\tau^*$ as a function of $\theta$ for several values of $q\in\left\{0.35,0.5,0.65\right\}$ while other values of the parameters are set to $(\mu,p)=(16,0.1)$. The range of $\theta$ is given by those values where $c_{d\to u}<0$ and $c_{d\to u}<c_{u\to d}$, for $q=0.35$ or $q=0.5$ the range is $\theta\in(\theta_*(\mu),\theta_0)$ with $\theta_0$ given in Figure~\ref{fig:s0Negqvary}(d)-(e), while for $q=0.65$ it is given by $\theta\in(\theta_*(\mu),1/2)$.}
 \label{fig:tauNeg}
\end{figure}

Finally, we can quantify, if it exists,
an eventual threshold value $\tau^*$ which is defined as the onset between propagation failure (uniform convergence to the down state) and propagation (either in the form of a stacked interface or a front interface). We numerically
computed such a threshold, as illustrated in Figure~\ref{fig:tauNeg}, for several representative values of the parameters. Results of Figure~\ref{fig:tauNeg} are reversed from the results of Figure~\ref{fig:tau} of the bottom-up case. More precisely, for small value of $q=0.35$, we remark that the critical threshold $\tau^*$ is smaller in the bottom-up case than in the top-down case, while for large value of $q=0.65$, we observe the opposite. Interestingly, for $q=0.5$, in both cases, the critical threshold $\tau^*$ takes  comparable values. These results can naturally be interpreted as follows. Since $q$ represents the strength of feedback connections it somehow measures the capacity of the network to propagate the up-state downwards into the network. As a consequence, large values of $q$ should facilitate the propagation in the top-down case and penalize propagation in the bottom-up case with respectively small and large values for the critical threshold $\tau^*$, which is precisely what is reported in Figure~\ref{fig:tauNeg}(c) and Figure~\ref{fig:tau}(c). 

\subsection{Biological interpretations}

Figure~\ref{fig:combineds0} summarizes the behavior of the network as a function of the parameter $q$, which quantifies the amount of top-down influence, and the strength of the stimulus $s_{0}$. Bearing in mind the limitations of our relatively simple cortical model, it is possible, however, to interpret our results from a broader perspective within the computational psychiatry domain \cite{ComPsyRev}. As mentioned above, previous work theorized the lack of propagation to higher regions in the case of unconscious perception \cite{dehIgn}, in which a weak perceptual stimulus would generate some activity in the lower regions, which fail to propagate to higher areas (i.e., Figure 9). Similarly, a weak endogenous stimulus would fail to propagate top-down and exert its influence on lower areas. On the other hand, stronger stimuli generate activity that travels through the network, either bottom-up (for low values of $q$, left side of the picture in orange) or top-down (higher $q$ values, right side in green), dominating the activity of the network. From a clinical perspective, both cases can be related to dysfunctional perceptions, as the network relies exclusively either on sensory information or on endogenous activity. Recently, some authors proposed that these two perceptual extremes characterize the schizophrenia-autism continuum \cite{tarasi}: individuals within the autistic spectrum tend to rely more on sensory information to drive their behaviors, \cite{autPC01}, whereas hallucinatory percepts in schizophrenia patients could be explained by uncontrolled and dominating top-down processes \cite{autPC02}. Remarkably, recent experimental work confirmed the presence of stronger top-down waves in the schizophrenia population compared to the healthy control \cite{alamia}, interpreting the traveling waves as signatures of predictive coding processes.

\section{Discussion}

\subsection{Summary of contributions}

We proposed a mathematical framework to explore the properties and stability of a class of continuous in time nonlinear neural network models comprising a hierarchy of processing areas (or "layers"), mutually connected according to the principles of predictive coding. More specifically, when considered on a bi-infinite domain, we systematically characterized propagation properties (bottom-up and top-down) by numerically computing the speed of an associated traveling wave. We even showed that for some non empty regions in parameter space, this wave speed identically vanishes indicating the presence of propagation failure phenomena related to the discrete nature of the problem under study. We also explored propagation properties for semi-infinite depth networks and the influence of time varying external inputs presented at the first layer of the network. We could demonstrate the existence of some threshold phenomena quantifying the \emph{strength} of the stimulus needed to generate a propagation either (upward or downward) in the network. We also explored the effect of simple time varying external input in the form of a flashed signal and quantified the minimal time of presentation required in order to generate a full propagation within the network. Finally, we could interpret our results from a broader perspective
within the computational psychiatry domain.

\subsection{Perspectives}

The generality and applicability of the work reported here, for both neuroscientific explanations and machine learning applications, is limited by several simplifying assumptions that could be relaxed or explored in future studies.

\subsubsection{Multi-population models}

A key simplifying assumption that we made was that each neuron in a layer is only connected to the corresponding neuron in the immediately preceding and following layer, with unit weight in each direction. This assumption was mathematically translated by the identity $\W^f=\W^b=\mathbf{I}_d$. Except in the linear case (see \cite{FFVR23}), going beyond this simplifying hypothesis is actually very challenging due to the nonlinear nature of the problem. What can be relaxed easily is the constraint on the normalization of the weights to be equal to one. If one considers, for instance, that both $\W^f$ and $\W^b$ are diagonal matrices, then model \eqref{model} still totally decouples, and its analysis boils down to the study of equations of the form
\bqs
v_j'(t)=\beta \left(\omega_fS(v_{j-1}(t))-v_j(t)\right) +\alpha S'(v_j(t))\omega_b \left( v_{j-1}(t)-\omega_bS( v_j(t))\right)+ \lambda  \left(\omega_b S(v_{j+1}(t))-v_j(t)\right),
\eqs
for two real weights $(\omega_f,\omega_b)\in\R^2$. When $\omega_f=\omega_b=1$, we recover the model that we have studied. Since we recover a scalar equation, the strategy that we have developed here can be extended to that case at the expense of dealing with two extra parameters $(\omega_f,\omega_b)$. An interesting avenue of future research will be to handle the case of two populations, that is $d=2$, which is already relevant when modeling left and right hemispheres of the visual pathway \cite{schwenk} for example. In this two-dimensional setting, one can further work under the assumption that both $\W^f$ and $\W^b$ are symmetric matrices; this already introduces six new parameters to the problem.

\subsubsection{Predictive coding with spike frequency adaptation}

In terms of dynamics, one of the limitations of the predictive coding model \eqref{model} stands from its intrinsic bistable nature: competition of stability between two stable steady states: a down-state and an up-state. As was shown in our numerous numerical simulations, once the network has stabilized to the up-state, it remains indefinitely in this high activity regime, which is not biologically plausible. Indeed, regimes of high neural activity can only be transient in actual biological cortical networks, and several regulatory mechanisms such as synaptic depression (which refers to the process by which resources are depleted in the presynaptic bouton \cite{bart,tsodyks,york}), or spike frequency adaptation (which is a cellular process that brings neurons back to their resting voltage after periods of high activity \cite{benda}), have been proposed and modeled in neural networks \cite{EFK,pinto} to account for such negative feedback mechanisms. Focusing on spike frequency adaptation and specializing to the one-population model \eqref{cont1d}, one can incorporate a linear adaptation variable, denoted by $a_j(t)$, 
\bqs
\left\{
\begin{split}
v_j'(t)&=(1-p-q) \left(S(v_{j-1}(t))-v_j(t)\right) +p S'(v_j(t)) \left( v_{j-1}(t)-S( v_j(t))\right)+ q  \left(S(v_{j+1}(t))-v_j(t)\right)-\gamma a_j(t),\\
 a_j'(t)&= \frac{ v_j(t)-a_j(t)}{\tau},
\end{split}
\right.
\eqs
with $\gamma,\tau>0$ two positive parameters which respectively model the strength of the adaptation mechanism and the timescale on which adaptation occurs. Usually, spike frequency adaptation occurs on a much longer timescale than the internal neural dynamics such that one considers only regimes where $\tau^{-1}\ll1$ so that the above system inherits a slow/fast dynamical structure \cite{FS}. Spatially homogeneous solutions, that are solutions of the form $(v_j(t),a_j(t))=(v(t),a(t))$ for all $j$, are now solutions of the following system of differential equations
\bqs
\left\{
\begin{split}
v'(t)&=F_p(v(t))-\gamma a(t),\\
 a_j'(t)&= \frac{ v(t)-a(t)}{\tau},
\end{split}
\right.
\eqs
whose stationary solutions $(v(t),a(t))=(v,a)\in\R^2$ are given by $F_p(v)=\gamma v$ and $v=a$. For small values of $\gamma$ with other values of the parameters such that $\Lambda\in\mathcal{P}$, we still expect to have three stationary homogeneous solutions with a down-state and an up-state which are stable for the dynamics of the above spatially homogeneous system. On the other hand, upon increasing $\gamma$ passed a given critical value, we expect to have a unique stationary homogeneous solution, corresponding to a stable down state. As a consequence, a working regime for the fully spatially extended system should be $\tau^{-1}\ll1$ and $\gamma\gg1$. This leaves promising directions for follow-up studies.

\section*{Acknowledgements} 

A.A. is supported by ERC Starting Grant OSCI-PRED 101075930. G.F. acknowledges support from the ANR via the projects: Indyana under grant agreement ANR-21-CE40-0008,  ChaMaNe under grant agreement ANR-19-CE40-0024 and from Labex CIMI under grant agreement ANR-11-LABX-0040. R.V. is supported by ERC Advanced Grant GLOW 101096017.  Our work has benefitted from the AI Interdisciplinary Institute ANITI, funded by the France 2030 program under the Grant agreement ANR-23-IACL-0002. Views and opinions expressed are however those of the author(s) only and do not necessarily reflect those of the European Union or the European Research Council. Neither the European Union nor the granting authority can be held responsible for them.

\appendix

\section{Comparison principles}

We state two useful lemmas.

\begin{lem}\label{lemCPZ}
We consider $t\mapsto\boldsymbol{v}(t)\in\mathscr{C}^{1}(\R_+,\ell^\infty(\Z))$ and $t\mapsto\boldsymbol{w}(t)\in\mathscr{C}^{1}(\R_+,\ell^\infty(\Z))$ which satisfy
\bqs
\forall t>0, \quad \forall j\in\Z, \quad v_j'(t)\leq \mathscr{N}\left(v_{j-1}(t),v_j(t),v_{j+1}(t)\right), \text{ and } w_j'(t)\geq \mathscr{N}\left(w_{j-1}(t),w_j(t),w_{j+1}(t)\right),
\eqs
together with a uniform in time bound $\|\boldsymbol{v}(t)\|_{\ell^\infty(\Z)}+\|\boldsymbol{w}(t)\|_{\ell^\infty(\Z)} \leq M$ for all $t\geq0$ and some constant $M>0$. If $v_j(0)\leq w_j(0)$ for all $j\in\Z$, then $v_j(t)\leq w_j(t)$ for all $t>0$ and $j\in\Z$.
\end{lem}

\begin{Proof}
We first introduce the following quantities:
\bqs
0<s_1:=\sup_{|u|\leq M}S'(u)<+\infty, \quad \text{ and } \quad 0<s_2:=\sup_{|u|\leq M}\left|S''(u)\right|<+\infty,
\eqs
together with
\bqs
0<L:=p(s_2(M+1)+s_1^2).
\eqs
Following \cite{CGW08}, for each $0<\epsilon<\frac{1}{3}$, we let
\bqs
T_\epsilon:=\sup\left\{\tau \geq0 ~|~ v_j(t) \leq w_j(t)+\epsilon e^{(L+2s_1)t}(1+\epsilon j^2), \quad \forall t\in[0,\tau], \quad j \in\Z\right\}.
\eqs
Our aim is to prove that $T_\epsilon=+\infty$. By contradiction, assume that $T_\epsilon<+\infty$. By boundedness and continuity of $\boldsymbol{v}$ and $\boldsymbol{w}$, there exists a finite integer $\ell\in\Z$ such that
\bqs
w_\ell(T_\epsilon)+\epsilon e^{(L+2s_1)T_\epsilon}(1+\epsilon \ell^2)-v_\ell(T_\epsilon) = 0 \leq w_j(t)+\epsilon e^{(L+2s_1)t}(1+\epsilon j^2)-v_j(t),
\eqs
for all $t\in[0,T_\epsilon]$ and $j\in\Z$. Let us also note that since $v_j(0)\leq w_j(0)$ for all $j\in\Z$, we necessarily have $T_\epsilon>0$. Next, we compute
\bqs
w'_\ell(T_\epsilon)-v'_\ell(T_\epsilon)= \underset{t\rightarrow T_\epsilon^-}{\lim}~ \frac{w_\ell(T_\epsilon)-v_\ell(T_\epsilon)-(w_\ell(t)-v_\ell(t))}{T_\epsilon-t} \leq - \epsilon (L+2s_1) e^{(L+2s_1)T_\epsilon}(1+\epsilon \ell^2).
\eqs
As a consequence, we have
\begin{align*}
 0&\leq w'_\ell(T_\epsilon)-v'_\ell(T_\epsilon) +\mathscr{N}\left(v_{\ell-1}(T_\epsilon),v_\ell(T_\epsilon),v_{\ell+1}(T_\epsilon)\right)-\mathscr{N}\left(w_{\ell-1}(T_\epsilon),w_\ell(T_\epsilon),w_{\ell+1}(T_\epsilon)\right) \\
 & \leq - \epsilon (L+2s_1) e^{(L+2s_1)T_\epsilon}(1+\epsilon \ell^2) - (1-p)(v_\ell(T_\epsilon)-w_\ell(T_\epsilon)) \\
 &+ (1-p-q)(S(v_{\ell-1}(T_\epsilon))-S(w_{\ell-1}(T_\epsilon))+q (S(v_{\ell+1}(T_\epsilon))-S(w_{\ell+1}(T_\epsilon))\\
 &~~~ -p\left(S'(v_\ell(T_\epsilon))S(v_\ell(T_\epsilon))-S'(w_\ell(T_\epsilon))S(w_\ell(T_\epsilon))\right)\\
 &~~~ +p S'(v_\ell(T_\epsilon))(v_{\ell-1}(T_\epsilon)-w_{\ell-1}(T_\epsilon))+p(S'(v_\ell(t))-S'(w_\ell(T_\epsilon))w_{\ell-1}(T_\epsilon).
\end{align*}
We now bound each term as follows:
\begin{align*}
(1-p-q)(S(v_{\ell-1}(T_\epsilon))-S(w_{\ell-1}(T_\epsilon)) & \leq (1-p-q)s_1 \epsilon e^{(L+2s_1)T_\epsilon}(1+\epsilon (\ell-1)^2),\\
q (S(v_{\ell+1}(T_\epsilon))-S(w_{\ell+1}(T_\epsilon)) & \leq qs_1 \epsilon e^{(L+2s_1)T_\epsilon}(1+\epsilon (\ell+1)^2),\\
p S'(v_\ell(T_\epsilon))(v_{\ell-1}(T_\epsilon)-w_{\ell-1}(T_\epsilon)) & \leq p s_1 \epsilon e^{(L+2s_1)T_\epsilon}(1+\epsilon (\ell-1)^2),
\end{align*}
where in the last inequality, we used the fact that $0<S'(u)<s_1$ for all $|u|\leq M$. For the other terms, we have
\begin{align*}
-p\left(S'(v_\ell(T_\epsilon))S(v_\ell(T_\epsilon))-S'(w_\ell(T_\epsilon))S(w_\ell(T_\epsilon))\right) & \leq p(s_2+s_1^2)\epsilon e^{(L+2s_1)T_\epsilon}(1+\epsilon \ell^2),\\
p(S'(v_\ell(t))-S'(w_\ell(T_\epsilon))w_{\ell-1}(T_\epsilon)& \leq p s_2M \epsilon e^{(L+2s_1)T_\epsilon}(1+\epsilon \ell^2).
\end{align*}
As a consequence, we have obtained that
\begin{align*}
 0&\leq w'_\ell(T_\epsilon)-v'_\ell(T_\epsilon) +\mathscr{N}\left(v_{\ell-1}(T_\epsilon),v_\ell(T_\epsilon),v_{\ell+1}(T_\epsilon)\right)-\mathscr{N}\left(w_{\ell-1}(T_\epsilon),w_\ell(T_\epsilon),w_{\ell+1}(T_\epsilon)\right) \\
 &\leq \epsilon s_1 e^{(L+2s_1)T_\epsilon}\left( -2(1+\epsilon \ell^2)+(1-q)(1+\epsilon (\ell-1)^2)+q(1+\epsilon (\ell+1)^2) \right) \\
 &=\epsilon s_1 e^{(L+2s_1)T_\epsilon}\left( -1 + \epsilon(-\ell^2+\ell(-2+4q)+1)  \right)\leq \epsilon s_1 e^{(L+2s_1)T_\epsilon}\left( -1 + 3\epsilon\right)<0,
 \end{align*}
 which is a contradiction. Thus, for all $t\geq0$, we have obtained 
 \bqs
 v_j(t) \leq w_j(t)+\epsilon e^{(L+2s_1)t}(1+\epsilon j^2), \quad \forall j \in\Z,
 \eqs
 then sending $\epsilon\rightarrow0$ we have reached the conclusion.
\end{Proof}

\begin{lem}\label{lemCPN}
We consider $t\mapsto\boldsymbol{v}(t)\in\mathscr{C}^{1}(\R_+,\ell^\infty(\N))$ and $t\mapsto\boldsymbol{w}(t)\in\mathscr{C}^{1}(\R_+,\ell^\infty(\N))$ which satisfy
\bqs
\forall t>0, \quad  \forall j\geq 1, \quad v_j'(t)\leq \mathscr{N}\left(v_{j-1}(t),v_j(t),v_{j+1}(t)\right), \text{ and } w_j'(t)\geq \mathscr{N}\left(w_{j-1}(t),w_j(t),w_{j+1}(t)\right),
\eqs
and $v_0(t)\leq w_0(t)$ for all $t\geq0$, together with a uniform in time bound $\|\boldsymbol{v}(t)\|_{\ell^\infty(\N)}+\|\boldsymbol{w}(t)\|_{\ell^\infty(\N)} \leq M$ for all $t\geq0$ and some constant $M>0$. If $v_j(0)\leq w_j(0)$ for all $j\geq1$, then $v_j(t)\leq w_j(t)$ for all $t>0$ and $j\geq1$.
\end{lem}

\begin{Proof}
The proof is a straightforward adaptation of the proof of the previous lemma with $T_\epsilon$ now defined as
\bqs
T_\epsilon:=\sup\left\{\tau \geq0 ~|~ v_j(t) \leq w_j(t)+\epsilon e^{(L+2s_1)t}(1+\epsilon j^2), \quad \forall t\in[0,\tau], \quad j \geq1\right\}.
\eqs
for $0<\epsilon<\frac{1}{3}$.
\end{Proof}

\end{document}